\documentclass[reqno,11pt,titlepage]{amsart}

\usepackage{fancyhdr}
\usepackage{wrapfig} 
\usepackage{graphicx}
\usepackage{amssymb}
\usepackage{amsmath}
\usepackage{amsfonts}
\usepackage{amsthm}
\usepackage{color}
\usepackage[colorlinks=true, allcolors=blue]{hyperref}
\usepackage[mathscr]{eucal}
\usepackage{enumerate}
\usepackage{times}
\usepackage{algorithm}
\usepackage{algorithmic}
\usepackage[skip=10pt]{caption} 

\usepackage{xpatch}

\usepackage{lipsum} 

\makeatletter
\xpatchcmd{\@setaddresses}
  {\addvspace\bigskipamount}
  {\penalty\z@\vskip\bigskipamount}
  {}{}
\makeatother
\usepackage{booktabs}
\usepackage{changepage}
\usepackage{pifont}

\usepackage{longtable} 

\usepackage[a4paper,margin=1.5cm]{geometry} 

\allowdisplaybreaks[1]

\newtheorem{Def}{Definition}

\newcommand{\beq}{\begin{equation}}
\newcommand{\eeq}{\end{equation}}
\newcommand{\Proof}{\begin{proof}}
\newcommand{\QED}{\end{proof} \noindent}

\newtheorem{Conj}{Conjecture}

\newcommand{\Z}{\mathbb{Z}}

\title[CayleyPy Growth]{CayleyPy Growth:\\ Efficient growth computations\\ and hundreds of new conjectures\\ on Cayley graphs\\(Brief version)}

\author[A.~Chervov]{A. Chervov}
\address{Institut Curie, CNRS UMR168, Paris, France}
\email[A.~Chervov]{al.chervov@gmail.com}

\author[D.~Fedoriaka]{D. Fedoriaka}
\address{University of Washington}
\email[D.~Fedoriaka]{fedimser@cs.washington.edu}

\author[E.~V.~Konstantinova]{E.V. Konstantinova}
\address{Three Gorges Mathematical Research Center, China Three Gorges University, \\ 8 University Avenue, Yichang 443002, Hubei Province, P.R. China\\
Sobolev Institute of Mathematics, Ak. Koptyug av. 4, Novosibirsk, 630090, Russia\\
Novosibirsk State University, Pirogova str. 2, Novosibirsk, 630090, Russia}
\email[Elena~V.~Konstantinova]{e\_konsta@ctgu.edu.cn, e\_konsta@math.nsc.ru}

\author[A.~Naumov]{A. Naumov}
\address{Independent Researcher}
\email[A.~Naumov]{anton.yr.naumov@gmail.com}

\author[I.~Kiselev]{I. Kiselev}
\address{Accenture}
\email[I.~Kiselev]{igor.kiselev@gmail.com}

\author[A.~Sheveleva]{A. Sheveleva}
\address{Moscow State University}
\email[A.~Sheveleva]{anastasiia.sheveleva@math.msu.ru, art2005art.463@gmail.com}

\author[I.~Koltsov]{I. Koltsov}
\address{Independent Researcher}
\email[I.~Koltsov]{ivankolt@gmail.com}

\author[S.~Lytkin]{S. Lytkin}
\address{Kazakh-British Technical University}
\email[S.~Lytkin]{smlytkin@gmail.com}

\author[A.~Smolensky]{A. Smolensky}
\address{Neapolis University Pafos, Cyprus}
\email[A.~Smolensky]{andrei.smolensky@gmail.com}

\author[A. Soibelman]{A. Soibelman}
\address{IHES}
\email[A. Soibelman]{asoibelman@gmail.com}

\author[F.~Levkovich-Maslyuk]{F. Levkovich-Maslyuk}
\address{City St George's, University of London}
\email[F.~Levkovich-Maslyuk]{fedor.levkovich@gmail.com}

\author[R.~Grimov]{R. Grimov}
\address{Independent Researcher}
\email[R.~Grimov]{grimovr@gmail.com}

\author[D.~Volovich]{D. Volovich}
\address{Hebrew University of Jerusalem}
\email[D.~Volovich]{d.e.volovich@gmail.com}

\author[A.~Isakov]{A. Isakov}
\address{ITMO University}
\email[A.~Isakov]{aoisakov@itmo.ru}

\author[A.~Kostin]{A. Kostin}
\address{Moscow Center for Advanced Studies}
\email[A.~Kostin]{anton.kostin@gmail.com}

\author[M.~Litvinov]{M. Litvinov}
\address{National Research University Higher School of Economics}
\email[M.~Litvinov]{litvinovmitch11@gmail.com}

\author[N.~Vilkin-Krom]{N. Vilkin-Krom}
\address{Ariel University}
\email[N.~Vilkin-Krom]{askuaov@gmail.com}

\author[A.~Bidzhiev]{A. Bidzhiev}
\address{Saint-Petersburg Pasteur Institute}
\email[A.~Bidzhiev]{alimbj09@gmail.com}

\author[A.~Krasnyi]{A. Krasnyi}
\address{Criteo AI Lab, Berlin, Germany}
\email[A.~Krasnyi]{artkrasnyy@gmail.com}

\author[M.~Evseev]{M. Evseev}
\address{University of Bologna}
\email[M.~Evseev]{mixnota@gmail.com}

\author[E.~Geraseva]{E. Geraseva}
\address{Faculty of Bioengineering and Bioinformatics, Lomonosov Moscow State University, 119234, Moscow, Russia}
\email[E.~Geraseva]{geraseva@fbb.msu.ru}

\author[L.~Grunwald]{L. Grunwald}
\address{Sobolev Institute of Mathematics}
\email[L.~Grunwald]{mathmanlily@gmail.com}

\author[S.~Galkin]{S. Galkin}
\address{PUC-Rio, Departamento de Matem\'{a}tica, Rua Marqu\^{e}s de S\~{a}o Vicente 225, G\'{a}vea, Rio de Janeiro, Brazil}
\email[S.~Galkin]{sergey@puc-rio.br}
\thanks{S.~Galkin is supported by CNPq grants PQ 315747 and PQ 308303, and Coordena\c{c}\~{a}o de Aperfei\c{c}oamento de Pessoal de N\'{i}vel Superior-Brasil (CAPES)-Finance Code 001.}

\author[E. Koldunov]{E. Koldunov}
\address{Researcher at Banner Stat}
\email{koldunov.eduard1@gmail.com}

\author[S.~Diner]{S. Diner}
\address{Independent Researcher}
\email[S.~Diner]{stanislav.diner@gmail.com}

\author[A.~Chevychelov]{A. Chevychelov}
\address{Independent Researcher}
\email[A.~Chevychelov]{heavy4evy@gmail.com}

\author[E.~Kudasheva]{E. Kudasheva}
\address{Independent Researcher}
\email[E.~Kudasheva]{evelinn\_a@mail.ru}

\author[A.~Sychev]{A. Sychev}
\address{MIREA — Russian Technological University}
\email[A.~Sychev]{ sychev.a.e@edu.mirea.ru, arseniysychev2015@gmail.com}

\author[A.~Kravchenko]{A. Kravchenko}
\address{Independent Researcher}
\email[A.~Kravchenko]{outsidenessx@gmail.com}

\author[Z.~Kogan]{Z. Kogan}
\address{Independent Researcher}
\email[Z.~Kogan]{zahar1991@gmail.com}

\author[A.~Natyrova]{A. Natyrova}
\address{Independent Researcher}
\email[A.~Natyrova]{natyrovaaltana@gmail.com}

\author[L.~Shishina]{L. Shishina}
\address{Independent Researcher}
\email[L.~Shishina]{ldsisina@gmail.com}

\author[L.~Cheldieva]{L. Cheldieva}
\address{Independent Researcher}
\email[L.~Cheldieva]{liuda.tarusina@gmail.com}

\author[V.~Zamkovoy]{V. Zamkovoy}
\address{Independent Researcher}
\email[V.~Zamkovoy]{zamkovoyvladislav@gmail.com}

\author[D.~Kovalenko]{D. Kovalenko}
\address{Independent Researcher}
\email[D.~Kovalenko]{caymon14@gmail.com}

\author[O.~Papulov]{O. Papulov}
\address{Independent Researcher}
\email[O.~Papulov]{ob.papulov@gmail.com}

\author[S.~Kudashev]{S. Kudashev}
\address{Independent Researcher}
\email[S.~Kudashev]{sergey0474@gmail.com}

\author[D.~Shiltsov]{D. Shiltsov}
\address{Independent Researcher}
\email[D.~Shiltsov]{da.shiltsov@gmail.com}

\author[R.~Turtayev]{R. Turtayev}
\address{Independent Researcher}
\email[R.~Turtayev]{u1am.n1t7@gmail.com}

\author[O.~Nikitina]{O. Nikitina}
\address{Independent Researcher}
\email[O.~Nikitina]{ol.ya.nik.dev@gmail.com}

\author[D.~Mamayeva]{D.Mamayeva}
\address{Independent Researcher}
\email[D.~Mamayeva]{dmamayeva367@gmail.com}

\author[S.~A.~Nikolenko]{S. A. Nikolenko}
\address{Independent Researcher}
\email[S.~A.~Nikolenko]{Nikolenko.Sergei@icloud.com}

\author[M.~Obozov]{M. Obozov}
\address{Research Center of the Artificial Intelligence Institute, Innopolis University}
\email[M.~Obozov]{obozovmark9@gmail.com}

\author[A.~Titarenko]{A. Titarenko}
\address{Independent Researcher}
\email[A.~Titarenko]{titarenkoav.84@mail.ru}

\author[A.~Dolgorukova]{A. Dolgorukova}
\address{Independent Researcher}
\email[A.~Dolgorukova]{an.dolgorukova@gmail.com}

\author[A.~Aparnev]{A. Aparnev}
\address{Inpdependent Researcher}
\email[A.~Aparnev]{apich238@gmail.com}

\author[O.~Debeaupuis]{O. Debeaupuis}
\address{Institut Curie, CNRS UMR168\\
Imagine Institute, INSERM UMR 1163, Paris, France}
\email[O.~Debeaupuis]{orianne.debeaupuis@curie.fr, orianne.debeaupuis@institutimagine.org, orianne.debeaupuis@gmail.com}

\author[Simo~Alami~C.]{Simo Alami C.}
\address{LIX, Ecole Polytechnique/CNRS, IP Paris}
\email[Simo~Alami~C.]{mohamed.alami-chehboune@polytechnique.edu}

\author[H.~Isambert]{H. Isambert}
\address{Institut Curie, CNRS UMR168, Paris, France}
\email[H.~Isambert]{Herve.Isambert@curie.fr}

\begin{document} 

\begin{abstract}
This is the third paper of the CayleyPy project applying artificial intelligence methods to problems in group theory. 
We announce the first public release of CayleyPy, an open-source Python library for computations with Cayley and Schreier (coset) graphs. Compared with state-of-the-art systems based on classical methods, such as GAP and Sage, CayleyPy handles significantly larger graphs and performs several orders of magnitude faster.

Using CayleyPy we obtained about 200 new mathematical conjectures on Cayley and Schreier graphs, with special regard to their diameters and growths, 
which we present in this paper.

For many Cayley graphs of symmetric groups $S_n$ we observe quasi-polynomial diameter formulas: a small set of quadratic or linear polynomials indexed by $n \mod s$, and conjecture that it is a general phenomenon.  These lead to efficient diameter computation, despite the problem being NP-hard in general.  We propose  refinement of the  Babai-type conjecture on diameters of $S_n$: 
$\frac12 n^2 + 4n$ upper bounds for the diameters in the standard undirected case,
as compared to prior conjectural bounds of $O(n^2)$.
We also provide explicit generator families, related to involutions in a simple “square-with-whiskers” pattern, 
which we conjecture to maximize the diameter; extensive (and in some cases exhaustive) search confirms this for all $n \leq 15$. 
We conjecture an answer to the celebrated open question raised by the "founding father of Soviet cybernetics" V.\,M.~Glushkov in~1968: the diameter of the directed Cayley graph generated by the left cyclic shift and the transposition of the first two elements is equal to $(3n^2-8n+9)/4$ for $n$ odd, and to $(3n^2-8n+12)/4$ for $n$ even.

For nilpotent groups we conjecture an improvement of J.\,S.~Ellenberg's results on the diameters of the upper unitriangular matrices over $\Z/p\Z$, presenting a phenomenon of linear dependence of the diameter on $p$. Moreover, the growth for nilpotent groups is conjectured to closely follow Gaussian distributions, that is, to exhibit a central limit phenomenon similar to the results of P.~Diaconis for $S_n$. 

Some of our conjectures are "LLM-friendly" --- they can be stated as sorting problems, which are easy to formulate for LLM, and their solutions can be given by an algorithm or by a Python code, which is easy to verify, so they can be used to test LLM's abilities to solve research problems. 
To benchmark various methods of path-finding on Cayley graphs we create more than 10 benchmark datasets in the form of Kaggle challenges, making benchmarking easy and public to the community. 
CayleyPy works with arbitrary permutation or matrix groups, and supports a pre-defined collection of more than a hundred generators
including puzzle groups.
Our code for direct growth computation outperforms similar functions on the standard computer algebra system GAP/SAGE up to 1000 times both in speed and in maximum sizes of the graphs that it can handle. 

\bigskip
\begin{center}\large
Project page:
\quad 
\url{https://github.com/CayleyPy/CayleyPy}
\quad
\includegraphics[width=.04\textwidth]{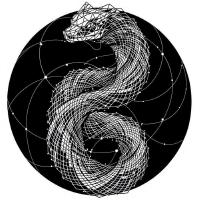}
\end{center}
\end{abstract}


\maketitle 

\setcounter{tocdepth}{2}
\small
\tableofcontents
\normalsize

\section{Introduction}
{\bf Brief.}
Deep learning methods have shown tremendous success in various fields. The goal of the CayleyPy project is to apply deep learning methods to group and graph theories providing ability to work with googol size graphs and groups, with a current focus on Cayley graphs of finite groups. 
In the first paper of the project, \cite{CayleyPyCube}, 
we proposed a new deep learning method for decomposing a given group element into a product of generators (or, equivalently, finding a path on the Cayley graph), which significantly outperforms all known analogs---in particular, the randomized Schreier--Sims algorithm implemented in the standard computer algebra system GAP/SAGE. In the second paper, \cite{CayleyPyRL}, 
we provided the first demonstration of how these methods can be used to advance concrete open mathematical conjectures such as  \href{https://oeis.org/A186783}{OEIS-A186783}. In a nutshell, the pipeline is as follows: computational experiments + observe pattern + try to prove it + iterate. 

The goal of the present paper is to demonstrate that this pipeline scales quite well for Cayley graph tasks. We analyze nearly half a hundred various Cayley/Schreier graphs, producing nearly 200 new mathematical conjectures, several new constructions and theorems. 

The other goal of the paper is to briefly present the first release of ``CayleyPy'', an AI-based Python open-source library, which provides a simple interface for mathematicians to make computational experiments with Cayley graphs more efficiently than other analogs, in particular, the standard computer algebra systems GAP/SAGE. The current release supports working with arbitrary generators provided by the user, and also includes built-in support for dozens of known Cayley graphs of $S_n, A_n$ and matrix groups. The focus of the release is direct growth computations by efficient versions of breadth first search algorithm,
which outperforms GAP up to 1000 times both in speed and achievable sizes of the groups. Visualizations of Cayley graphs, random walks on them, and beam-search algorithms are also supported. The current version also provides limited support for the deep learning component, which will be polished in the next versions. The project links:

\begin{itemize}
    \item 
    \href{https://github.com/cayleypy/cayleypy}{Project}: \verb|https://github.com/cayleypy/cayleypy|;
    \item 
    \href{https://cayleypy.github.io/cayleypy-docs/api.html}{Documentation}: \verb|https://cayleypy.github.io/cayleypy-docs/api.html|.
\end{itemize}

\subsection{Context. Cayley graphs as a testbed for RL methods}
The narrow context is the study of Cayley graphs which are of fundamental importance for mathematics.
The broader context: graph path-finding represents a much wider class of ``planning/reinforcement learning'' problems
important in various fields. In a nutshell one needs to perform a sequence of elementary actions such that the final output achieves a
certain goal. In the group theory setting these elementary actions are generators of the group (edges of the Cayley graph),
in games like Go or Chess these are the ``game moves'', 
in robotics these are the elementary moves of a robot manipulator, 
in LLM these are tokens generated one by one,
in automatic theorem proving in mathematics these are arguments generated one by one, which should lead to a correct proof.
All these problems are approached by the same technique of reinforcement learning and its variations.
The setup of reinforcement learning is practically identical to the graph setup; one has the following dictionary:
states corresponds to nodes of the graph, actions correspond to edges, rewards correspond to weights of edges, and thus the goal of optimizing cumulative reward
corresponds to minimizing the path length, i.e., optimal graph path-finding. 

From that perspective, Cayley graphs represent an ideal testbed for the more general problem mentioned above, because it is extremely easy to create Cayley graphs, just by picking sets of matrices or permutations.
There are many open research problems which successful approaches can help to resolve and thus write names
into history of science. One can easily choose between tasks with different levels of complexity, from elementary to very difficult.
The actions (moves) are easy to describe as matrix multiplication or as permutation applications, which work very fast on GPU and CPU, thus the training data is practically unlimited. For example, solving Rubik's cube or other puzzles is a prototypical example
of the Cayley graph path-finding --- developing optimal solvers for $4\times4\times4$ Rubik's cube and higher, as well as estimating their diameters
(``God's numbers''), are classical hard problems unresolved for almost half a century, representing other difficult problems on Cayley graphs. 

Moreover, the present time is characterized by growing interest in and number of applications of deep learning methods
to mathematics: machine learning has been emerging as ``a tool in theoretical science''~\cite{douglas2022machine}.
In recent years, this has led to several noteworthy applications to mathematical problems:~\cite{lample2019deep,davies2021advancing, bao2021polytopes, 
romera2024mathematical,
coates2024machine,alfarano2024global, charton2024patternboost,shehper2024makes,swirszcz2025advancing,hashemi2025transformers,he2024ai}.
Seewoo Lee created a repository that collects papers in AI for mathematics,  \href{https://seewoo5.github.io/awesome-ai-for-math/}{
Awesome AI for Math}.

\subsection{Problem statements on Cayley graphs from the mathematical perspective}
Cayley graphs are fundamental in group theory~\cite{gromov1993geometric},\cite{tao2015expansion}, and have various applications: bioinformatics
~\cite{Pevzner1995human2mice, Pevzner1999cabbage2turnip,  bulteau2019parameterized}; processor interconnection networks~\cite{akers1989group, cooperman1991applications,heydemann1997cayley}; coding theory and cryptography~\cite{hoory2006expander,zemor1994hash,petit2013rubik}; quantum computing~\cite{ruiz2024quantum,sarkar2024quantum,dinur2023good, acevedo2006exploring,gromada2022some}, etc.

There are many open conjectures in the subject and making progress in their understanding is a fundamental challenge in the field. Two of these that are quite well-known, easy to formulate, wide open and most relevant to us are: 
\begin{itemize}
    \item 
    {\bf Babai-like conjecture:} for any choices of generators the diameter of $S_n$ is $O(n^{2})$ (see, e.g., ~\cite{helfgott2014diameter}, \cite{helfgott2019growth}, \cite{helfgott2015random}); 

    \item 
    {\bf Diaconis conjecture \cite{diaconis2013some}:}  the mixing time for random walks is $O(n^{3}\log n)$ (again for any choices of generators).
\end{itemize}

An important characteristic of a Cayley graph $G$ is its \textit{growth}
--- the vector of sizes of \textit{spheres} (or \textit{layers}) containing all elements at the same distance
(length of the shortest path) from some fixed element $g_0 \in G$.
It is easy to see that for Cayley graphs (in strict sense) growth does not depend on the choice of $g_0$.
The growth of a graph $G$ contains a lot of information on it: for example, the diameter is just the length of the growth vector minus one.
It is suggestive to view growth as an unnormalized probability distribution over $\mathbb N$, as in P.~Diaconis' works. 
Then the diameter is just the maximum of a random variable. And it is important to understand its other characteristics:
its mean, mode, moments, etc. Ideally, the goal is to understand from what family of distributions it comes,
and, hopefully, to observe some universality phenomena, such as the distribution approaching something known for large values of $n$.
For example, the Gaussian normal distribution trivially arises in the case of abelian groups, 
while in the case of $S_n$ with generators close to commutative --- like Coxeter's neighbor transpositions $(i,i+1)$ --- the appearance of the Gaussian normal approximation has been demonstrated by P.~Diaconis (in some sense
based on classical works in statistics by Kendall, Mann, Whitney, Wilcoxon). 

So, having some elements in, say, the permutation group $S_n$ (or in some other group), one constructs
the Cayley graph, and there is a set of natural questions and lines of investigation: 
\begin{itemize}
  \item What group is obtained?
  \item What is its diameter?
  \item Growth statistical characteristics: mean, mode, moments, what distribution does it follow (or at least asymptotically as $n\to \infty$)?
  \item Algorithm: is there an effective/polynomial algorithm which decomposes a given element into a product of generators (optimally/sub-optimally)? 
  \item Antipodes (``super-flips''): is there an explicit description of the longest elements? 
  \item Can one explicitly describe the word-metric (i.e., number of generators in the decomposition of an element --- the length of the shortest path on the Cayley graph)? 
  \item What can be said about the graph's spectrum?
  \item What is the mixing time? 
\end{itemize}

For the vast majority of Cayley graphs, many fundamental questions remain open and
they represent important research problems in the field. 

On the one hand, there are theoretical barriers that limit the prospect for complete solutions of the problems outlined above.
Finding the shortest paths on generic finite Cayley graphs is an NP-hard problem~\cite{even1981minimum} (even P-space complete \cite{jerrum1985complexity}). And it is NP-complete for many specific group families, such as $N\times N \times N$ Rubik's Cube groups~\cite{demaine2017solving} and others~\cite{bulteau2015pancake}. Determining the diameters of generic groups is NP-hard (again~\cite{even1981minimum}).

On the other hand, there are positive results for many generators and it is a huge and active field of research to study questions similar to the above. For example, Coxeter's generators $(i,i+1)$ represent an extreme case where almost all of the above questions have well-known and beautiful answers, see e.g.  the \href{https://en.wikipedia.org/wiki/Bubble_sort}{bubble sort} algorithm which provides an optimal decomposition.
It is not the only example and  progress can be achieved in many other cases, as we aim to show.

\subsection{Main Contributions}
The aim of the present paper is to make progress in the direction of fundamental problems described above
(i.e. understanding of various properties of Cayley graphs)
with the help of the new tool which we are developing: the AI-based Python open-source library CayleyPy 
which allows one to make computational experiments orders of magnitude more effectively than
standard computer algebra systems GAP/SAGE. 

\begin{itemize}
  \item We present the first release of the CayleyPy library with documentation and tutorial notebooks,
  which provides the ability for non-professionals in software engineering
  to use the power of the effective implementations and GPU computations
  for several tasks for Cayley graphs, mainly growth computations and elements decompositions (path-finding). For example, our direct growth computation outperforms similar functions in the standard computer algebra system GAP/SAGE up to 1000 times in speed and achievable sizes of the graphs.
  \item We generate around 200 conjectures on various properties of Cayley graphs,
  which are backed up by extensive computational experiments on nearly 50 Cayley graphs.
  The conjectures are summarized in Tables~\ref{tab:long1} and~\ref{tab:long2}.
  All these generators are included in CayleyPy as named generators, and the results of hard computations are available on GitHub.
  More than 250 notebooks with computational experiments are \href{https://www.kaggle.com/datasets/alexandervc/cayleypy-development-3-growth-computations/code}{publicly available} on the Kaggle platform, where it is easy to reproduce them. 
  \item In particular, we propose the following:
  \begin{itemize}
      \item We conjecture that the diameters of many $S_n$-Cayley graphs are quasi-polynomials (quadratic/linear) in $n$ (i.e. several polynomials depending on $n$ modulo some $s$), allowing one to find them rather efficiently, which is surprising since it is NP-hard in general. Additionally,  distances from identity to constructive elements 
      (i.e. their "\href{https://en.wikipedia.org/wiki/Word_metric}{word metric}") are observed to be quasi-polynomials, 
      and that is conjectured to be a general phenomenon.
      \item The improvement of the L.~Babai-type conjecture for $S_n$. Specifically, the diameters are bounded by $n^2/2 + 4n$ 
      for standard undirected Cayley graphs of $S_n$,
      improving on the earlier conjectural $O(n^2)$ bounds. Moreover, we present explicit families of generators for $S_n$ which conjecturally provide the largest (or nearly largest) diameters. They are related to involutions and follow a rather simple and beautiful pattern (that we call ``square-with-whiskers''). They were found by an extensive (and in some cases exhaustive) search of the generators with maximum diameter for $n\leqslant 15$, the outcomes of which we also present. 
      \item For nilpotent groups we conjecture an improvement of J.\,S.~Ellenberg's results on the diameter of upper unitriangular  matrices over $\Z/p\Z$, presenting the phenomenon of linear dependence of the diameter on $p$. Moreover, the growth for nilpotent groups is conjectured to follow Gaussian distributions (a central limit phenomenon, similar to certain results of P.~Diaconis for $S_n$). 
      \item We present a conjectural answer to the question raised by ``\href{https://en.wikipedia.org/wiki/Victor_Glushkov}{one of the founding fathers of Soviet cybernetics}'' V.\,M.~Glushkov in~1968, which has been much studied but not yet resolved: the diameter of the directed Cayley graph generated by the left cyclic shift and the transposition $(1,2)$ is equal to $(3n^2-8n+9)/4$ for odd $n$, and to $(3n^2-8n+12)/4$ for even $n$.
   \end{itemize}  
   \item Some of our conjectures are ``LLM-friendly'' --- they can be stated as sorting problems and thus are easy to formulate for LLMs. Moreover, solutions can be given by an algorithm or a Python code which is easy to verify, so they can be used to test abilities of LLMs to solve research problems. 
   \item To benchmark various methods of path-finding on Cayley graphs we create more than 10 benchmark datasets in the form of Kaggle challenges, making benchmarking easy and public to community. 
\end{itemize}

\clearpage
\section{Background and related works}
\subsection{Reminder on Cayley graphs, diameters, growth }
\href{https://en.wikipedia.org/wiki/Cayley_graph}{Cayley} (and more general \href{https://en.wikipedia.org/wiki/Schreier_coset_graph}{Schreier coset}) graphs represent a fundamental concept in mathematics.

{\bf Non-technical definition.} A (generalized) Cayley graph is defined by any set of vectors $V=\{v_i\}$ and any set of matrices $\mathcal M=\{M_k\}$ (``moves''), as follows: the nodes of the graph correspond to $v_i$, and there is an edge between $i$ and $j$ if there exists a matrix $M_k$ from our set $\mathcal M$ such that $ v_i = M_k v_j $. 

In the strict mathematical sense, Cayley and Schreier graphs are distinguished by imposing certain conditions on $v_i$ and $M_k$. 
While a detailed distinction is not critical for the scope of this paper, we mention the formal requirements for completeness.  
For both cases: $M_k$ are invertible matrices (this corresponds to the ``group'' part);
the set of vectors is closed under the $M_k$ action, that is, $M_k v_i \in V$ for any $M_k$ and $v_i$ (i.e. "no forbidden moves''); and the graphs are connected. Another requirement for Cayley graphs is the absence of fixed vectors for any $M_k$ in $V$,
which implies that $V$ actually coincides with the set generated by all $M_k$, i.e. with the group generated by $M_k$.

Clearly the concept is rather general and represents a more general idea of state-transition graph:
for any system which has some set of states and transitions between them, these states can be taken as nodes,
with edges connecting those states for which the transition exists.
Several examples are: positions in the Go/Chess game with edges corresponding to moves between positions,
in reinforcement learning the states correspond to graph nodes and the actions correspond to edges,
for Rubik's cubes and other puzzles the nodes correspond to all positions of the puzzle and the
edges correspond to moves, in bioinformatics the gene sequences correspond to nodes while ``mutations''
define edges, etc. 

For any graph $G$ the \textit{distance} between two vertices $u, v \in G$ is the length of the shortest path between $u$ and $v$:
$$
d(u, v) = \min \{n \colon u \leftrightarrow v_1 \leftrightarrow v_2 \leftrightarrow \ldots \leftrightarrow v_{n-1} \leftrightarrow v\} \ .
$$

\begin{Def} \textbf{Diameter} of the graph $G$ is the distance between the two most distant nodes:
$\mathrm{diam}(G) = \max\limits_{u, v \in G} d(u, v)$.
\end{Def}

In the context of puzzles, the diameter of the corresponding state-transition graph is known as "God's number".
It represents the number of moves needed to solve the hardest configuration, i.e. the "worst case performance". If the graph represents a communication network, God's number corresponds to  the largest latency in the network.

Estimating the diameters of various Cayley and other graphs is typically a difficult problem.
It was proven to be NP-hard for general Cayley graphs \cite{even1981minimum},
and is very difficult to solve in particular cases. For example, it took 40 years to determine "God's number''
for the standard Rubik's cube \cite{Rokicki2014Diameter}, and it is still unknown for its variations and higher versions. 
It is the subject of many deep mathematical conjectures like the L.~Babai's conjecture discussed above.

\begin{Def} 
\textbf{Growth} $\gamma_G$ of the graph $G$
is the vector of sizes of nonempty \textbf{layers} containing all elements at the same distance from some fixed element $v_0 \in G$:
$$
\gamma_G(m)= \vert \{v\in G\colon d(v, v_0) = m\}\vert, \quad m \in \mathbb N \cup \{0\}, \quad v_0 \in G.
$$
\end{Def}

Cayley graphs (in strict mathematical sense) are highly symmetric, and growth does not depend on the choice of $v_0$, thus being a natural characteristic of the graph and the group itself.
Growth encapsulates a wealth of important information; in particular, $\mathrm{diam}(G) = \mathrm{len}(\gamma_G)
 - 1$. Making progress in understanding growth for various Cayley graphs is an important research problem. 

\subsection{Related works }
There is an enormous volume of work on Cayley graphs, some notable monographs are~\cite{gromov1993geometric} and~\cite{tao2015expansion}.
There are also various applications: in bioinformatics
~\cite{Pevzner1995human2mice, Pevzner1999cabbage2turnip,  bulteau2019parameterized} for estimation of the evolutionary distances; in processor interconnection networks~\cite{akers1989group, cooperman1991applications,heydemann1997cayley}
(as previously mentioned, there the diameter corresponds to the worst latency and one is interested in choosing the graphs with the least diameter); in coding theory and cryptography~\cite{hoory2006expander,zemor1994hash,petit2013rubik}; in quantum computing~\cite{ruiz2024quantum,sarkar2024quantum,dinur2023good, acevedo2006exploring,gromada2022some}, etc.

{\bf Diameters estimation.} Huge interest in mathematics is attracted to the problem of the diameters estimation. 
An excellent survey \cite{glukhovzubov1999lengths} covers many developments from 60s and up to 2000,
providing more than 20 concrete examples of generators of permutations realizing the diameters (part of the CayleyPy collection of generators is borrowed from there). 
It also covers many Soviet school developments not well known worldwide.
Unfortunately, this survey does not seem to be translated to English. 
The survey \cite{konstantinova2008some} additionally covers some extra topics such as Hamiltonicity and special graphs of interest to  molecular biologists (e.g. sorting by reversals).
In recent years, much work has been done on Babai's conjecture and its variations.
This conjecture predicts in some sense that the diameters are ``small'': 
bound by $\log^c(|G|)$, where $|G|$ is the number of elements in the group, and $c$ is some universal constant.
See, for example, \cite{BGT2011,BGT2012, Breuillard2015, Eberhard2023} and 
H.~Helfgott's surveys ~\cite{helfgott2014diameter,helfgott2019growth,helfgott2015random}. In addition,  
\href{https://www.math.auckland.ac.nz/~conder/SODO-2012/Seress-SODO2012.pdf}{A.~Seress's slides} and S.~Eberhand's informal blog-posts (\href{https://seaneberhard.com/2023/08/15/babais-conjecture-for-generating-sets-containing-transvections/}{2023}, \href{https://seaneberhard.com/2020/11/23/talk-at-tau-about-babais-conjecture-in-high-rank/}{2020Talk}, \href{https://seaneberhard.com/2020/05/21/babais-conjecture-for-at-least-three-random-generators/}{3random}) provide a non-technical introduction. 

{\bf Decomposing elements into products of generators (path-finding on Cayley graphs).}
Another relevant line of research is the algorithmic problem of decomposing an element of the group into a product of generators --- in other words, solving Rubik's cube or other puzzles, or Cayley graph path-finding, or the sorting problem in computer science (all these formulations are equivalent).
Some important milestones:
\begin{itemize}
    \item  The \href{https://en.wikipedia.org/wiki/Schreier%E2%80%93Sims_algorithm}{Schreier–Sims algorithm} \cite{sims1970computational} can in principle work for arbitrary permutation groups.
    Its improved randomized version by Donald Knuth~\cite{knuth1991efficient} is implemented in GAP/SAGE. However, this algorithm is known to be impractical for large groups (e.g., of order $10^40$), as the outputs "are usually exponentially long''~\cite{fiat1989planning}.
    \item It is NP-hard to {\bf optimally} decompose elements for generic finite groups \cite{even1981minimum}, improved to 
    \href{https://en.wikipedia.org/wiki/PSPACE-complete}{P-space complete} in \cite{jerrum1985complexity}
    \item (1998-now) {\bf Optimal} decomposition is NP-complete for many concrete families of generators like Rubik's cube~\cite{demaine2017solving}, \href{https://en.wikipedia.org/wiki/Pancake_sorting}{pancakes}~\cite{bulteau2015pancake}, etc.
    \item Optimal decomposition nevertheless can be achieved for groups of non extra-huge sizes: \cite{korf1997finding}
    proposed the general method of ``pattern databases'' and provided the first demonstration of the possibility to solve the $3\times3\times 3$ Rubik's cube ($4.3\times 10^{19}$ states) optimally. The solver was very slow, but it is currently improved to several cubes per second (see "\href{https://en.wikipedia.org/wiki/God%27s_algorithm}{God's algorithm}''). A big challenge is achieving optimal solution for higher group sizes like $10^{30}$--$10^{40}$; hopefully, it might be resolved by machine learning  which is one of the CayleyPy project goals.
    \item Some examples of algorithms for particular generator families are known. 
    \cite{BafnaPevzner1998reversals} presented a surprising breakthrough by showing that despite usual reversal sorting being NP-complete, for {\bf signed} reversals they found an {\bf optimal polynomial} algorithm. This made possible effective computations of evolutionary distance in biology and stimulated many further developments (survey: \cite{bulteau2019parameterized}). 
    Classical algorithms include bubble sort -- for Coxeter generators (optimal), \href{https://en.wikipedia.org/wiki/Pancake_sorting}{pancake sorting} algorithm (suboptimal) \cite{GatesPapadimitriou1979PrefixReversal}, Rubik's cube solvers (suboptimal), etc.
    \cite{larsen2003sl2} proposed an algorithm for $SL_2(\Z/p\Z)$ with complexity $ O(\log(p) \log(\log p))$, near optimal $O(\log p)$.
    Participants of the Kaggle Challenge \href{https://www.kaggle.com/competitions/santa-2023}{Santa 2023} proposed methods which can effectively solve some puzzles with sizes up to $10^{1000} $ (like Rubik's cube $33\times33\times33$). They used a remarkable idea: first finding "small support'' elements expressed via original generators, then using these new small support generators one can, for example, run beam search just with Hamming distance as guiding heuristics. However the generality of such an approach is unclear --- it is unknown and not even investigated by mathematical community, to the best our knowledge, for what permutation groups those "small support'' elements can be effectively found. It is known that a restriction on the support of even a single generator of the form $\operatorname{supp}\leq 0.63n $ implies a polynomial bound on the diameter \cite{bamberg2014bounds}, see also
     \href{https://www.math.auckland.ac.nz/~conder/SODO-2012/Seress-SODO2012.pdf}{A~Seress's slides}.
\end{itemize}

After the deep learning revolution it became natural to try deep learning methods on this problem.  However, systematic investigations applying deep learning across diverse families of Cayley graphs appear to have been limited prior to our project.
Nevertheless for the specific case of $3\times3\times3$ Rubik's cube there were two notable works which have demonstrated that deep learning methods can effectively solve it: the \href{https://deepcube.igb.uci.edu/}{DeepCube} 
series of papers~\cite{mcaleer2019solving, agostinelli2019solving, khandelwal2024towards,agostinelli2024q}, and later the Efficient Cube:~\cite{takano2023selfsupervision}.  Some others~\cite{brunetto2017deep,johnson2021solving,amrutha2022deep,noever2022puzzle,chasmai2022cubetr,bedaywi2023solving,pan2021fourier} proposed several approaches, but did not achieve a solution. One noteworthy idea~\cite{pan2021fourier} is combining neural networks with the representation theory of the symmetric group ---  a neural net predicts the coefficients of the non-abelian Fourier transform for the distance function. The rationale is to observed sparsity (bandlimitedness) of the Fourier transform of the common distance functions on $S_n$~\cite{swan2017harmonic}.

Recently \cite{douglas2025diffusion} (partially motivated by CayleyPy)  proposed a novel diffussion-like approach to that problem, considered several classes of permutation and matrix groups and also demonstrated that the deep learning approach works successfully not only on Rubik's cube, but on general Cayley graphs up to certain sizes.
Another recent work \cite{ZiarkoBortkiewiczZawalskiEysenbachMilos2025CRTR} provided solution of 333 Rubik's cube based on 
approach of latent representations. 

The recent seminal achievement from S.~Gukov's team~\cite{shehper2024makes} develops an effective AI-based path-finding approach for an infinite group of Andrews--Curtis moves and resolves the Akbulut--Kirby conjecture that remained open for 39 years.

AI methods are also used for path-finding on other types of graphs, such as planning movements in obstacle-rich environments and road networks \cite{Leskovec22pathfinding, Yakovlev23}.

\clearpage
\section{Results Summary}

\subsection{Brief features of CayleyPy }
\subsubsection{Main features }
CayleyPy is an AI-based open-source Python library which can work with googol size graphs,
with the current main focus on mathematical tasks for Cayley graphs of finite groups. 
The key current goals are: using the AI approach (\cite{CayleyPyCube, CayleyPyRL})
to solve path-finding tasks on Cayley graphs, or in other words to decompose a given group element
into a product of generators  (solution of Rubik's cube is an illustrative example).

More generally one is interested in understanding various properties of Cayley graphs:
their diameters, growth, spectrum, random walks mixing time, etc. 
CayleyPy provides a framework to handle all these tasks in a simple manner accessible to non-experts in programming,
but utilizing the power of efficient code, algorithms and GPU accelerators. 
For graphs of not too large size (up to trillions) one can compute 
the diameter and growth effectively, for smaller sizes one can also compute the spectrum and 
generate visualizations, and so on.

CayleyPy can work with arbitrary permutation and matrix groups, which the user defines as an input.
Moreover, it  also supports a large collection (more than a hundred) of predefined generators
of permutation groups (including various puzzles) and matrix groups. 

The main focus of the current release is not the AI-component which will be polished in the next releases, 
but rather presenting the main features for working with the Cayley/Schreier graphs:
defining them by arbitrary permutations or matrices (or using predefined graphs), computations of growth, diameters,
spectrum, visualizations, random walks, etc. 

The main outcome of that stage of the project is that we were able to generate hundreds of new 
mathematical conjectures on Cayley/Schreier  graphs via extensive computational experiments
with CayleyPy. Thus, we demonstrate that effective computational tools can advance discoveries in pure mathematics.

\clearpage
\subsubsection{Examples and Tutorials }

The \href{https://www.kaggle.com/code/fedimser/cayleypy-demo/notebook}{tutorial notebook}
provides a user-friendly introduction on how to work with CayleyPy. 
For the sake of completeness, let us present some examples on using CayleyPy based on it. An example of a simplified version of AI-based path-finding can be found in 
\href{https://www.kaggle.com/code/fedimser/beam-search-with-cayleypy}{notebook}
(it will be reworked in future releases of CayleyPy). 

Let us first list the main resources of CayleyPy:
\begin{itemize}
    \item 
    \href{https://github.com/cayleypy/cayleypy}{Project}:
    {\verb|https://github.com/cayleypy/cayleypy|} 
    \item 
    \href{https://cayleypy.github.io/cayleypy-docs/api.html}{Documentation page}: {\verb|https://cayleypy.github.io/cayleypy-docs/api.html|}
    \item
    \href{https://www.kaggle.com/code/fedimser/cayleypy-demo}{Basic tutorial notebook}: 
    \verb|https://www.kaggle.com/code/fedimser/cayleypy-demo|
    \item 
    More than \href{https://www.kaggle.com/datasets/alexandervc/cayleypy-development-3-growth-computations/code}{250 notebooks} with computational experiments using CayleyPy:\\
    {\verb|https://www.kaggle.com/datasets/alexandervc/|}\\
    {\verb|cayleypy-development-3-growth-computations/code|}
\end{itemize}

\begin{figure}[ht]
  \centering
  \includegraphics[width=0.8\linewidth]{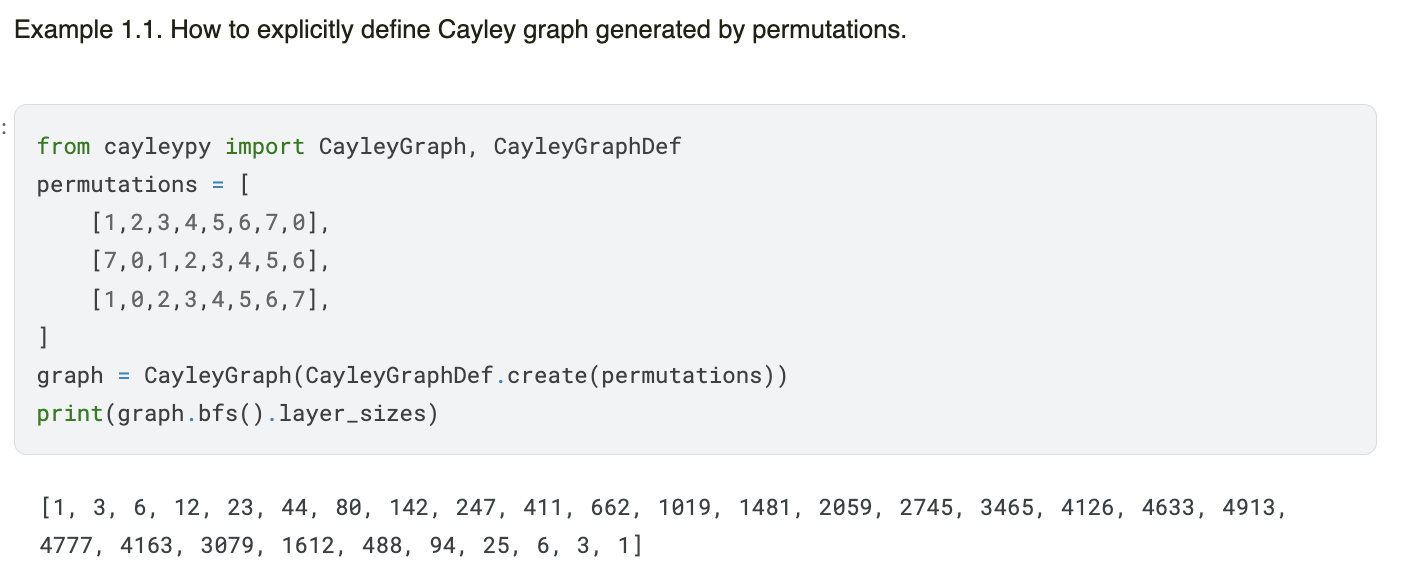}
  \caption{Defining a Cayley graph with given permutations 
  (one-line notation) 
  and computations of the growth.}
  \label{fig:tutotialEx1}
\end{figure}

Figure \ref{fig:tutotialEx1} presents an example on how to define a Cayley graph and how to compute its growth.

\begin{figure}[ht]
  \centering
  \includegraphics[width=0.8\linewidth]{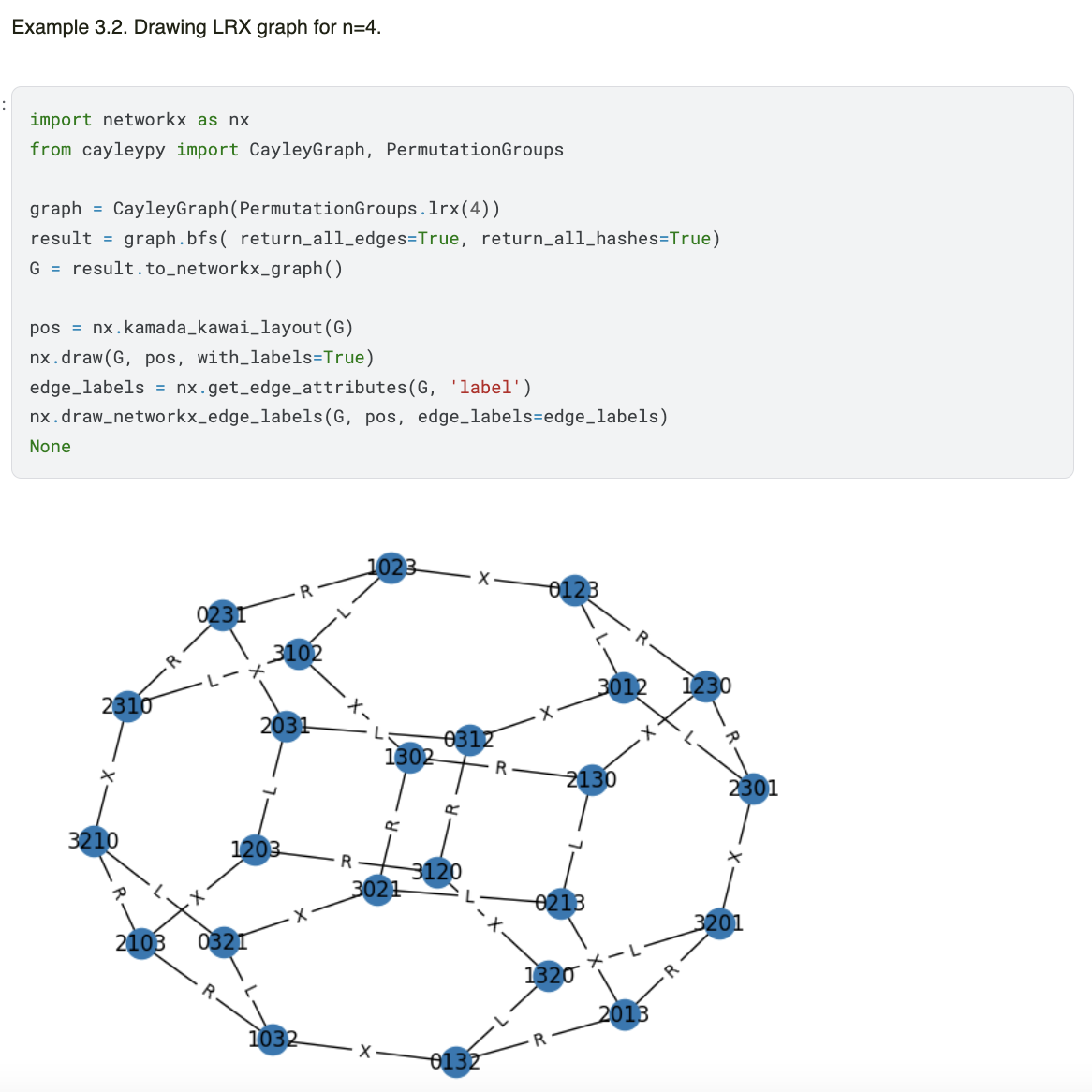}
  \caption{Defining a Cayley graph from the predefined collection ("LRX" generators), computing its adjacency matrix, 
  visualization.}
  \label{fig:tutotialEx2}
\end{figure}

Figure \ref{fig:tutotialEx2} presents an example on how to define a Cayley graph from the predefined collection,
compute its adjacency matrix and visualize it.

\begin{figure}[ht]
  \centering
  \includegraphics[width=0.8\linewidth]{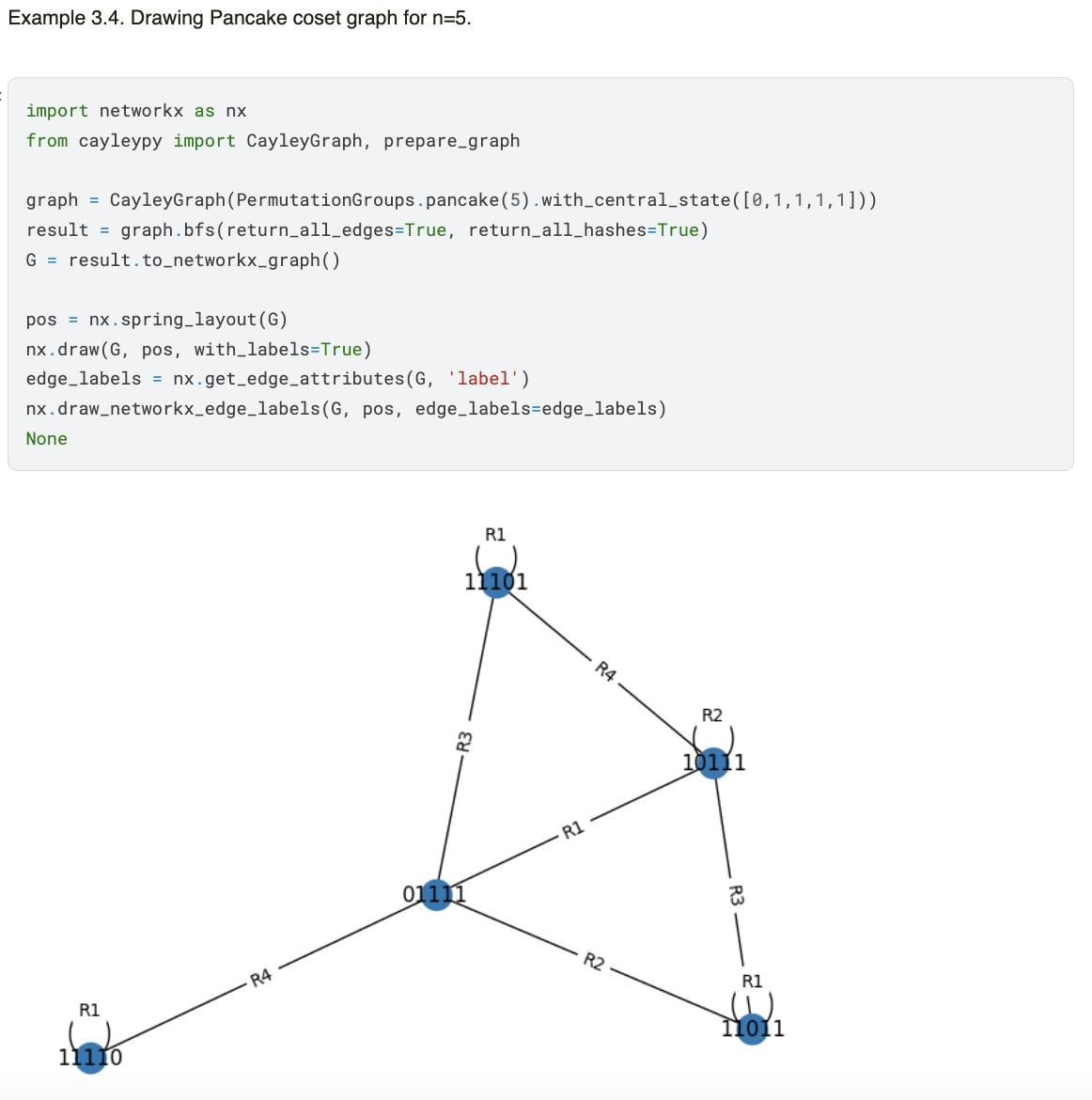}
  \caption{Defining a Schreier coset graph by specifying the initial vector ("central state" -- the vector whose stabilizer defines a subgroup $H$ for the factor set $G/H$),  computing the adjacency matrix, 
  visualization.}
  \label{fig:tutotialEx3}
\end{figure}

Figure \ref{fig:tutotialEx3} presents an example on how to define a Schreier coset graph.
It is not defined by the subgroup explicitly, but instead by specifying the vector whose stabilizer is the desired subgroup. We call that vector the "central state". The graph is constructed by applying generators first to it,
then to the previously obtained vectors, and so on, i.e. the standard BFS (\href{https://en.wikipedia.org/wiki/Breadth-first_search}{breadth first search}) algorithm with an efficient (not so straightforward) implementation. 

Further features of CayleyPy: 
\begin{enumerate}
    \item A collection of predefined permutation group Cayley graphs:  \href{https://cayleypy.github.io/cayleypy-docs/generated/cayleypy.PermutationGroups.html#cayleypy.PermutationGroups}{link}.
    \item Groups originating from puzzles:
    \href{https://cayleypy.github.io/cayleypy-docs/generated/cayleypy.Puzzles.html#cayleypy.Puzzles}{link},
    \href{https://cayleypy.github.io/cayleypy-docs/generated/cayleypy.GapPuzzles.html}{link}.
    \item Predefined matrix groups: \href{https://cayleypy.github.io/cayleypy-docs/generated/cayleypy.MatrixGroups.html#cayleypy.MatrixGroups}{link}. 
    \item Random walks: \href{https://cayleypy.github.io/cayleypy-docs/generated/cayleypy.algo.RandomWalksGenerator.html#cayleypy.algo.RandomWalksGenerator}{link}.
    \item Components for AI-based path-finding (to be improved in future releases): 
    \href{https://cayleypy.github.io/cayleypy-docs/generated/cayleypy.algo.BeamSearchAlgorithm.html#cayleypy.algo.BeamSearchAlgorithm}{Beamsearch}, \href{https://cayleypy.github.io/cayleypy-docs/api.html}{ML-models}. 
\end{enumerate}

\clearpage
\subsubsection{Benchmarks for growth computations}
CayleyPy utilizes GPU out of the box, for example,  
even on old GPU P100 it can compute growth for $S_{11}$ with Coxeter generators  in 0.6 seconds,
while GAP takes 2352 seconds, so CayleyPy is more than 1000 times faster in that example.
Performance depends on group size and number of generators. The tables below demonstrate that for large groups starting from $S_8$
CayleyPy is typically 10-100 times faster than GAP, even when using CPU. Moreover, it can support larger groups. Below are the results 
obtained and easily reproducible on the Kaggle cloud where we can perform computations up to $S_{13}$. 
For $S_{14}$, $S_{15}$ we use more powerful machines. The group $S_{14}$ requires around $40-100$ GB RAM and $4-20$ hours of computations.
The group $S_{15}$ requires $500-1000$ GB RAM, and computations take several days.  (For $S_{15}$ we mainly worked with a small number of generators like 3, while a larger number of generators (e.g. 15) may take months.) 

Tomas Rokicki and Lucas Garron's program  \href{https://github.com/cubing/twsearch}{"Twsearch"} (\href{twsearch}{example on Kaggle})
apparently is faster than our growth computations for CPU. However our code supports  GPU out of the box, and apparently
can achieve better timing using modern GPU. 
Also our framework seems to be more user-friendly, supports directed Cayley graphs and achieves computations for large groups like $S_{15}$, which, apparently, is not yet achieved by Twsearch. 

CayleyPy supports several algorithms for growth computation all based on BFS (\href{https://en.wikipedia.org/wiki/Breadth-first_search}{breadth first search}), but different in internal
data representation.
The \href{https://cayleypy.github.io/cayleypy-docs/generated/cayleypy.algo.bfs_bitmask.html#cayleypy.algo.bfs_bitmask}{bfs bitmask} uses bit-wise encoding with 3 bits per any state,
and it is more memory efficient. It allows one to work with $S_{13}$ requiring only $3-8$ GB RAM. We use it for computations for $S_{13},S_{14},S_{15}$.
See also the notebooks: \href{https://www.kaggle.com/code/fedimser/memory-efficient-bfs-on-caley-graphs-3bits-per-vx}{algorithm},
 \href{https://www.kaggle.com/code/alexandervc/cayleypy-timing4-bfs-bitmask-memory-effecient}{benchmarks}.

In the tables below, the label "(bm)" identifies uses of the "bfs-bitmask" algorithm.

\begin{table}[ht]
\caption{Growth computations. Time in seconds for CayleyPy and GAP using CPU on Kaggle cloud, 32 GB RAM.
Different types of generators.}\label{tab:timing1}

\begin{tabular}{|c|c|c|c|c|c|c|c|c|c|c|}
\hline
 & \multicolumn{2}{|c|}{LRX} &  \multicolumn{2}{|c|}{Coxeter}  & \multicolumn{2}{|c|}{Transpositions} & \multicolumn{2}{|c|}{Pancake}  & \multicolumn{2}{|c|}{Reversals}  \\
n & GAP  & CayleyPy & GAP  & CayleyPy & GAP  & CayleyPy & GAP  & CayleyPy & GAP  & CayleyPy \\
\hline
6 & 1.97 & 0.02 & 0.01& 0.038& 0.03 & 0.02& 0.01& 0.02& 0.04& 0.02\\
\hline
7 & 2.04 & 0.03 & 0.12& 0.052& 0.34& 0.03& 0.08& 0.02& 0.41& 0.05\\
\hline
8 & 2.66 & 0.05& 1.32 & 0.079& 3.92& 0.09& 0.80& 0.04& 4.24& 0.1\\
\hline
9 & 8.75 & 0.16& 13.99& 0.293& 52.98& 1.03& 8.90& 0.26& 54.45& 1.089\\
\hline
10 & 75.72 & 1.24& 172.38& 3.54& 737.42 & 14.14& 112.37& 3.77& 762.54 & 14.8\\
\hline
11 & 884& 22.5& 2352& 41& 10891&   296   & 1559& 72& 11138& 354\\
\hline
12 & 12024& 658  & 34470&  706 & - & 3331(bm) & 35387& 4923 & - & 26657  \\
\hline
13 & - & 2670(bm) & - & 7308(bm) & -  & >12h & - & 17694(bm) & - & >12h \\
\hline
\end{tabular}
\end{table}

\begin{table}[ht]
\caption{CayleyPy growth computations on different hardware. Time in seconds: CPU (32G), GPU (16G), and advanced CPU (330G), on Kaggle cloud. Coxeter generators.}\label{tab:timing2}

\begin{tabular}{|c|c|c|c|c|}
\hline
\multicolumn{5}{|c|}{Coxeter CayleyPy} \\
\hline
n & CPU & GPU T4 & GPU P100 & CPU (at TPU v3-8) \\
\hline
4 & 0.013 & 0.013 & 0.014 & 0.010 \\
\hline
5 & 0.029 & 0.026 & 0.026 & 0.022 \\
\hline
6 & 0.038 & 0.043 & 0.043 & 0.037 \\
\hline
7 & 0.052 & 0.068 & 0.067 & 0.067 \\
\hline
8 & 0.079 & 0.151 & 0.138 & 0.087 \\
\hline
9 & 0.293 & 0.110 & 0.110 & 0.140 \\
\hline
10 & 3.539 & 0.199 & 0.154 & 0.641 \\
\hline
11 & 41/41(bm) & 1.085 & 0.601 & 10.8/12(bm) \\
\hline
12 & 705/512(bm) & 214(bm)  & 207(bm) & 187/151(bm)\\
\hline
13 & 7308(bm) & 3004(bm)  & 2867(bm) & 2613/2180(bm) \\
\hline
\end{tabular}
\end{table}

Notebooks: 
GAP testing has been performed in this \href{https://www.kaggle.com/code/avm888/group-elements-decomposition-gap-limits/}{notebook}. 
CayleyPy: \href{https://www.kaggle.com/code/alexandervc/cayleypy-timing2-coxeter}{Coxeter}, 
\href{https://www.kaggle.com/code/alexandervc/cayleypy-timing3-several-graphs}{other generators},
\href{https://www.kaggle.com/code/alexandervc/cayleypy-timing4-bfs-bitmask-memory-effecient}{bfs bitmask algorithm}.


\subsubsection{Benchmarks for AI path-finding}
The main focus of the present release is not the AI-component,
but for the sake of completeness let us reproduce the comparison with GAP
from our previous paper \cite{CayleyPyRL}. AI methods easily outperform GAP, and in particular produce the optimal solution
even for graph sizes where GAP fails  to find any solutions at all. We were also able to find non-optimal solutions beyond $S_{100}$
for these particular generators, which corresponds to graphs of more than googol size that are far beyond GAP's capabilities.

\begin{table}[ht]
    \centering
    \begin{tabular}{|c|c|c|c|c|c|}
        \hline
        $n$  &  GAP length   & $\frac{n(n-1)}2$            &   DD length    \\
        \hline
        9    &  41     & 36               &  36   \\
        \hline
        10   &    51   & 45               &  45       \\
        \hline
        11   &     65      & 55             &  55     \\
        \hline
        12   &    78       &  66          & 66     \\
        \hline        
        13   &    99       &  78           &   78       \\
        \hline        
        14   &   111       &  91           &   91     \\
        \hline        
        15     &   268      &  105         &  105    \\
        \hline        
        16     &    2454    &  120         &  120     \\
        \hline        
        17   &     380     &  136         & 136    \\
        \hline        
        18   &     20441   &   153        & 153    \\
        \hline        
        19   &      3187      & 171        & 171    \\
        \hline      
        20   &     217944    & 190        & 190      \\
        \hline 
        21   &   -           & 210     & 210        \\
        \hline           
    \end{tabular}
    \captionsetup{skip=10pt} 
    \caption{Comparison of GAP and diffusion distance (DD) method for the conjecturally longest element of LRX Cayley graph for $S_n$ with expected length $n(n-1)/2$.}
    \label{tab:example}
\end{table}

For $S_{20}$, the  GAP timing is 41min 18s,  while our methods can find results much faster.
For example, the diffusion distance method with a simple neural network requires only 3m 48s of calculations on GPU P100 (\href{https://www.kaggle.com/code/alexandervc/lrx-cayleypy-rl-mdqn?scriptVersionId=224270083}{notebook version 423}). The GAP benchmarks have been performed \href{https://www.kaggle.com/code/avm888/group-elements-decomposition-gap-limits}{here}.

\clearpage
\subsection{Diameter quasi-polynomiality conjecture }

Having discussed the main features of our library CayleyPy, let us now turn to our main mathematical results. We start with the following definition. A function 
\( f \colon \mathbb{N} \to \mathbb{N} \) 
is a \href{https://en.wikipedia.org/wiki/Quasi-polynomial}{quasi-polynomial} if there exist polynomials 
\( p_{0}, \dots, p_{s-1} \) 
such that 
\[
   f(n) = p_{i}(n) \quad \text{when } i \equiv n \pmod{s}.
\]
The polynomials \( p_{i} \) are called the \emph{constituents} of \( f \).

The following conjecture  generalizes  results of extensive computations we performed (see the next section):

\begin{Conj}(Extremely optimistic). For any generators of $S_n$ (or $A_n$) which can be constructed by an algorithm
with say polynomial complexity in $n$
(e.g. a Python function which takes as input $n$ and outputs generators in time polynomial in $n$) the diameter of the Cayley graph 
will be given by some quadratic or linear quasi-polynomial in $n$ (at least for $n$ large enough). Even more optimistically, the leading terms of all constituents  coincide.

\end{Conj}

Experimentally we see an even more general phenomenon: not only the distance to the longest element (i.e. the diameter) is quasi-polynomial, but also the distance to many other elements. Thus, it is tempting to propose the following:

\begin{Conj}(Extremely optimistic). For any generators of $S_n$ (or $A_n$) and additionally elements $g_n \in S_n$ (or $A_n$) which can be constructed by an algorithm with say polynomial complexity in $n$
(e.g. a Python function which takes as input $n$ and outputs generators jointly with elements $g_n$ in time polynomial in $n$) 
the distance from identity to $g_n$  (i.e. "\href{https://en.wikipedia.org/wiki/Word_metric}{word metric}" of $g_n$) 
will be given by some quadratic or linear quasi-polynomial in $n$,  at least for $n$ large enough. Even more optimistically, the leading terms
of all constituents  coincide.
\end{Conj}

The two conjectures do not seem to imply each other in general, however in practice we often see that the longest elements
("antipodes", "superflips") can be often described quite effectively. In such cases the second conjecture implies the first one.

Generalizations to other groups are also plausible,  e.g. to Coxeter's groups and the \href{https://en.wikipedia.org/wiki/Generalized_symmetric_group}{generalized symmetric group} (examples will be given below). Even for  matrix groups one may hope to have a similar quadratic/linear quasi-polynomial dependence on $n$ (i.e. on the "rank").  
And similarly for Schreier coset graphs, for example  $S_{n}/(S_{\lfloor n/\ell\rfloor}\times S_{n-\lfloor n/\ell\rfloor})$ (i.e. "Grassmanians over the field with one element"),
and other similar cosets like flags over $F_1$. However, in that case we typically choose some node and compute the most distant node to it (not exactly the diameter), and in contrast to a Cayley graph that distance ("God's number" in puzzle's terminology) can depend on the choice of the node. We expect quasi-polynomiality for all choices of the starting node.
Another direction of generalization is weighted Cayley graphs, in particular
circular Cayley graphs, i.e. permutations factorized by cyclic shifts 
(see e.g. \cite{alon2025circular}), which  represent weighted Cayley graphs with weights of cyclic shifts set to zero.

If the conjecture (or any of its weaker versions) is indeed true, it would be quite surprising, 
because computing the diameter is known to be NP-hard \cite{even1981minimum} in general,
and extremely hard in many particular cases, like for example the Rubik's cube where it took 40 years to achieve it
\cite{Rokicki2014Diameter}. Nevertheless, there seems to be no immediate contradiction since 
determining the actual quasi-polynomials and the modulo $s$ might be NP-hard. But from a practical side it opens a rather effective way to compute diameters -- 
one needs to compute diameters for some values of $n$ and then just try to fit a quasi-polynomial.
Computations up to $n\le 15$ are often achievable by e.g. effective implementations of brute force BFS (breadth first search) 
algorithms such as those provided in CayleyPy (for cosets above we can even compute up to $n\le 42$). 
Further sizes are hard to access by direct methods, however there are ways around this discussed below. 

Taking into account the amount of examples  confirming the conjecture, it is natural to believe that
it is true in one or another way. It might be that one needs to restrict  the class of generators. 
An example quite interesting to us is to consider three (or more) generators which are involutions
and which come from some natural construction for all $n$ (like many examples in the present text) -- will  the
conjecture be true in this case? 

All known to us examples of the explicitly computed diameters are indeed quasi-polynomial,
although they may be written in a slightly different way in the literature, e.g. in terms of  rounding functions floor or ceil which are just particular simple cases of quasi-polynomials. In our analysis we observe more complicated examples with the modulo $s$ equal to e.g. 4, 6, or 8. 

{\bf Examples (known)}. For Coxeter generators ( $(i,i+1), 1\le i<n$) of $S_n$ the diameter is $n(n-1)/2$  -- just polynomial (and similarly for all Coxeter's groups). As a small modification one can add the transposition $(1,n)$ to the Coxeter generators -- this gives the diameter  $\left\lfloor \tfrac{n^2}{4}\right\rfloor$
\cite{vanZuylen2016cyclic}, which is $n^2/4$ for even $n$, and $(n^2-1)/4$ for odd $n$, i.e. quasi-polynomial with modulo
2.  (A similar  (but very recent \cite{alon2025circular}) circular version gives $\left\lfloor \tfrac{(n-1)^2}{4}\right\rfloor$). 
For all transpositions $(i,j)$ the diameter is $n-1$, while for transpositions of the form $(1,i)$ (star graph)  it is 
 $\left\lfloor \tfrac{3(n-1)}{2}\right\rfloor$, which is $3(n-1)/2$ for odd $n$, and $(3n-2)/2$ for even -- a linear quasi-polynomial. 

{\bf Examples (conjectural)}. Consider the LRX generators: $L$ -- left cyclic shift, $R$ -- right cyclic shift, $X =(1,2)$. The diameter is conjectured to be $n(n-1)/2$ (\href{https://oeis.org/A186783}{OEIS-A186783}), strong support for this is presented in our previous work \cite{CayleyPyRL}.
The LX case with only two of these generators ($L$ and $X$) has been first considered by V.M.Glushkov \cite{glushkov1968completeness}
and studied by his school (survey: \cite{glukhovzubov1999lengths} pages 18-21). 
Our conjecture is that the diameter is $(3n^2-8n+9)/4$ for $n$ odd, and $(3n^2-8n+12)/4$ for $n$ even, discussed below in detail.

Consider consecutive 4-cycles: $(i,i+1,i+2,i+3),i\le 1\le n-3$, and their action on the coset $S_{n}/(S_{\lfloor n/2\rfloor} \times S_{n-\lfloor n/2\rfloor})$ which is just the action on  binary vectors with 0 and 1 having $\lfloor n/2\rfloor$ zeros. Choose the vector with first $\lfloor n/2\rfloor$ zeros 
as the starting vector and compute the most distant element with respect to that initial vector (it is not exactly the diameter
in general, but its analog depending on choice of initial node, and can be called "God's number" like in puzzles). 
We expect that for $n\ge 12$ the God's number is given by quasi-polynomials modulo 6: $n^2/12+1, n \equiv 0 \mod 6$, $(n^2+4n-5)/12, n \equiv 1 \mod 6$, etc.

\subsection{Table of results: hundreds of new conjectures on Cayley graphs}

We conducted extensive computations of growth for a large number of Cayley and Schreier coset graphs, up to $n\le 15$ and $n\le 42$ respectively. 
The obtained results and conjectures are summarized in the tables discussed below, which also include results known in the literature. 

We analyzed not only diameters but other growth characteristics as a probability distribution -- mean, mode, variance,
skewness, kurtosis, and also tried to fit the distribution by some known shape like Gaussian or Gumbel. We also looked for antipodes (longest elements, or super-flips) and explored the spectra of the Cayley graphs. In some cases we observe that growth by itself might 
have a closed analytical formula -- given by Stirling or Fibonacci numbers, or e.g. coinciding with 
some known sequence such as \href{https://oeis.org/A367270}{OEIS-A367270} ("Growth/F-la" column of the table). 
For diameters in most (but not all) cases we were able to fit the data by quasi-polynomials, while for some cases, apparently, the
available data is not enough. But  for mean diameters and other characteristics of growth
the formulas are not quasi-polynomials in general, for example in the literature there are results of the form $n-\log(n)$ for mean diameters. Nevertheless we expect that our numerical fits for the data provide approximations 
to the leading terms of these characteristics. They are obtained as a fit for small values of $n$ and can be 
considered as conjectures for large values.

Notation used in the table:
\begin{enumerate}
    \item \ding{117} -- conjecture obtained by CayleyPy project, \textcolor{blue}{\ding{117}} -- proved by CayleyPy
    \item $+$ -- information is known and can be found in the main text (too big to fit into table)
    \item $*$ -- conjecture from the literature 
    \item ? -- no information, neither the literature nor our experiments suggest a clear pattern 
    \item notation like $1|2$ indicates 1 for even $n$ and $2$ for odd 
    \item $+I$ -- some quasi-polynomial typically of zero degree 
    \item "Group" -- information on the generated group;  just '$+$' means information is known but does not fit into the table
    \item "Growth/PDF" -- what continuous distribution fits growth for large $n$ (Gaussian, or Gumbel, etc)
    \item "Growth/F-la" -- explicit formulas for the growth 
    \item "Antipode" --  information on longest elements, i.e. if there is an explicit description or 
    if the number is known (and simple to fit into table) we indicate it, or simply write '$+$'   
    \item "Algorithm" -- 
    is there a known algorithm to decompose an element into a product of these generators.
    The superscript $O$ indicates that an optimal algorithm is known, notation like $NP/2$ means
    that the optimal decomposition was proven to be NP-hard, but there are polynomial approximations up to a factor of 2.
    \item "Metric" -- is there an explicit expression for the word metric for given generators, for example for Coxeter generators
    it is the number of inversions.
    \item "Spectrum" -- information on the spectrum. We use labels "Int" for integer spectrum, "Wig" for spectrum approaching a Wigner semi-circle law for large $n$,
    "Uni" for an almost uniform spectrum.
\end{enumerate}

{\bf Generators naming.} The names of the generators either correspond to their naming in CayleyPy library (\href{https://cayleypy.github.io/cayleypy-docs/generated/cayleypy.PermutationGroups.html#cayleypy.PermutationGroups}{link}),
or can be hopefully be easily guessed e.g. $(i,i+1,i+2)$ denotes generators consisting of consecutive 3-cycles sending $i\to i+1 \to i+2 \to i $, for $i<=n-3$. The full version of manuscript presents separate chapter on each  generator and describing expriments and results in details. But for the sake of completeness let us explicitly describe some namings, provided links leads to CayleyPy documentation for these generators: 
\begin{enumerate}
    \item "\href{https://cayleypy.github.io/cayleypy-docs/generated/cayleypy.PermutationGroups.html#cayleypy.PermutationGroups.coxeter}{Coxeter}" - standard Coxeter generators of $S_n$, i.e. neighbor transpositions $(i, i+1)$
    \item "\href{https://cayleypy.github.io/cayleypy-docs/generated/cayleypy.PermutationGroups.html#cayleypy.PermutationGroups.cyclic_coxeter}{Cyclic Coxeter}" - add cyclicity: add $(1,n)$ generator to Coxeter's ones  (studied in \cite{jerrum1985complexity}, and see references in \cite{alon2025circular})
    \item "\href{https://cayleypy.github.io/cayleypy-docs/generated/cayleypy.PermutationGroups.html#cayleypy.PermutationGroups.lx}{LX}/\href{https://cayleypy.github.io/cayleypy-docs/generated/cayleypy.PermutationGroups.html#cayleypy.PermutationGroups.lrx}{LRX}" - here "L" is left cyclic shift, "R" - right cyclic shift, "X" transposition of the first two elements,
    studied a lot: starting at least from V.M.Glushkov paper 1968 (see section \ref{sect:Glushkov}).
    \item "\href{https://cayleypy.github.io/cayleypy-docs/generated/cayleypy.PermutationGroups.html#cayleypy.PermutationGroups.larx}{LARX}" - X  is again transposition of the first two elements, but cyclic shifts here are NOT maximal,  shifting positions from the second to the last (first position is stable). Acronym from "left Almost right and X".
    Respectively LARX+I - adding inverse to the cycle, i.e. right shift. 
    \item "LSL" - Long and Sub-Long cycles, i.e. full cyclic shift and cyclic shift rotating only positions from the second to the last keeping first one stable. Respectively LSL+I - adding inverses
    \item "\href{https://cayleypy.github.io/cayleypy-docs/generated/cayleypy.PermutationGroups.html#cayleypy.PermutationGroups.three_cycles}{3-cyc}" - all 3-cycles 
    \item "\href{https://cayleypy.github.io/cayleypy-docs/generated/cayleypy.PermutationGroups.html#cayleypy.PermutationGroups.three_cycles_01i}{(0ij)}" - 3-cycles having form (0ij)
    \item "\href{https://cayleypy.github.io/cayleypy-docs/generated/cayleypy.PermutationGroups.html#cayleypy.PermutationGroups.consecutive_k_cycles}{(i,i+1,i+2)}" - 
    consecutive 3-cycles sending $i\to i+1 \to i+2 \to i $, for $i<=n-3$, "(i,i+1,i+2)I" - same with inverses,
    similar "(i,i+1,i+2,i+3)", "(i,i+1,i+2,i+3,i+4)", consecutive 4,5-cycles, etc. all supported in CayleyPy by the same
    function \href{https://cayleypy.github.io/cayleypy-docs/generated/cayleypy.PermutationGroups.html#cayleypy.PermutationGroups.consecutive_k_cycles}{consecutive\_k\_cycles}
    \item "(i,i+1,i+2)C" - consecutive 3-cycles with cyclically wrapped across "n", that means include
    cycles (n-1,n,1), (n,1,2), similar "(i,i+1,i+2,i+3)C" ,"(i,i+1,i+2,i+3,i+4)C" - 4-cycles, 5-cycles, etc.
    all supported in CayleyPy by the same function:
    \href{https://cayleypy.github.io/cayleypy-docs/generated/cayleypy.PermutationGroups.html#cayleypy.PermutationGroups.wrapped_k_cycles}{wrapped\_k\_cycles}, "I" at the end means inverses are added
    \item "\href{https://cayleypy.github.io/cayleypy-docs/generated/cayleypy.PermutationGroups.html}{Pref.cyc}" - prefix cycles i.e. $(1, 2, ..., i), 1<i\le n$
    \item "\href{https://cayleypy.github.io/cayleypy-docs/generated/cayleypy.PermutationGroups.html}{Down.cyc}" - downcycles: $(i, i+1, ..., j), i < j$
    \item "Inc.3cyc", "Inc.4cyc"... 
    supported in CayleyPy: increasing cycles $ (i_1, i_2, ..., i_k), i_1 < i_2 < ... <i_k$
    \href{https://cayleypy.github.io/cayleypy-docs/generated/cayleypy.PermutationGroups.html#cayleypy.PermutationGroups.increasing_k_cycles}{increasing\_k\_cycles}
    \item 
    \href{https://cayleypy.github.io/cayleypy-docs/generated/cayleypy.PermutationGroups.html#cayleypy.PermutationGroups.pancake}{Pancake} - \href{https://en.wikipedia.org/wiki/Pancake_sorting}{prefix (pancake) sorting} generators \cite{GatesPapadimitriou1979PrefixReversal}
    \item "\href{https://cayleypy.github.io/cayleypy-docs/generated/cayleypy.PermutationGroups.html#cayleypy.PermutationGroups.transposons}{Transposons}",
    "\href{https://cayleypy.github.io/cayleypy-docs/generated/cayleypy.PermutationGroups.html#cayleypy.PermutationGroups.full_reversals}{Reversals}","\href{https://cayleypy.github.io/cayleypy-docs/generated/cayleypy.PermutationGroups.html#cayleypy.PermutationGroups.signed_reversals}{signed Reversals}" 
    biologically motivated generators described in: \ref{sec:bio-generators}
    \item "\href{https://cayleypy.github.io/cayleypy-docs/generated/cayleypy.PermutationGroups.html#cayleypy.PermutationGroups.rapaport_m1}{RapaportM1}",
     "\href{https://cayleypy.github.io/cayleypy-docs/generated/cayleypy.PermutationGroups.html#cayleypy.PermutationGroups.rapaport_m2}{RapaportM2}" - generators from \cite{RapaportStrasser1959Cayley}, see also \cite{glukhovzubov1999lengths}
    \item "\href{https://cayleypy.github.io/cayleypy-docs/generated/cayleypy.PermutationGroups.html#cayleypy.PermutationGroups.cubic_pancake}{3Pancake S1,S2,...S7}" - 
    cubic Pancake graphs considered by Bass and Sudborough~\cite{BS03}, and "S7" in \cite{Gun22}. 
    \item "\href{https://cayleypy.github.io/cayleypy-docs/generated/cayleypy.Puzzles.html#cayleypy.Puzzles.globe_puzzle}{Globes}"
    - introduced in Kaggle Santa 2023 challenge see e.g. "\href{https://www.kaggle.com/code/marksix/visualize-allowed-moves}{notebook}",
    similar to \href{https://www.jaapsch.net/puzzles/master.htm}{"Masterball"} puzzle, but not exactly it. 
    
\end{enumerate}

{\footnotesize 
\begin{longtable}{|@{}p{1.7cm}|p{0.4cm}|@{}p{1.7cm}|@{}p{1.1cm}
|@{}p{0.4cm}|@{}p{1cm}|p{1cm}|p{0.7cm}|p{0.7cm}|p{0.7cm}|p{0.9cm}|@{}p{0.8cm}
|p{0.5cm}|@{}p{0.7cm}|@{}p{1cm}|}

\caption{Summary of properties of Cayley graphs}\label{tab:long1}\\

\hline
\textbf{Gene-} & \textbf{Gro-} & \textbf{Dia-} & \multicolumn{7}{|c|}{\textbf{Growth}} & \textbf{Anti-}  & \textbf{Algo-} & \textbf{Met-} & \textbf{Spec-} & \textbf{Mixing}  \\ 
\textbf{rators} & \textbf{up} & \textbf{meter} & \textbf{PDF} & \textbf{F-la} &  \textbf{Mean}  & \textbf{Mode} & \textbf{Var} & \textbf{Skew} & \textbf{Kurt} &  \textbf{podes} & \textbf{rithm} &  \textbf{ric} & \textbf{trum}  & \textbf{Time} \\
\hline
\endfirsthead

\hline
\textbf{Gene-} & \textbf{Gro-} & \textbf{Dia-} & \multicolumn{7}{|c|}{\textbf{Growth}} & \textbf{Anti-}  & \textbf{Algo-} & \textbf{Met-} & \textbf{Spec-} & \textbf{Mixing}  \\ 
\textbf{rators} & \textbf{up} & \textbf{meter} & \textbf{PDF} & \textbf{F-la} &  \textbf{Mean}  & \textbf{Mode} & \textbf{Var} & \textbf{Skew} & \textbf{Kurt} &  \textbf{podes} & \textbf{rithm} &  \textbf{ric} & \textbf{trum}  &  \textbf{Time} \\
\hline
\endhead

\hline
\multicolumn{15}{|r|}{\textit{Continued on next page}} \\
\hline
\endfoot

\hline
\multicolumn{15}{|r|}{\textit{End of table}} \\
\hline
\endlastfoot

Coxeter & $S_n$ & $\frac{n(n-1)}{2}$ & Gauss & + & $\frac{n(n-1)}4$  & $\frac{n(n-1)}4$ 
& \hfil $+$  & \mbox{ $\to 0$}   & \mbox{ $\to 0$}   & 1 & Bubble  & $+$   &   \hfil ?   & \hfil ?    \\ \hline
Cyclic Coxeter&  $S_n$ & $\left\lfloor \tfrac{n^2}{4}\right\rfloor$
  & Gauss\ding{117} & ? & $0.17(n^2-n+1)$\ding{117} &  \mbox{$\approx$Mean\ding{117}}  & \hfil ?    & \mbox{ $\to 0$\ding{117}} &  \mbox{ $\to 0$\ding{117}}  & \mbox{1$|$2\ding{117}}  &  \hfil ?   & +  & \mbox{Wig\ding{117}} & \hfil ?    \\ \hline
LRX & $S_n$ & $\frac{n(n-1)}{2}*$ & \mbox{Gumbel\ding{117}} & ?  & $\approx 0.38n^2-n$\ding{117}  &  $\approx 0.39n^2-n$\ding{117}   
& \hfil ?  &  $\to -0.7$\ding{117}  & $\to 3.3$\ding{117}    &  \hfil \ding{117} & \hfil \mbox{\ding{117}}   & ? & Uni\ding{117}  &  $>n^3$\ding{117}  \\ \hline
LX-Glushkov & $S_n$  & \raisebox{-5pt}{\mbox{ $\frac{ 3n^2 - 8n  + 9|12}{4}$\ding{117}}} 
&  \mbox{Gumbel\ding{117}} & Fib /?  & $\approx 0.57n^2-2n$\ding{117}   & $\approx 0.57n^2-1.6n$\ding{117}   & \hfil ?  &$\to -0.7$\ding{117}   &  $\to 0.5$\ding{117} & 
 \hfil \raisebox{-5pt}{ *  }  &   \hfil \raisebox{-5pt}{ \ding{117} } & \hfil ?   & \hfil ?  &  \hfil ?   \\ \hline
 %
LARX  & $S_n$   & \mbox{$\frac{n^2-2|5}{2}$\ding{117}}   
&  \hfil ?  &  ?  
&   $0.4n^2 - 0.7n$ \ding{117}   & $\approx$Mean\ding{117}   & \hfil ?  
&  \hfil ? & \hfil ? 
&   \mbox{+\ding{117}$|$?} & \hfil ? 
&  \hfil ? & \hfil ? & \hfil ?  \\ \hline
%
LARX+I  &  $S_n$   &   \mbox{$\frac{n(n+6)-12|19}{4}$\ding{117}}  
&   \hfil ? &  \hfil ?   
&  {$\approx \frac{n(n+1)}{4}$\ding{117}}   & $\approx$Mean\ding{117}  & \hfil ?  
&  \hfil ? & \hfil ? 
& $+$\ding{117} & \hfil ? 
&  \hfil ? &  \hfil ? & \hfil ?  \\ \hline
LSL & $S_n$ &  \mbox{$\tfrac{n(n-3)}{2}+3$\ding{117}}  &  \hfil ?  & ?    
& $\approx  0.4n^2-1.5n$\ding{117}  &   \mbox{$\approx$Mean\ding{117}}  & \hfil ? 
&   \hfil ? &  \hfil ? 
& $+$\ding{117} & \hfil ?  &  \hfil ?   
&  \hfil ?  & \hfil ?    \\ \hline
%
LSL+I & $S_n$  & \mbox{$\tfrac{n(n+4)}{4}-3|\frac{17}4$\ding{117}}   &  \hfil ?  & ?    
&  $\approx 0.2n^2$ & $\approx 0.2n^2$   & \hfil ?  & \hfil ?   & \hfil ?  
& \hfil ?  & \hfil ? &   \hfil ? 
& \hfil ?  &  \hfil ?  \\ \hline
%
3-cyc & $A_n$ &   \raisebox{-5pt}{$\left\lfloor \tfrac{n}{2}\right\rfloor$}  &  \hfil ?  & $+-$    
& $\approx D-0.5|1$\ding{117}  & $D-0|1$\ding{117}  &  \hfil ? 
&  \hfil ? & \hfill ?  
&   $+$   & $+$  
&  $+$  & Int  &  \hfil ? \\ \hline
%
%
(0ij) &  $A_n$  &  $\left\lfloor \tfrac{3(n-1)}{4}\right\rfloor$  &  \hfil ?  & \hfil ?    & $\approx 0.55n $\ding{117}   & $\approx Mean $\ding{117}   & \hfil ? 
&  \hfil ?  & \hfil ?  
&  \hfil ?  & \hfil ?  
&  \hfil ?  & \hfil ?  & \hfil ?   \\ \hline
%
 (01i)  &   $A_n$  & $\tfrac{3n-5}2+\tfrac{i^n+(-i)^n}4$\ding{117}    
 & \hfil ? &  \hfil ?    
&  $\approx n-2$\ding{117} &  $n-1$\ding{117}  & \hfil ?  
&  \hfil ?  & \hfil ?  
& \hfil ?  & $+$\textcolor{blue}{\ding{117}}  
& \hfil ?   & \hfil ?  & \hfil ?   \\ \hline
(01i)I&  $A_n$  & \raisebox{-5pt}{$\left\lfloor \tfrac{3n-6}{2}\right\rfloor$ \textcolor{blue}{\ding{117}}}  & \hfil ? & \hfil ? & 
{$\approx n+1.25\ln n+...$\ding{117}} & \mbox{$\approx$Mean\ding{117}}  & \ding{117}  & \ding{117}  & \ding{117} & \ding{117} &  \textcolor{blue}{\ding{117}}$^{O}$  & \textcolor{blue}{\ding{117}} & Int\textcolor{blue}{\ding{117}} & \hfil nlog(n)\ding{117}    \\ \hline
%
\mbox{(i,i+1,i+2)} &  $A_n$   
&  \mbox{ $ \left\lfloor \tfrac{n^2+1}{4}\right\rfloor${\ding{117}}}  &  Gauss\ding{117}  & ?    & \mbox{$\approx D/2$\ding{117}}   & \mbox{$\approx D/2$\ding{117}}    
& $\approx \tfrac{n^3}{100}$ &   \mbox{ $\to 0$\ding{117}} 
& \mbox{ $\to 0$\ding{117}}  
&  \hfil ?  & \hfil ? 
&  \hfil ? & \hfil ?  &  \hfil ? \\ \hline
%
\mbox{(i,i+1,i+2)I} &  $A_n$   & \mbox{$\left\lfloor \tfrac{n(n-1)}{4}\right\rfloor$\textcolor{blue}{\ding{117}}} &  Gauss{\ding{117}}   & ?    
& \mbox{$\approx D/2$\ding{117}}   &  \mbox{$\approx D/2$\ding{117}}   
& $\approx \tfrac{n^3}{100}$\ding{117}  
&   \mbox{ $\to 0$\ding{117}}  &    \mbox{ $\to 0$\ding{117}}  
& $+$\ding{117}  & \hfil ? &  \hfil ? & \hfil ?  &  \hfil  ?   \\ \hline
%
\mbox{(i...i+3)} &  $S_n$   &   $\approx 0.3n^2$\ding{117}
 &  \hfil ?  & ?    
 & $\approx 0.16n^2$\ding{117}  &  $\approx 0.15n^2$\ding{117}  
 &   \hfil ?&  $\approx 0$\ding{117} & $\approx 0$\ding{117}  
 &  \hfil ? &  \hfil ? &  \hfil ?  
 & \hfil ? &   \hfil ? \\ \hline
\mbox{(i...i+3)I} &  $S_n$   &  $(n(n-1)+ 1|1|4)/6$\ding{117} 
 &  \hfil ?  & \hfil ?    
 &  $\approx 0.036n^2$\ding{117}   &  $\approx 0.045n^2$\ding{117}   
 & \hfil ?  &  \hfil ?  &   \hfil ?    
 &  \ding{117}   & \hfil ?  &  \hfil ?  
 & \hfil ?  & \hfil ?   \\ \hline
%
%
 \mbox{(i,i+1,i+2)C} &  $A_n$   & \mbox{ $\approx \tfrac{n(n+2)}{8}+I$}{\ding{117}}
 &  \hfil ?  & ?    
 & $\approx 0.085n^2$\ding{117}   &   $\approx 0.065n^2$\ding{117}  
 & \hfil ? &\hfil ?   & \hfil ?  
 &    & \hfil ? & \hfil ?   
 & \hfil ? & \hfil ?  \\ \hline
 \mbox{(i,i+1,i+2)CI} &  $A_n$   &  \mbox{$\left\lfloor \tfrac{n^2}{8}\right\rfloor$\textcolor{blue}{\ding{117}}}  
 &  \hfil ?  &  \hfil ?    
 & $\approx 0.08n^2$\ding{117}   & $\approx 0.086n^2$\ding{117}   
 & \hfil ?   & \hfil ?    & \hfil ?    
 & $+$\ding{117} & \hfil ?   &  \hfil ?   
 & \hfil ?   & \hfil ? \\ \hline
 \mbox{(i...i+3)C} &  $S_n$   &   \hfil ?
 &  \hfil ?  & \hfil ?    & \hfil ?  & \hfil ?  
 & \hfil ? & \hfil ?  & \hfil ?  
 &   \hfil ? & \hfil ? & \hfil ?  
 & \hfil ? & \hfil ?  \\ \hline
 \mbox{(i...i+3)CI} &  $S_n$   &  \hfil ?
 &  \hfil ?  & \hfil ?    & \hfil ?  &  \hfil ? 
 & \hfil ? &  \hfil ? &  \hfil ? 
 &  \hfil ?  & \hfil ? &  \hfil ?  
 & \hfil ?  &\hfil ?   \\ \hline
%
%
%
%
Pref.cyc  &   $S_n$ & \hfil $n-1$    
&    \hfil ? &   \hfil ?   
&  $n-1.7$\ding{117} & $n-1$\ding{117}  & 
&  $\to -1.25$\ding{117}  &   $\to 1.55$\ding{117} 
&  \hfil ? &\hfil ? 
&  \hfil ? &\hfil ?  &  \hfil ? \\ \hline
%
Pref.cyc+I  &  $S_n$ &  \mbox{ $n-1${\ding{117}}}   
&   \hfil ?  &  \hfil ?   
&   & $n-2$\ding{117}  & 
&  \hfil ? &  \hfil ?  
&  \hfil ? &  \hfil ? 
&   \hfil ? &  \hfil ?  &   \hfil ?  \\ \hline
%
Down.cyc   & $S_n$ &    \mbox{ $n-1${\ding{117}}}    &  \hfil ?  & Stir ling    
&  \mbox{$\approx 0.86n$\ding{117}}   &   \mbox{$\approx$Mean\ding{117}}   
&  \hfil ? &   \hfil ?  & \hfil ? 
&  $+$\ding{117}  &  \hfil ? &   \hfil ? 
&  \hfil ? &  \hfil ?  \\ \hline%
%
Down.cyc+I   & $S_n$ &    \mbox{ $n-1${\ding{117}}}    &  \hfil ?  & \hfil ?   
&  \mbox{$\approx 0.75n$\ding{117}}   &   \mbox{$\approx$Mean\ding{117}}   
&  \hfil ? &   \hfil ?  & \hfil ?  
& $+${\ding{117}} &  \hfil ? &   \hfil ? 
&  \hfil ? &  \hfil ?  \\ \hline
Inc.3cyc & $A_n$ &  \mbox{ $n-1|2${\ding{117}}}  &  \hfil ?  & ?    & 
\mbox{$\approx 0.5n$\ding{117}} &  \mbox{$\approx$Mean\ding{117}}    
&  \hfil ?   
&  \hfil ?  &  \hfil ? 
&   $+-$\ding{117}  & \hfil ?   & \hfil ?    & \hfil ?   
&  \hfil ?  \\ \hline
Inc.4cyc & $S_n$ &  \mbox{ $n-1${\ding{117}}}  &  \hfil ?  & ?    & 
\mbox{$\approx 0.5n$\ding{117}} &  \mbox{$\approx$Mean\ding{117}}    
&  \hfil ?   
&  \hfil ?  &  \hfil ?  
& \hfil ?  &\hfil ?  
& \hfil ?  &  \hfil ? &  \hfil ?  \\ \hline
%
%
%
RapaportM1 & $S_n$  & \mbox{$\left\lceil \tfrac{3n}{2}\right\rceil${\ding{117}}}   &  \hfil ? &   \hfil ?  
& \mbox{$\approx 1.4n$\ding{117}} & \mbox{$\approx$Mean\ding{117}} 
&   \hfil ?  &  \hfil ?   &   \hfil ?  &  \hfil ?  & \hfil ? & \hfil ? & \hfil ?  & \hfil ?   \\ \hline
%
RapaportM2 &  $S_n$ &  \mbox{  $\approx \tfrac{n^2+n}{2}$ {\ding{117}}}  &  \hfil ?  & ?    &   \mbox{$\approx \tfrac{n^2-3n}{2}$}{\ding{117}}  & \mbox{$\approx$Mean\ding{117}}  
&  \hfil ?  &  $\to -0.6$\ding{117} & $\to 0.5$\ding{117}  
& ? & \hfil ? &  \hfil ?  & \hfil ? 
&  \hfil ?  \\ \hline
%
Globes n/1& + &  \mbox{ $\left\lceil \tfrac{5n+12}{4}\right\rceil${\ding{117}}}  &  \hfil ?  & ?    
& $\approx n+1${\ding{117}}  &  $\approx n+1${\ding{117}}   
& \hfil ?  & \hfil ?    & \hfil ?  
& \hfil ? & \hfil ? 
& \hfil ?  & \hfil ?   & \hfil ?   \\ \hline
3Pancake$S_1$  & $S_n$  & \mbox{$\tfrac{3n(n+2)-48+I}{8}${\ding{117}}}   
&  \hfil ?  & \hfil ?   
&   \mbox{ $\tfrac{ }{ }${\ding{117}}}    &   & \hfil ? 
&  \hfil ? & \hfil ? 
&  & \hfil ? 
& \hfil ?  & \hfil ?  &  \hfil ? \\ \hline
3Pancake$S_2$  & $S_n$  & \mbox{ $ ${\ding{117}}}   
&  \hfil ?  & \hfil ?   
&  \mbox{ $\tfrac{ }{ }${\ding{117}}}   &   & \hfil ? 
&  \hfil ? & \hfil ? 
&  & \hfil ? 
& \hfil ?  & \hfil ?  &  \hfil ? \\ \hline
3Pancake$S_3$  & $S_n$  & \mbox{ $ ${\ding{117}}}   
&  \hfil ?  & \hfil ?   
&  \mbox{ $\tfrac{ }{ }${\ding{117}}}   &   & \hfil ? 
&  \hfil ? & \hfil ? 
&  & \hfil ? 
& \hfil ?  & \hfil ?  &  \hfil ? \\ \hline
3Pancake$S_4$  & $S_n$  & \mbox{ $ ${\ding{117}}}   
&  \hfil ?  & \hfil ?   
&  \mbox{ $\tfrac{ }{ }${\ding{117}}}   &   & \hfil ? 
&  \hfil ? & \hfil ? 
&  & \hfil ? 
& \hfil ?  & \hfil ?  &  \hfil ? \\ \hline
3Pancake$S_5$  & $S_n$  & \mbox{ $ ${\ding{117}}}   
&  \hfil ?  & \hfil ?   
&  \mbox{ $\tfrac{ }{ }${\ding{117}}}   &   & \hfil ? 
&  \hfil ? & \hfil ? 
&  & \hfil ? 
& \hfil ?  & \hfil ?  &  \hfil ? \\ \hline
3Pancake$S_6$  & $S_n$  & \mbox{$\tfrac{3n(n+2)-48}{8}${\ding{117}}}   
&  \hfil ?  & \hfil ?   
&  \mbox{ $\tfrac{ }{ }${\ding{117}}}    &   & \hfil ? 
&  \hfil ? & \hfil ? 
&  & \hfil ? 
& \hfil ?  & \hfil ?  &  \hfil ? \\ \hline
3Pancake$S_7$  & $S_n$  & \mbox{ $ ${\ding{117}}}   
&  \hfil ?  & \hfil ?   
&   \mbox{ $\tfrac{ }{ }${\ding{117}}}   &   & \hfil ? 
&  \hfil ? & \hfil ? 
&  & \hfil ? 
& \hfil ?  & \hfil ?  &  \hfil ? \\ \hline
%
Pancake& $S_n$  &  \mbox{ $\approx 1.2n$? }  &  \hfil ?  & ?    & $<\tfrac{17n}{12}$  & $\approx n/2$*  &  $\approx 0.2n$?  &  ? & ? 
& $+-$ & $NP/2?$ &  \hfil $-$   &  \hfil ?  & \hfil ?  \\ \hline
%
Reversals & $S_n$  &  $n-1$  &  \hfil ?  & ?    
&  $\approx n/2$ & $\approx Mean$?  & $\approx 0.05 n$ ? 
& \hfil ?    & \hfil ?  
& $+$ & $NP/1.5$ & \hfil $-$  
& \hfil ? & \hfil ?  \\ \hline
%
sReversals &  $B_n$  &  \mbox{ $n-1$  }  &  \hfil ?  & ?    
& $\approx n-\tfrac{n \log n}{2} $  
&  $\approx Mean$?   &   $\approx 0.05 n$  &  \hfil ?   & \hfil ?   
& $+-$ & $+$ & $+-$  
& \hfil ? & \hfil ?  \\ \hline
%
Transposons &  $S_n$ & $\left\lceil \tfrac{n+1}{2}\right\rceil$*  &  \hfil ?  & ?    
&  $\Theta(n)$ &    $\Theta(n)$   & &  \hfil ?  
&  \hfil ?  
& $-$ & $NP/1.375$ & \hfil $-$  
& \hfil ? & \hfil ?  \\ \hline
%
%
%
$(i,j)$ &  $S_n$  & $n-1$     &  \hfil ?  & Stir ling    
& $\approx n-\ln n$  & $\approx n-\ln n$  &  $\approx \ln n$  
&   $\approx \tfrac{1}{\sqrt {\ln n}} $  &  $\approx \tfrac{1}{ {\ln n}} $  
& $+$ & $+$ & $+$  & Int  &  $n \ln n$ \\ \hline
%
$(1,i)$(Star) & $S_n$  &$\left\lfloor \tfrac{3(n-1)}{2}\right\rfloor$  &  \hfil ?  & ?    
& $\approx n-\ln n$  & $\approx n-\ln n$   & $\approx \ln n$ &   &  
& $+$ & $+$ & $+$  
& Int  &  \hfil ? \\ \hline
G-star transp.& $S_n$ & $+$  &  \hfil ? &  \hfil ? &  \hfil ? &  \hfil ? &  \hfil ? &  \hfil ? &  \hfil ? &  \hfil ? &  \hfil ? & \hfil ?  &  \hfil ? & \hfil ?    \\ \hline
%

\end{longtable}
}


\clearpage

For the table for the Schreier coset graphs we use the same notation modulo two exceptions:
\begin{enumerate}
    \item "Coset" -- type of the coset,  currently we mainly work with "Bin", which is $S_{n}/(S_{\lfloor n/2\rfloor}\times S_{n-\lfloor n/2\rfloor})$ 
    ("Grassmannian over $F_1$"), 
    or, in simple language,  nodes are vectors with only  0,1 coordinates and  with $\lfloor n/2\rfloor$ zeros
    \item "God's number" -- instead of diameter we consider the largest distance to some selected node in the graph.
    We take vectors with zeros first, then units (i.e. sorted vectors)
    as our default choice for "Bin" for this initial state. 
\end{enumerate}
{\footnotesize 
\begin{longtable}{|@{}p{1.7cm}|p{0.4cm}|@{}p{1.7cm}|@{}p{1.1cm}
|@{}p{0.4cm}|@{}p{1cm}|p{1cm}|p{0.7cm}|p{0.7cm}|p{0.7cm}|p{0.9cm}|@{}p{0.8cm}
|p{0.5cm}|@{}p{0.7cm}|@{}p{1cm}|}

\caption{Summary of properties of Schreier graphs}\label{tab:long2}\\

\hline
\textbf{Gene-} & \textbf{Co-} & \textbf{God's} & \multicolumn{7}{|c|}{\textbf{Growth}} & \textbf{Anti-}  & \textbf{Algo-} & \textbf{Met-} & \textbf{Spec-} & \textbf{Mixing}  \\ 
\textbf{rators} & \textbf{set} & \textbf{number} & \textbf{PDF} & \textbf{F-la} &  \textbf{Mean}  & \textbf{Mode} & \textbf{Var} & \textbf{Skew} & \textbf{Kurt} &  \textbf{podes} & \textbf{rithm} &  \textbf{ric} & \textbf{trum}  & \textbf{Time} \\
\hline
\endfirsthead

\hline
\textbf{Gene-} & \textbf{Co-} & \textbf{God's} & \multicolumn{7}{|c|}{\textbf{Growth}} & \textbf{Anti-}  & \textbf{Algo-} & \textbf{Met-} & \textbf{Spec-} & \textbf{Mixing}  \\ 
\textbf{rators} & \textbf{set} & \textbf{number} & \textbf{PDF} & \textbf{F-la} &  \textbf{Mean}  & \textbf{Mode} & \textbf{Var} & \textbf{Skew} & \textbf{Kurt} &  \textbf{podes} & \textbf{rithm} &  \textbf{ric} & \textbf{trum}  & \textbf{Time} \\
\hline
\endhead

\hline
\multicolumn{15}{|r|}{\textit{Continued on next page}} \\
\hline
\endfoot

\hline
\multicolumn{15}{|r|}{\textit{End of table}} \\
\hline
\endlastfoot

%
LRX  &  Bin  
&  \mbox{ $ \tfrac{n(3n-4)+32-2(n\%4)}{16}$ {\ding{117}} }  
&   \hfil ?  & ?    
&   $\approx 0.16n^2$ &  $\approx 0.16n^2$ & \hfil ? 
&   $\approx -.5$\ding{117} &   $\to 3$\ding{117}  
&   $+-$\ding{117}  &  $+$\ding{117}   & \hfil ?   
& \hfil ? & \hfil ?   \\ \hline
Pancake   & Bin    
& $n-1|2$\ding{117}  
&  \hfil ?  &  $+$\ding{117}    
&  $\approx n/2$\ding{117} & $\lfloor n/2 \rfloor-I(n\%4)$\ding{117}  & \hfil ?   
& $\to 0$\ding{117} &  $\to 0$\ding{117}   
& $+$\ding{117}  & $+$\textcolor{blue}{\ding{117}} 
& \hfil ?   
&  \hfil ? & \hfil ?  \\ \hline
Transposons  &  Bin  
&  \mbox{ $\left\lfloor \tfrac{n}{2} \right\rfloor${\ding{117}}} 
&  Gauss\ding{117}  & $+$\ding{117}  
& $n/4$\ding{117}  &   $(n+1)//4$\ding{117} & \hfil ? 
&  $\to 0${\ding{117}}  &  $\to 0${\ding{117}}     
&  $+${\ding{117}} & \hfil ?  & \hfil ?   
&  \hfil ? & \hfil ?  \\ \hline
Reversals  &  Bin  
& \multicolumn{13}{c|}{ Growth coincides with transposons }  \\ \hline
RapaportM1  &  Bin  
&  $n-1$\ding{117}  
&  \hfil ?  & ?    
&  $\approx \tfrac{3n-4}{4} $\ding{117}  &  $\approx Mean$\ding{117}  
&  \hfil ?  
&   \hfil ?  &   \hfil ?    
&  $+-$\ding{117} & \hfil ?  & \hfil ?   
&  \hfil ? & \hfil ?  \\ \hline
RapaportM2  &  Bin  
&  $\tfrac{n(n+1)}{4}-3|4.75 $\ding{117}  
&  \hfil ?  & ?    
&  $\approx 0.21n $\ding{117}   &  $\approx $Mean\ding{117}    
& \hfil ? 
&  $\approx 0.6 $\ding{117}   &   $\approx 0 $\ding{117}    
&  $+$\ding{117}   & \hfil ?  & \hfil ?   
&  \hfil ? & \hfil ?  \\ \hline
 \mbox{(i,i+1,i+2)}  & Bin    
&   $1/8n^2+I(\%4)$\ding{117} 
&  \hfil ?  & ?    
&  $\approx 0.06n^2$\ding{117} & $\approx 0.06n^2$\ding{117}   
& \hfil ? 
&  $\to 0$\ding{117}  &    $\to 0$\ding{117}  
&  $+-$\ding{117} & \hfil ?  & \hfil ?   
&  \hfil ? & \hfil ?  \\ \hline
 \mbox{(i...i+3)}  & Bin    
&  $n^2/2+I(\%6)$\ding{117} 
&  \hfil ?  & ?    
&   $\approx 0.04n^2$\ding{117}  &  $\approx 0.04n^2$\ding{117}  
& \hfil ? 
&  $\to 0$\ding{117} & $\to 0$\ding{117}    
&    $+-$\ding{117}  & \hfil ?  & \hfil ?   
&  \hfil ? & \hfil ?  \\ \hline
 \mbox{(i...i+4)}  &  Bin   
&  $9/16n^2+I(\%8)$\ding{117} 
&  \hfil ?  & ?    
&  $\approx 0.03n^2$\ding{117}  & $\approx 0.03n^2$\ding{117}   
& \hfil ? 
&  $\to 0$\ding{117}   &  $\to 0$\ding{117}     
&  $+-$ \ding{117}  & \hfil ?  & \hfil ?   
&  \hfil ? & \hfil ?  \\ \hline
3Pancake$S_1$  & Bin& \mbox{$\tfrac{3n^2+36}{16}+I${\ding{117}}}   
&  \hfil ?  & \hfil ?   
&   \mbox{ $\tfrac{ }{ }${\ding{117}}}    &   & \hfil ? 
&  \hfil ? & \hfil ? 
&  & \hfil ? 
& \hfil ?  & \hfil ?  &  \hfil ? \\ \hline
3Pancake$S_2$  & Bin & \mbox{ $\approx \tfrac{n(n-2)}{4}  ${\ding{117}}}   
&  \hfil ?  & \hfil ?   
&  \mbox{ $\tfrac{ }{ }${\ding{117}}}   &   & \hfil ? 
&  \hfil ? & \hfil ? 
&  & \hfil ? 
& \hfil ?  & \hfil ?  &  \hfil ? \\ \hline
3Pancake$S_3$  & Bin  & \mbox{ $\tfrac{n(n+14)-I}{8}${\ding{117}}}   
&  \hfil ?  & \hfil ?   
&  \mbox{ $\tfrac{ }{ }${\ding{117}}}   &   & \hfil ? 
&  \hfil ? & \hfil ? 
&  & \hfil ? 
& \hfil ?  & \hfil ?  &  \hfil ? \\ \hline
3Pancake$S_4$  & Bin  & \mbox{ $\approx \tfrac{5n^2}{32}${\ding{117}}}  
&  \hfil ?  & \hfil ?   
&  \mbox{ $\tfrac{ }{ }${\ding{117}}}   &   & \hfil ? 
&  \hfil ? & \hfil ? 
&  & \hfil ? 
& \hfil ?  & \hfil ?  &  \hfil ? \\ \hline
3Pancake$S_5$  & Bin  &  \mbox{ $\approx \tfrac{47n^2}{128}${\ding{117}}}  
&  \hfil ?  & \hfil ?   
&  \mbox{ $\tfrac{ }{ }${\ding{117}}}   &   & \hfil ? 
&  \hfil ? & \hfil ? 
&  & \hfil ? 
& \hfil ?  & \hfil ?  &  \hfil ? \\ \hline
3Pancake$S_6$  & Bin  & \mbox{ $\tfrac{3n^2}{16}+I${\ding{117}}}   
&  \hfil ?  & \hfil ?   
&  \mbox{ $\tfrac{ }{ }${\ding{117}}}    &   & \hfil ? 
&  \hfil ? & \hfil ? 
&  & \hfil ? 
& \hfil ?  & \hfil ?  &  \hfil ? \\ \hline
3Pancake$S_7$  & Bin & \mbox{ $\tfrac{5n^2+58n - 64 +I}{72}${\ding{117}}}   
&  \hfil ?  & \hfil ?   
&   \mbox{ $\tfrac{ }{ }${\ding{117}}}   &   & \hfil ? 
&  \hfil ? & \hfil ? 
&  & \hfil ? 
& \hfil ?  & \hfil ?  &  \hfil ? \\ \hline
%

%

\end{longtable}
}

\subsection{Largest diameters found (and known) for  \texorpdfstring{$n\leq 15$}{n<=15}  }

We conducted extensive search for generators producing large diameters for small $n$, and
remarkable patterns showed up. These will be discussed in the next section,
here we just present the diameters and the corresponding generators,
and organization of the experiments. But already here it is worth to highlight that all generators we found
are related to involutions. 

To the best of our knowledge these are the largest diameters
known so far. We consider both standard  undirected  and directed  cases (meaning that the set of
generators is not necessarily closed under inverses). For  the latter case the same diameters up to $n\le 7$
were  found in \cite{nagy2016computational} (table 4), but our results extend to  larger $n$. 

\begin{Conj} The diameters provided below are the largest possible or at least within say 5\% from them.
\end{Conj}
It is impossible to make an exhaustive search even for these values of $n$,
so we cannot exclude the chance that large diameters exist,
though it seems unlikely to us that they will be significantly larger
than presented here. In any case we hope our results may stimulate further research.  

The last column of the tables represent a pattern of generators -- see section \ref{sect:new-gens}. It is surprising that  most of the generators we found do not look random, but follow beautiful patterns described below.

\begin{center}
\begin{tabular}{|c|p{1.5cm}|p{6.5cm}|p{5.5cm} |}
\hline
\multicolumn{4}{|c|}{Maximal diameter for Cayley graph of the group $S_n$ (undirected graph)} \\
\hline
n & Maximal diameter & Example of a set of generators & Pattern \\
\hline
$3$ & $3$ & $(12)$, $(02)$  & Custom \\
\hline
$4$ & $6$ & $(02)$, $(031)$, $(013)$  & Koltsov3 type=1, d=1, k=0 \\
\hline
$5$ & $10$ & $(23)$, $(04321)$, $(01234)$ & Koltsov3 type=1, d=1, k=1 \\
\hline
$6$ & $16$ & $(01)(23)(45)$, $(45)$, $(12)(34)$  & Koltsov3 type=1, d=1, k=0  \\
\hline
$7$ & $30$ & $(23)(56)$, $(14)(26)(35)$, $(06)(23)(45)$  & Koltsov3 type=2, k=1 \\
\hline
$8$ & $39$ & $(45)(67)$, $(23)(46)(57)$, $(07)(13)(26)$  & Custom  \\
\hline
$9$ & $52$ & $(08)(15)(37)$, $(012)(354)(687)$, $(021)(345)(678)$  & Lytkin's triangles  \\
\hline
$10$ & $77$ & $(01)(23)(45)(68)(79)$, $(45)(67)(89)$,  $(16)(24)(38)$  & Koltsov3 (experimental), k=2 \\
\hline
$11$ & $85$ & $(01)(23)(45)(67)(89)$, $(13)(24)(56)(78)(9,10)$, $(14)(23)$ & Koltsov3 type=2, k=1 \\
\hline
$12$ & $95$ & $(01)(23)(45)(67)(89)(10,11)$, $(12)(35)(46)(78)(9,10)$, $(36)(45)$  & Koltsov3 type=2, k=3 \\ 
\hline
$13$ & $111$ & $(01)(35)(67)(89)(10,11)$, $(1234)(56)(78)(9,10)(11,12)$, $(1432)(56)(78)(9,10)(11,12)$ 
& Sheveleva2, k=2 (inverse closure) \\ 
\hline
$14$ & $132$ & $(01)(23)(45)(67)(89)(10,11)(12,13)$, $(12)(34)(56)(78)(9,10)(11,12)$, $(45)$ & Koltsov3 type=1, d=1, k=4 \\ 
\hline
$15$ & $148$ & $(01)(23)(45)(67)(89)(10,11)(12,13)$, $(12)(34)(56)(78)(9,10)(11,12)(13,14)$, $(58)(67)$ & Koltsov3 type=2, k=5 \\ 
\hline
\end{tabular}
\end{center}
\noindent
\\


\begin{center}
\begin{tabular}{|c|c|c|c|}
\hline
\multicolumn{4}{|c|}{Maximal diameter for Cayley graph of the group $S_n$ (oriented graph)} \\
\hline
n & Maximal diameter & Example of a set of generators & Sheveleva2 params (n,k)\\
\hline
$3$ & $3$ & $(01)$, $(02)$ &\\
\hline
$4$ & $7$ & $(01)$, $(123)$ &\\
\hline
$5$ & $14$ & $(01)(23)$, $(0314)$ &\\
\hline
$6$ & $18$ & $(01)(23)(45)$, $(012)(34)$ &\\
\hline
$7$ & $34$ & $(01)(23)(45)$, $(052)(146)$ &\\
\hline
$8$ & $44$ & $(01)(23)(45)$, $(1736)(25)$ & $(8,2)$\\
\hline
$9$ & $61$ & $(01)(23)(45)$, $(3647)(12)(58)$ & $(9,3)$\\
\hline
$10$ & $83$ & $(01)(23)(45)(67)(89)$, $(185)(237)(469)$ &\\ 
\hline
$11$ & $93$ & $(01)(23)(45)(67)(89)$, $(1528)(47)(6,10)$ & type 2\\ 
\hline
$12$ & $106$ & $(01)(23)(45)(67)(89)$, $(1,11,9,10)(04)(28)(36)$ & $(12,4)$\\ 
\hline
$13$ & $147$ & $(01)(23)(45)(67)(89)$, $(1,11,2,12)(34)(56)(78)(9,10)$ & $(13,2)$\\ 
\hline
\end{tabular}
\end{center}

For the directed case (i.e. not inverse closed generators)
the search has been organized as follows. For small $n=3,4$
we considered several possible numbers of generators -- from $2$ to $5$, and it was observed that
the largest diameter is observed for $2$ generators, which is not surprising
since having fewer generators means fewer words of a given length, and hence potentially larger diameters.
We continue the search with 2 generators for $n$ up to $13$. 
We perform an exhaustive search up to $n\le 9$, and use randomized search afterwards.
We did not search all possible pairs, but relied on the simple fact that conjugating all generators
produces an isomorphic graph, so say the first generator can be taken as a unique representative of a conjugacy class,
with a loop over conjugacy classes. This of course significantly reduces the search space. For the randomized search we sampled the second generator from each conjugacy class uniformly. 
To reduce the search further for $n\ge 12$, we mostly considered conjugacy classes for the first generator
to be involutions, and used guesses from previously observed patterns. 

Notebooks: \href{https://www.kaggle.com/code/ani3221/max-diam-classes-transpositions}{1},
\href{https://www.kaggle.com/code/fedmug/max-diameter}{2}.


\clearpage
\subsection{Distribution of diameter over conjugacy classes pairs -- involutions strikes }

\begin{figure}[t]
  \centering
  \includegraphics[width=1.2\linewidth]{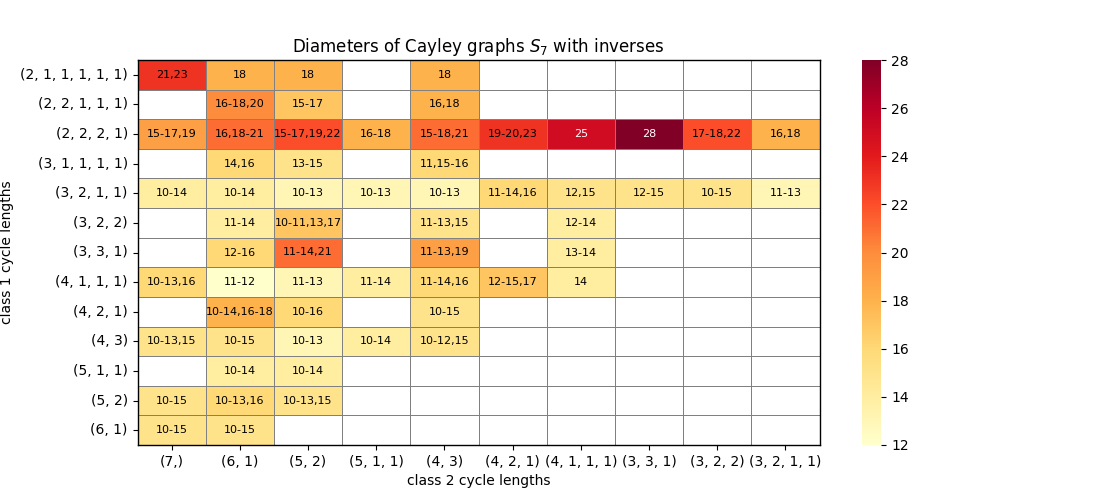}
  \includegraphics[width=1.2\linewidth]{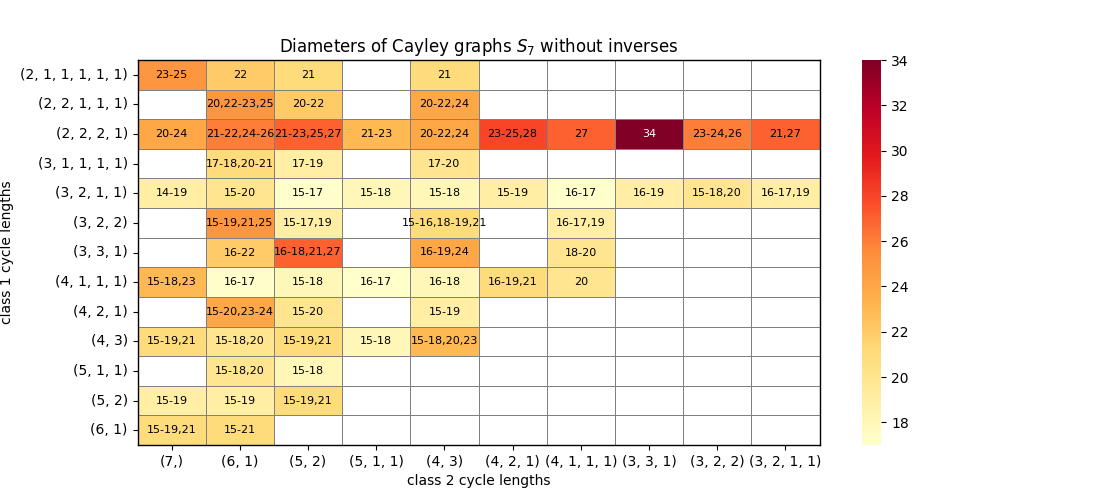}
  \caption{Diameters of all possible Cayley graphs for $S_7$ generated by two permutations with/without their inverses.}
  \label{fig:diameters_7}
\end{figure}

Figure \ref{fig:diameters_7} represents the dependence of diameters on conjugacy classes of generators, 
obtained by exhaustive search over all combinations of two generators, with the generators taken from one of the listed classes. The presented figure is for $S_7$.  The values in each cell represent all values of the diameters found for the corresponding pair of classes and  the heatmap coloring shows the maximal diameter found. We leave the cells blank if  $S_n$ is not generated by any pair of elements from the above classes. 
Due to symmetry we are showing only
the relevant part of data, i.e. we do not show duplicate entries for generator choices $(g_1,g_2)$ and $(g_2,g_1)$.

One can see that large diameters appear when one of the classes is an involution and the other one has rather short cycles in its decomposition.
This is natural to expect, because if longer cycles are present it means that individual generators would have higher order and allow one to create
more words like $XX...XX$. The more words of fixed length $k$ one has, the less the diameter can be, since possible
words will exhaust the finite space of group elements earlier.
So naively it is easy to expect that generators related to as small degree as possible are good candidates. 
Exhaustive search confirms it for many cases. So we propose:

\begin{Conj} The generating set with largest diameter for $S_n$  contains an involution (at least for infinite number of $n$),
 both for directed and undirected cases.
\end{Conj}

We performed a similar exhaustive search up to $n\le 9$ and a randomized search up to $n=12$ --
for the pairs of generators and for  both undirected and directed cases. The results are consistent with the conjecture.
Similar heatmap plots for other $n$ can be found in the full version of the manuscript,  or up to $n\le 9$ in the \href{https://www.kaggle.com/code/fedmug/diameter-visualizations}{notebook}.

\clearpage

\subsection{Graphical representation for any generators -- pattern "square with whiskers" }

Here we describe a very simple  visual representation of elements (e.g.  generators) of permutation groups,
and present a pattern "square with whiskers" which corresponds to most of the largest diameters found for $n\le 15$.

{\bf Step 1. Single permutation.} Each permutation defines a directed graph on $n$ nodes in a natural and obvious way -- 
if $p(i)=j$ let us connect $i\to j$. (In  other words, we take a permutation matrix
and consider it as adjacency matrix of a directed graph). 
Clearly, if the permutation is an involution, then orientation of edges is unnecessary -- if $i\to j$, then $j\to i$
(in matrix language -- the permutation matrix is symmetric). 

{\bf Step 2. Many permutations -- use colors. } Consider several permutations and just use the same
construction but with different colors to represent edges coming from different permutations.

Thus for any sequence of permutations (generators) we constructed a directed multi-colored graph on $n$ nodes. 

The code for the visualization can be found e.g. in this 
\href{ https://www.kaggle.com/code/olegpushs/draw-edges }{notebook}.

\begin{figure}[t]
  \centering
  \includegraphics[width=0.5\linewidth]{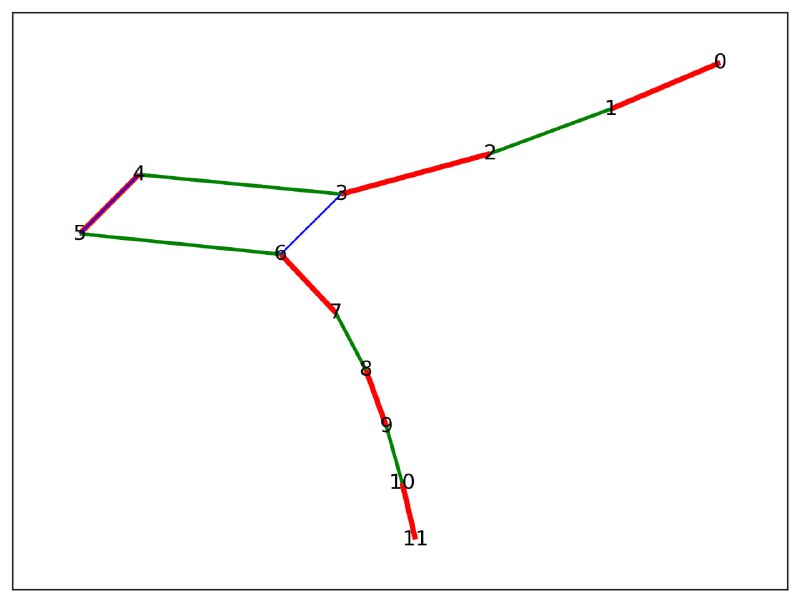}
  \caption{Three generators of $S_{12}$ represented graphically, showing the pattern "square with whiskers". 
  Three involutions, so the edges are undirected.}
  \label{fig:SquareWhiskers12.jpg}
\end{figure}

The figure \ref{fig:SquareWhiskers12.jpg} presents an example of this visualization
and also presents an example of the pattern which we call "square with whiskers" -- there is one 4-cycle (square)
and two branches going out of its corners. 

\begin{Def} We will say that the generators follow a "square with whiskers" pattern if the underlying undirected  graph (forgetting colors and multi-edges) is of that type -- one 4-cycle (square) and two branches.
\end{Def}

\begin{Conj} The generators with maximal (or nearly) diameter for $S_n$ and $A_n$  follow the "square with whiskers" pattern  (at least for infinite number of $n$). We expect that to be true for both undirected and directed Cayley graphs.
\end{Conj}

The computations up to $n\le 15$ described above support this conjecture.

\clearpage

\subsection{New generators with large diameter: `Sheveleva2` and `Koltsov3involutions` (directed and undirected cases) }\label{sect:new-gens}

Here we present new families of generators "Sheveleva2" (directed case) and "Koltsov3involutions" (standard undirected case),
which were found by generalizing the patterns seen on the large diameters obtained for $n\le 15$,
by the participants of the CayleyPy project Anastasia Sheveleva and Ivan Koltsov. We expect these 
generators to produce the largest (or near) diameters for infinite number of $n$. 

\begin{figure}[ht]
    \centering
    \begin{minipage}{0.45\textwidth}
        \centering
        \includegraphics[width=\linewidth]{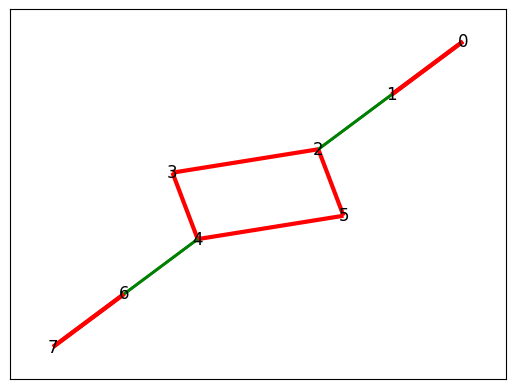}
        \caption{Sheveleva2 $n = 8, k = 3$}
        \label{fig:n = 8, k = 3}
    \end{minipage}
    \hfill
    \begin{minipage}{0.45\textwidth}
        \centering
        \includegraphics[width=\linewidth]{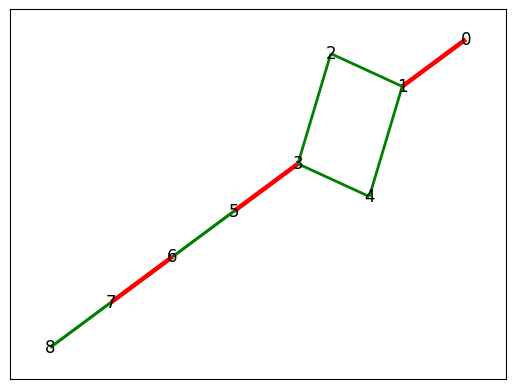}
        \caption{Sheveleva2 $n = 9, k = 2$}
        \label{fig:n = 9, k = 2}
    \end{minipage}
\end{figure}

\begin{figure}[ht]
    \centering
    \begin{minipage}{0.45\textwidth}
        \centering
        \includegraphics[width=\linewidth]{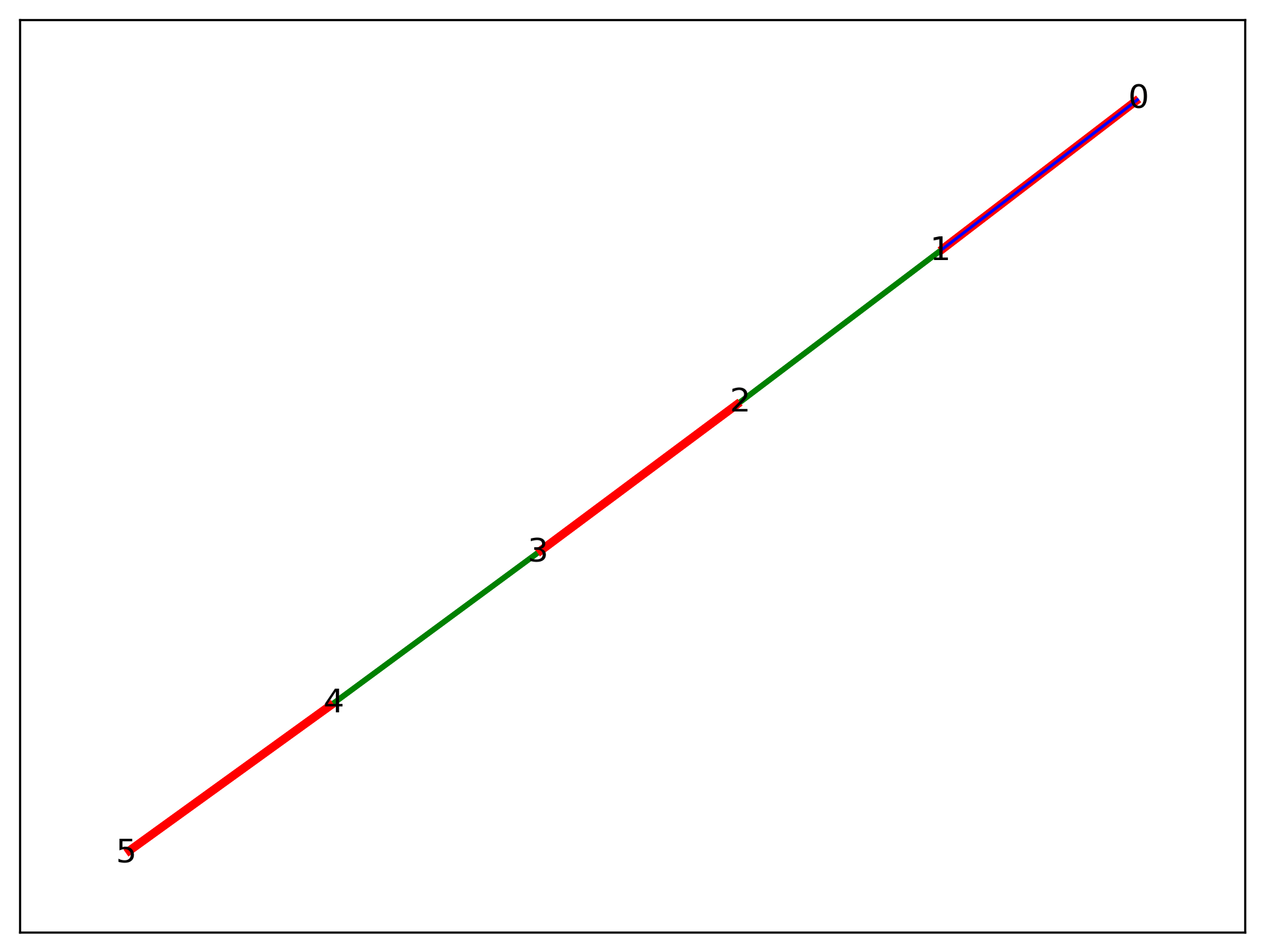}
        \caption{Koltsov3involutions type 1 $n = 6, k = 0, d=1$}
        \label{fig:Koltsov1}
    \end{minipage}
    \hfill
    \begin{minipage}{0.45\textwidth}
        \centering
        \includegraphics[width=\linewidth]{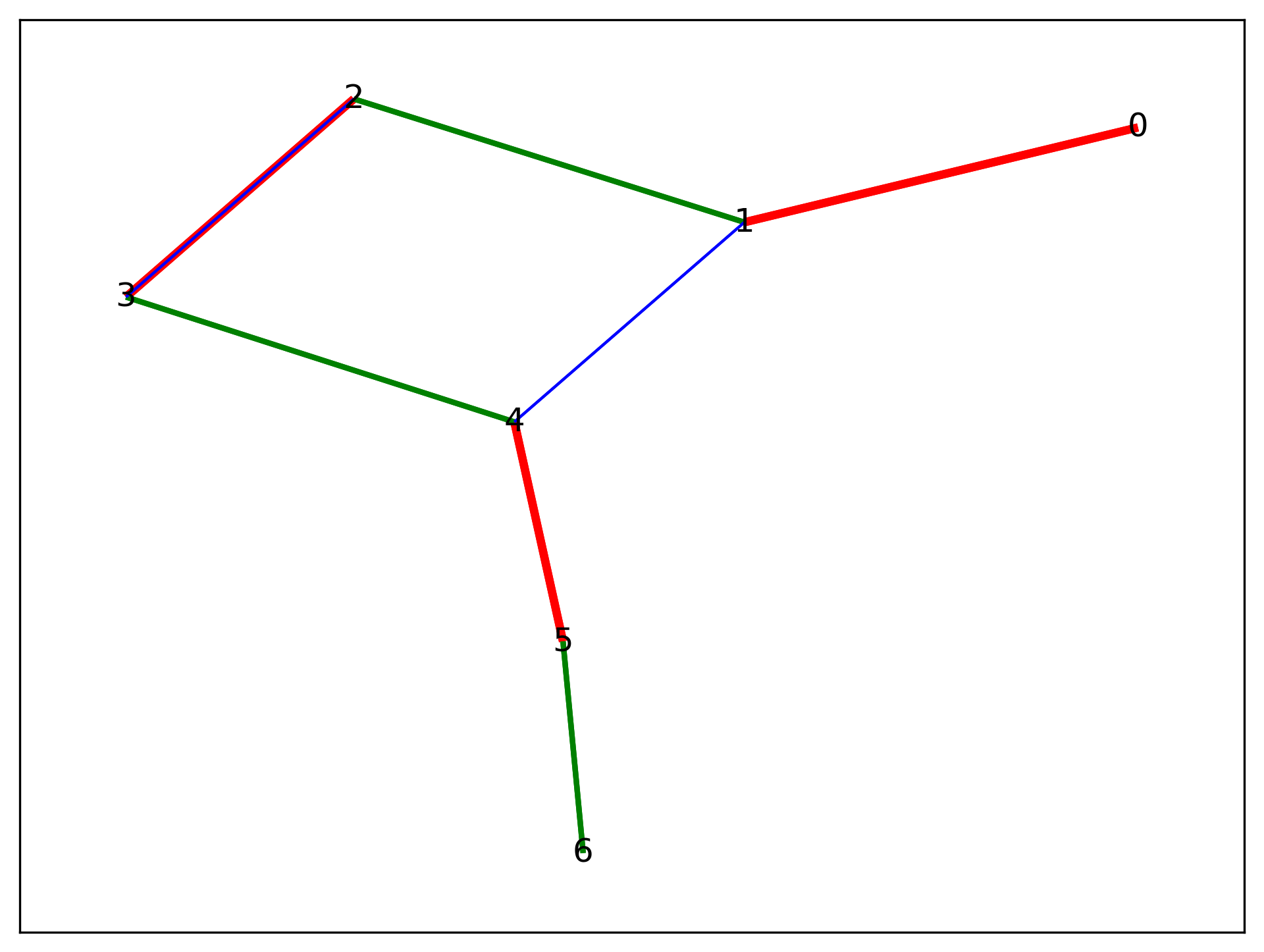}
        \caption{Koltsov3involutions type 2 $n = 7,k=2$}
        \label{fig:Koltsov2}
    \end{minipage}
\end{figure}

"Sheveleva2" generators are a particular case of the "square with whiskers" pattern described in the previous section.
Figures \ref{fig:n = 8, k = 3}, \ref{fig:n = 9, k = 2} show examples of these patterns. There are only two generators, corresponding to two colors in the figures.  One is an involution, the other one contains a single 4-cycle ("square"), multiplied by transpositions. We give formal description below, but it might be easier  to read off the pattern from the picture.

In the definition we work with permutations of $(0,...,n-1)$, not as usually in mathematics $(1,...,n)$.

\begin{Def} The set of "Sheveleva2" generators consists only of 2 elements, and they have two parameters $n,k$, where $n$ is the length of permutations and $k$ corresponds to the "position of the square". 
One generator is always an involution (generator "A"), the other one is a product of transpositions and of one 4-cycle 
(generator "S"  -- contains a square in the graphical representation). 
They are defined by the following simple construction. One starts with the transposition (0,1) which 
will be assigned to the first generator, then the transposition (1,2) is assigned to the second one, then transposition (2,3) to the first one,
then (3,4) to the second one, and so on sending transpositions $(i,i+1)$ to one of the two generators. This is repeated
until one reaches the index $k$ (position of "square"), here one creates the 4-cycle $(k-1,k,k+1,k+2)$ and assigns it to 
the current generator, then one continues what was done at the beginning: creating transpositions $(i,i+1)$ and assigning them to the first or second generator
alternately. When we reach $n$ the process stops and generators are obtained by multiplying all permutations assigned to them. 
\end{Def}

The set of generators "Koltsov3involutions" consists of 3 involutions. 
We give several variations of the construction differing by the parameter "type". 
But the backbone of the construction is similar and simple in all cases:
one starts with two involutions both of which are products of $(i,i+1)$, but one is for all even $i$, 
whereas the other one is for odd $i$. We will denote them as "I"(interleaver)$=\prod_{i=0,2,4,...}(i,i+1)$ and "K"(conjugate to interleaver)$=\prod_{i=1,3,5,...}(i,i+1)$. 
Now we need to define the third involution "S" (swap) which is a short product of transpositions
that differs depending on the parameter "type". For type 1, "S" is just 1 transposition,
for type 2 it is a product of 2 transpositions.

\begin{Def} The set  "Koltsov3involutions type 1" depends on two additional parameters, $k,d$ (in addition to $n$, the length of the permutation) and
is defined simply by involutions "I","K", and "S", with the latter being just the transposition $(k,k+d)$.
\end{Def}

The set "Koltsov3involutions type 2" follows the "square with whiskers" pattern described above. 
\begin{Def} The set "Koltsov3involutions type 2" depends on a single additional parameter $k$ (in addition to $n$, the length of the permutation) and is
defined simply by involutions "I","K" described above and by taking "S"=$(k,k+3)(k+1,k+2)$. 
\end{Def}
The addition of the "S" generator to "I","K" creates a "square" in the graphical representation above.

Figures \ref{fig:Koltsov1}, \ref{fig:Koltsov2} present visualisations of these generators. 

\begin{Conj} Generators "Koltsov3involutions" and "Sheveleva2" provide maximal (or near) diameters
for standard undirected and directed cases respectively (at least for infinite number of $n$). 
\end{Conj}
The computations up to $n\le 15$ described above support this conjecture.

\clearpage
\subsection{Refinement of the Babai-like conjecture for \texorpdfstring{$S_n$}{Sn} }

The Babai conjecture is a fundamental conjecture on Cayley graphs of finite simple groups.
It states that for any simple group and any choices of the generators, the diameters
are bounded by $log^c(|G|)$, for some universal constant $c$.  
It is still open, but there was a substantial progress  e.g. by T.~Tao and his collaborators:
\cite{BGT2011,BGT2012, Breuillard2015}, and more recent works such as \cite{Eberhard2020,Eberhard2023,Eberhard2021}. (S.~Eberhand's  informal blog-posts  
\href{https://seaneberhard.com/2023/08/15/babais-conjecture-for-generating-sets-containing-transvections/}{2023},
\href{https://seaneberhard.com/2020/11/23/talk-at-tau-about-babais-conjecture-in-high-rank/}{2020Talk},
\href{https://seaneberhard.com/2020/05/21/babais-conjecture-for-at-least-three-random-generators/}{3random} provide a non-technical introduction). 

For $S_n$ and $A_n$, the variant of the conjecture simply predicts that the diameter is less than $O(n^2)$ or, in stronger variants, than $n^2$. 

We propose the following refinement of that  conjecture:

\begin{Conj} For any choice of the generators 
for $S_n$ and $A_n$ the diameter of 
the Cayley graph is bounded by 
 $n^2/2 + 4n$ (in the standard case when the generating set includes inverses).
\end{Conj}
We expect variations of the conjecture to be true, but currently we are less sure about the constants in these cases:
\begin{enumerate}
    \item $Cn^2/4 + O(n)$ for not necessarily inverse closed generators, with some $C:3/4\le C \le 1$
    \item $Cn^2 + O(n) $ for Schreier coset  graph for $S_{n}/(S_{\lfloor n/2 \rfloor}\times S_{n-\lfloor n/2 \rfloor}) $ with inverse closed generators, with some $C:1/4\le C \le 1/2$, most likely surprisingly $C=1/2$
    \item $Cn^2 + O(n) $ for "circular permutations" (\cite{alon2025circular}, with some $C:1/4\le C \le 1/2$
\end{enumerate}

The previously stated versions of the conjecture contained bounds $O(n^2)$ or just $n^2$. 
Our refinement consists in 1) reducing the quadratic coefficient and suggesting 
that it depends on whether the generators are inverse-closed or not, and on the coset structure (determining how exactly it 
it depends on the coset is an interesting question); 2) we optimistically propose that the subleading term is linear
and in particular not more than  $4n$ in the standard case. 

The conjecture that diameters of $S_n$ and  $A_n$ are $O(n^2)$ remains widely open; in fact, even a polynomial bound is not yet proven.   
Currently the best rigorous results are due to H.A.Helfgott and S. Seress \cite{helfgott2014diameter} 
(see also Helfgott's surveys: ~\cite{helfgott2014diameter}, \cite{helfgott2019growth}, \cite{helfgott2015random} ) and read:
\[
\operatorname{diam}(\mathrm{Sym}(n)) \;=\; 
\exp\!\bigl(O((\log n)^4 \log\log n)\bigr).
\]

{\bf Heuristics.} (Naive, but instructive (and fun).)
The general Babai conjecture
states that, for any simple group and any choice of the generators, the diameters
are bounded by $\log^c(|G|)$, for some universal constant $c$ (which is somewhat technical to our story and it may be simpler 
to think of it as $c=1$ for the moment). Typical examples of finite simple groups come from matrices over finite fields ${\mathbb Z}/p{\mathbb Z}$. The
number of all matrices over ${\mathbb Z}/p{\mathbb Z}$ is obviously $p^{n^2}$, and taking the logarithm we get $n^2 \log(p)$. Thinking of $S_n$ as a kind of  "field with one element" version of matrices  (so canceling $\log(p)$ by taking $p=1+\epsilon$ and leading term in $\epsilon$) we get exactly $n^2$ --  as in the standard form of
the conjecture. Now let us take a small detail into account: we are dealing not with all matrices but only the invertible ones,
so the counting is: $\prod_i (p^{n}-p^i)=\prod_i ((1+\epsilon)^{n}-(1+\epsilon)^i)$, and by the same argument we get $n(n+1)/2$.
This, surprisingly, corresponds to some examples with a quite large diameter, it might be even the largest.
So if one believes in both Babai conjecture and $O(n^2)$-conjecture (both widely regarded to be true),
and in addition that they should be naturally compatible, it would be difficult to believe that the coefficient in front of
$n^2$ is $1$. Instead, $1/2$ is a natural choice. 

{\bf Known results.} To the best of our knowledge there are no examples which contradict the $n^2/2 + 4n$ limit bound on the diameter that we conjecture.
Moreover,  we performed extensive computational experiments summarized in tables \ref{tab:long1}, \ref{tab:long2}, with more than half a hundred generators, and again they confirm our proposals. 
We are also not aware of any proven examples  with diameters  larger than  $ n(n+1)/2$ for infinite number of $n$.
At the same time, our computational experiments suggest that larger diameters may be possible, and our current fit is $n^2/2 + 4n$ for reasons discussed  below. 

\begin{figure}[ht]
  \centering
  \includegraphics[width=0.9\linewidth]{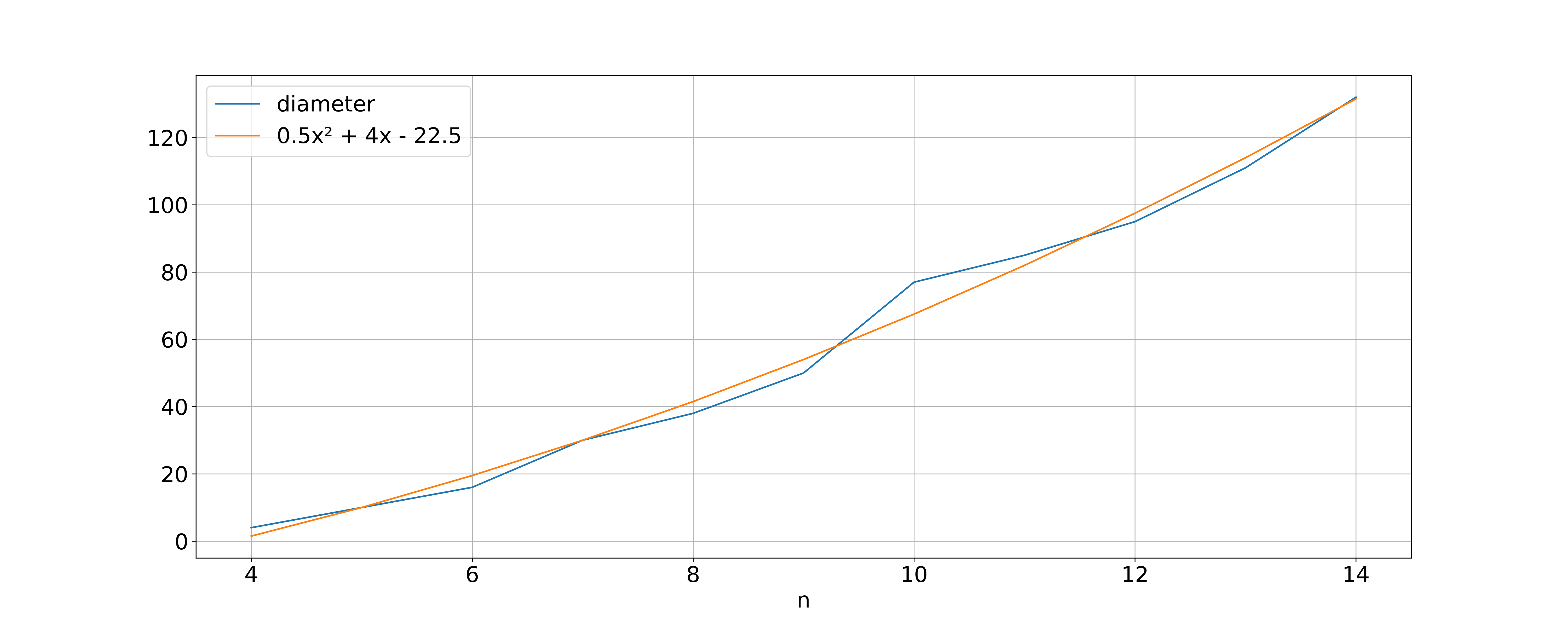}
  \caption{Maximal diameters found up to $n\le 14$ and quadratic fit: $ n^2/2 +4n - 22.5$ .}
  \label{fig:max_diam}
\end{figure}

{\bf Numerical argument.} 
As discussed in previous sections, we have found maximal known
diameters up to $n\le 14$ and the generators providing them followed nice patterns, which is hardly accidental. Now making a quadratic fit and some small corrections  we get a polynomial $n^2/2 +4n - 22.5$ (see figure \ref{fig:max_diam}). 
So we see a {\bf perfect match} -- the fit obtained from data $n\le 14$, and our knowledge on various examples of generators 
with known/conjectured diameters for all $n$, both point to $n^2/2 + Cn$, with $C$ from $1/2$ to $4$, so we safely take
maximum $C=4$ (apparently it is some overestimation). 
It is worth noting that obtaining data till $n=14$ is important, as making a fit on less data (e.g. below $n=9$) would give different results, overestimating the leading coefficient. 

{\bf Sanity check} -- our conjecture predicts a bound of 150 for $n=15$ 
(it was obtained without $n=15$, so this case can be used for validation), 
but we already demonstrated that 
our generators can produce a diameter of 148 , which is very close to the prediction,
thus confirming the expectations. 

Unfortunately, for examples of directed Cayley graphs the situation is less clear -- we know an example with $3n^2/4$,
but analysis for $n\le 14$ suggests large values around $n^2$, so there is no perfect match as above. But analysis up to $n\le 14$ 
can overestimate the result -- that is seen if we take values like $n\le 10$. We hope to resolve these issues in the future.

\begin{figure}[ht]
    \centering
    \includegraphics[width=0.95\linewidth]{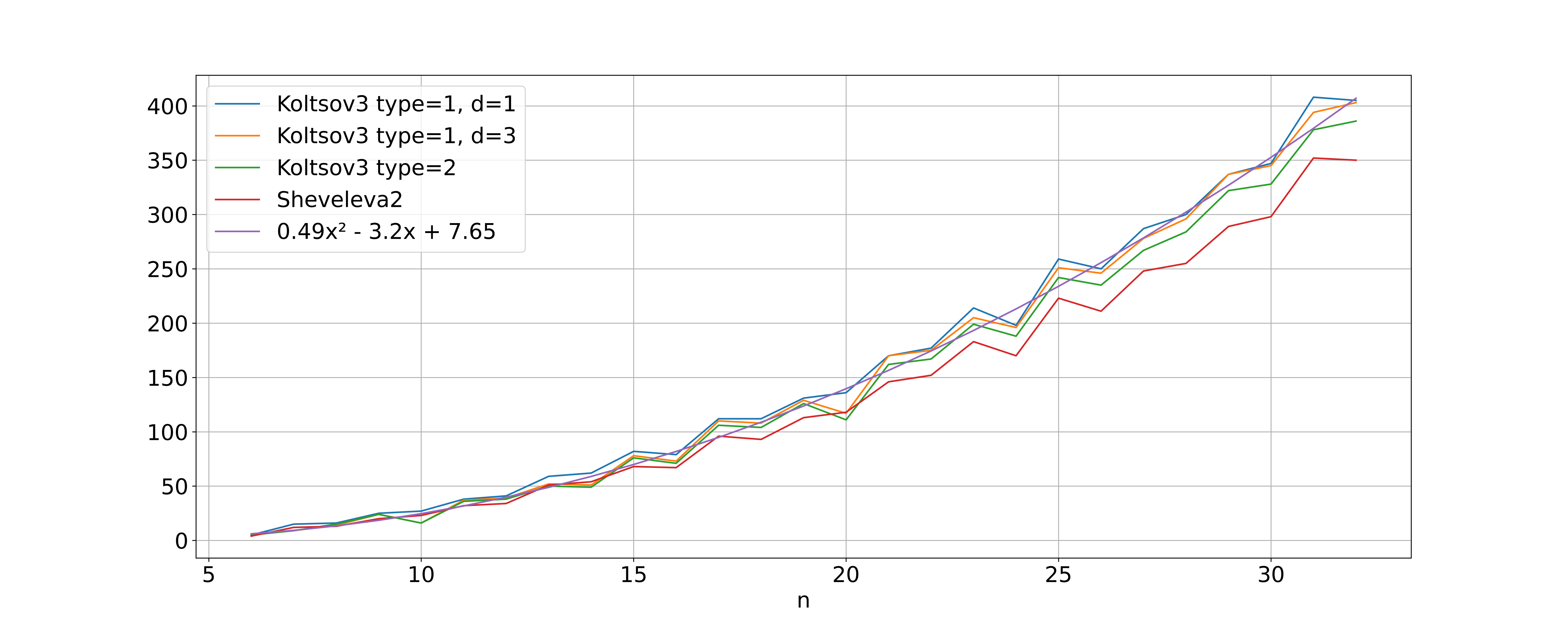}
    \caption{God's numbers of Schreier coset graphs $S_n/(S_{\lfloor n/2 \rfloor}\times S_{n-\lfloor n/2 \rfloor})$ (binary strings). }
    \label{fig:box_w_legs_cosets}
\end{figure}

{\bf Schreier coset graphs.} Figure \ref{fig:box_w_legs_cosets} provides similar example of God's number 
for the Schreier coset graphs $S_n/(S_{\lfloor n/2 \rfloor}\times S_{n-\lfloor n/2 \rfloor})$ with generators introduced 
in the previous section (maximal result over parameter $k$). In simple words, these cosets can be described as constructed by applications of $S_n$ to a vector with $\lfloor n/2 \rfloor$ zeros and $n-\lfloor n/2 \rfloor$ ones (taking a sorted vector as the initial state).
Since these graphs are smaller, we can compute up to $n\le 32$. 
Numerical fit suggests the bound $n^2/2 - 3.2n$, which is surprising because these graphs are noticeably smaller
than full Cayley graphs and it would be natural to see a decrease in the leading coefficient, moreover such a decrease is indeed 
observed in many examples typically from $n^2/2+O(n)$ to $n^2/4+O(n)$.
But here numerics suggest that leading coefficient is the same as for the full Cayley graph: $n^2/2+O(n)$, only the linear term changes. Currently we lack a theoretical proof that such a result is indeed achievable for all $n$.  We hope to get more insights on this soon and clarify the issue.

If this were true, it would suggest an interesting perspective to attack Babai-like conjectures
via cosets like  $S_n/(S_{\lfloor n/2 \rfloor}\times S_{n-\lfloor n/2 \rfloor})$ which are much smaller and typically more tractable
than full Cayley graphs. And it might be that diameters in the full case are not much larger (only linear difference, 
not quadratic).

\clearpage
\subsection{Bounds on the diameters of unitriangular groups}
\begin{figure}[ht]
    \centering
    \includegraphics[width=0.75\linewidth]{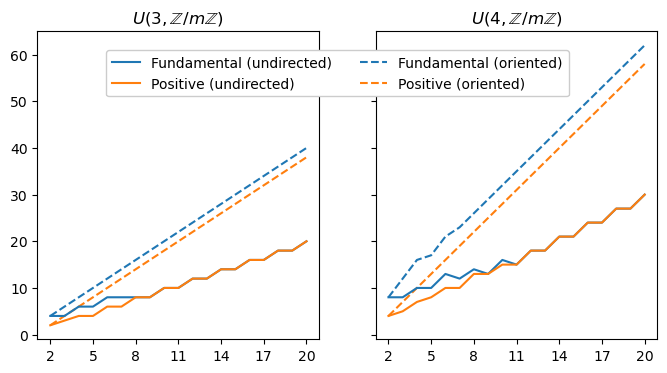}
    \caption{Diameters for unitriangular groups, depending on the modulus.}
    \label{fig:unitriangular_diameters}
\end{figure}

J.~Ellenberg in \cite{ellenberg1993diameter}, and later jointly with J.~Tymoczko  \cite{ellenberg2005unipotent}, established that the diameters of $U(n,\mathbb{Z}/p\mathbb{Z})$, the unitriangular groups over prime fields, are within a constant factor of $np+n^2\log(p)$ (see also enjoyable introductions in J.~Ellenberg's blog: \href{https://quomodocumque.wordpress.com/2010/07/03/in-which-i-publish-my-first-paper/}{1},
\href{https://quomodocumque.wordpress.com/tag/diameter/}{2}).
Here the group is considered with respect to the generators $E\pm E_{i,i+1}$ (fundamental roots). Experimentally, we see that starting from large enough $m$, the diameters of $U(n,\mathbb{Z}/m\mathbb{Z})$ are in fact linear in $m$, see Fig.~\ref{fig:unitriangular_diameters}. In particular, for $n=3$ and $m=8,\ldots,50$ the diameter equals $2\lfloor\frac{m}{2}\rfloor$, while for $n=4$ and $m=12,\ldots,31$ it is $3\lfloor\frac{m}{2}\rfloor$.

Another natural choice of generators is elementary transvections corresponding to all positive roots, i.e. $E\pm E_{ij}$ for $i<j$. One can also consider the directed Cayley graphs corresponding to monoid generators of the form $E+E_{i,i+1}$ or $E+E_{ij}$, respectively, thus providing four basic options: fundamental/positive and undirected/oriented.
\begin{Conj}
For large enough $m$ (depending on $n$) the diameter of $U(n, \mathbb{Z}/m\mathbb{Z})$ is a follows:
 \begin{itemize}
    \item Fundamental, undirected: $(n-1)\lfloor\frac{m}{2}\rfloor$;
    \item Fundamental, oriented: $(n-1)m+O(1)$;
    \item Positive, undirected: $(n-1)\lfloor\frac{m}{2}\rfloor$;
    \item Positive, oriented: $(n-1)m-2$.
\end{itemize}
\end{Conj}
Note that the conjectural diameters in both undirected cases are equal to each other and to the diameter of $(\mathbb{Z}/m\mathbb{Z})^{n-1}$ under natural generators. Since the latter is the abelianization of $U(n,\mathbb{Z}/m\mathbb{Z})$, one immediately gets that these values serve as a lower bound (for all $m$).

 We observe experimentally that  a similar phenomenon occurs for other unipotent groups, in particular, the maximal unipotent subgroups of symplectic groups. However, in the oriented case for fundamental roots the sequence, while very close to a linear, does not display a periodic deviation.

Computational evidence also suggests the following
\begin{Conj}
The diameter of $U(n,\mathbb{Z}/2\mathbb{Z})$ with respect to the fundamental root transvections is equal to the number of metacyclic groups of order $2^n$ (\href{https://oeis.org/A136184}{OEIS A136184}).  
\end{Conj}

Notebooks: \href{https://www.kaggle.com/code/alexandervc/growth-ut4-analysis}{UT4},
\href{https://www.kaggle.com/code/alexandervc/growth-ut5-6-7-8-etc-analysis}{UT5etc},
\href{https://www.kaggle.com/code/alexandervc/growth-in-finite-groups-heisenberg-uptriag-etc}{various groups}.

\subsection{Distribution of distances in finite nilpotent groups}
\begin{figure}[ht]
    \centering
    \includegraphics[width=0.75\linewidth]{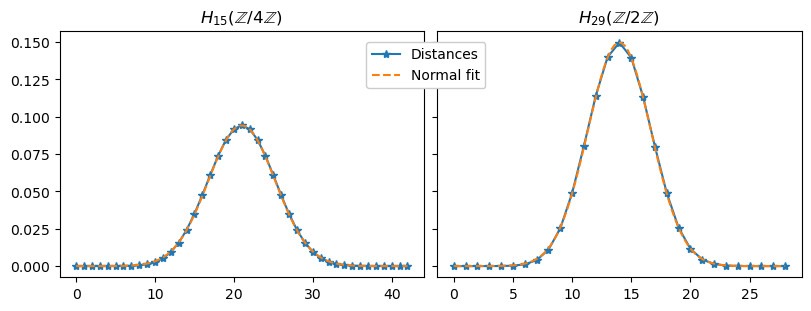}
    \caption{Distances distribution in higher-dimensional Heisenberg groups.}
    \label{fig:heisenberg_distrib}
\end{figure}

When sampling uniformly from a finite group, the distance to the identity (with respect to some fixed set of generators) can be considered as a random variable. It is easy to see that for $(\Z/m\Z)^n$ this distribution is asymptotically normal (as $n$ grows), being the sum of independent discrete uniform random variables (\href{https://www.kaggle.com/code/alexandervc/growth-commutative-groups-z-p-d}{notebook}). Experimental evidence (Fig.~\ref{fig:heisenberg_distrib}) suggests the following generalization.
\begin{Conj}
The distribution of distances in the higher-dimensional modular Heisenberg group $H_{2d+1}(\Z/m\Z)$ is approximately normal.
\end{Conj}
It is unclear whether this is also the case for $U(n,\Z/m\Z)$, which poses the following research direction: under what conditions do the (suitably normalized) growth distributions of a sequence of nilpotent groups converge to a normal distribution?

Notebooks: \href{https://www.kaggle.com/code/alexandervc/growth-commutative-groups-z-p-d}{Abelian},
\href{https://www.kaggle.com/code/alexandervc/growth-in-finite-groups-heisenberg-uptriag-etc}{various groups}.

\clearpage
\subsection{Benchmark Kaggle Challenges}

We create about a dozen benchmark datasets for Cayley graph path-finding methods and  we frame them as Kaggle challenges.
It provides a public and easy way to benchmark various methods,  adds a gamification to benchmarks,
and hopefully can attract more participants. 
Each challenge is devoted to a separate set of Cayley graph generators.
The organization of challenges is rather simple. For each challenge we provide a list 
of permutations (just vectors of integers numbers,  typically thousand of them) -- the states which are to be solved 
(file "test.csv"), 
and provide a set of generators ("allowed moves") of the permutation group.
The solution ("submission.csv" file) should consist of sequences of the names of generators 
(for each permutation  in "test.csv"),
such that application of these generators brings the input vector to the sorted form ($[0,1,2,3,4,5...]$, i.e. the identity of the group). 
In other words, the task is to decompose an input permutation into a product of given generators. 
The score of the solution is the sum of the lengths of the provided solutions sequences. 
The lower the score the better. 
The Kaggle platform provides automatic support for organization of such challenges,
and scores of all of the submitted solutions become visible on the "leaderboard" page. All of them are evaluated by the same procedure, on the same data -- which makes a fair comparison 
of different approaches possible . 
Our challenges are almost identical to the official  Kaggle challenge \href{https://www.kaggle.com/competitions/santa-2023}{Santa 2023},
with the only difference that each challenge is devoted to a separate Cayley graph.
All challenges are focused on real open mathematical problems,
for some of them we hope that the final (optimal) solution is within reach,
for the others it might be beyond the level of current technology, but still it is worth
fixing what best results are achievable currently. 

\begin{figure}[ht]
    \centering
    \includegraphics[width=0.85\linewidth]{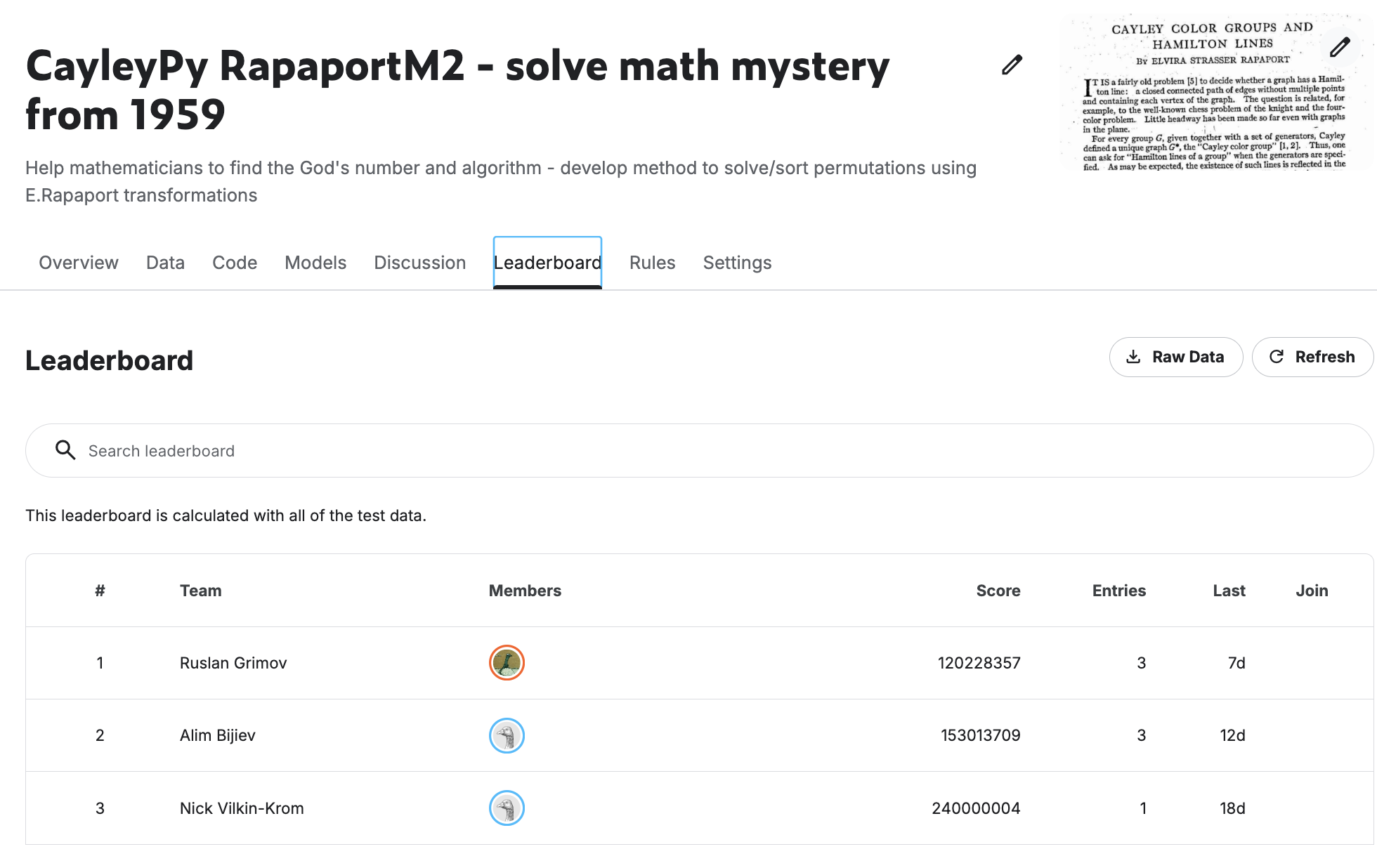}
    \caption{Kaggle challenge to benchmark Cayley graph path-finding and to determine mean diameter of the Rapaport M2 generators.
    The leaderboard score is the path length, a measure of optimality (the lower the better).}
    \label{fig:Kaggle1}
\end{figure}

Figure \ref{fig:Kaggle1} presents a leaderboard of 
\href{https://www.kaggle.com/competitions/cayleypy-rapapport-m2/leaderboard}{one the challenges}
devoted to solution of the Rappaport-M2 Cayley graphs. 
Solvers for these generators are not yet known. 

In total we created more than a dozen of Kaggle benchmark challenges with various complexity: 
\href{https://www.kaggle.com/competitions/CayleyPy-pancake}{Pancake sorting},
\href{https://www.kaggle.com/competitions/cayleypy-transposons}{Transposons},
\href{https://www.kaggle.com/competitions/cayleypy-reversals}{Reversals},
\href{https://www.kaggle.com/competitions/cayleypy-glushkov}{Glushkov problem},
\href{https://www.kaggle.com/competitions/cayleypy-rapapport-m2}{RapapportM2},
\href{https://www.kaggle.com/competitions/cayley-py-444-cube}{Rubik's cube 444},
\href{https://www.kaggle.com/competitions/cayley-py-professor-tetraminx-solve-optimally}{Professor Tetraminx},
\href{https://www.kaggle.com/competitions/cayley-py-megaminx}{Megaminx},
\href{https://www.kaggle.com/competitions/cayleypy-christophers-jewel}{Christophers jewel},
\href{https://www.kaggle.com/competitions/cayleypy-ihes-cube}{SuperCube from IHES},
\href{https://www.kaggle.com/competitions/lrx-binary-in-search-of-gods-number}{LRX-binary},
\href{https://www.kaggle.com/competitions/lrx-discover-math-gods-algorithm}{LRX},
\href{https://www.kaggle.com/competitions/lrx-oeis-a-186783-brainstorm-math-conjecture}{LRX longest}.


\subsection{LLM ResearchArena -- CayleyBench}
Path-finding for permutation groups is exactly the sorting problem in computer science.
Indeed the input is a vector of integers (permutation) and the goal is to bring it to the 
sorted form (i.e. vector $(0,1,2,3,4,...)$ -- identity of the permutation group).
The only restriction is that the allowed transformations are restricted to the generators of the group.
Classical examples include: \href{https://en.wikipedia.org/wiki/Bubble_sort}{bubble sort}
(provides an optimal algorithm for decomposition for Coxeter's generators $(i,i+1)$),
 \href{https://en.wikipedia.org/wiki/Pancake_sorting}{prefix (or pancake) sorting} (with a known suboptimal algorithm  by \cite{GatesPapadimitriou1979PrefixReversal}), etc. 

Existence of the algorithm and its complexity estimation  is often the main tool to prove upper bound  on the diameter - to show constructively that any element can be decomposed in certain number of generators (allowed operations). 

Clearly modern LLM would easily produce a code for bubble sorting, prefix (pancake) sorting
and other standard algorithms, which are existing the literature. 

But the question is -- can LLM invent algorithms for the cases not described in the literature ? 
I.e. we take some generators of $S_n$ and the task is to produce a code which would perform sorting
using only these allowed transformations. 
The output is a code and it is easy to verify it and estimate number of operations required.
If LLM would be able to invent new algorithms not known in the literature -- that would be
a solution for open mathematical problems some of them open for dozens of years - 
for example like V.M. Glushkov's problem on LX generators goes back to 1968. 

Unfortunately our current experiments show unsatisfactory results.
Out of the box results for all standard LLM which we tried
are unsatisfactory - LLM fail to produce just correct code for any cases we tried, outside well-known.
Even small modification for known generators are not solved by LLM. 

To summarize. The task to produce a sorting algorithm with restriction to given permutations generators is 
a natural framework to test LLM's abilities to solve research problems:
\begin{enumerate}
    \item Tasks are easy to formulate 
    \item Answer (code) is easy to verify
    \item There are plenty tasks of that kind and generation of more tasks is not a problem
    \item Tasks solutions would advance mathematical research and resolve problems some open for dozens years
    \item Vice versa mathematical results on estimation of diameters, mean diameters provides quality control 
    for the algorithmic outputs of LLM
\end{enumerate}

We expect that many of such sorting problems are within reach for modern AI-based technologies. 
CayleyPy supports more than half hundred different generators and provides their analysis for them,
we plan to combine it with LLM in future releases. 

\clearpage
\section{Case study: LX/LRX (Sigma/Tau) -- transposition and the longest cycle/s, V. M. Glushkov 1968 problem, Sigma-Tau (Nijenhuis, Wilf 1975) problem } \label{sect:Glushkov}
\subsection{Definition, diameter conjecture and related works}
{~}\\
{\bf Definition and background.}
One of the simplest families of generators of $S_n$ consists of just two elements -- a transposition (e.g. $(1,2)$)
and the longest cycle (cyclic shift, e.g. the left cyclic shift $L=(1,2,3,...)$). One can also add the right cyclic shift
and get a set of inverse closed generators. We will  denote these by  "LX" (respectively "LRX" when the right shift is included), following 
\href{https://oeis.org/A186783}{OEIS-A186783} which presents a conjecture that the diameter for LRX generators
is $n(n-1)/2$. This conjecture is still open, currently best upper/lower bounds and further results were proposed in \cite{CayleyPyRL}.
The "LX" case has an even  longer history. It appeared in influential work by  
"\href{https://en.wikipedia.org/wiki/Victor_Glushkov}{one of the founding fathers of Soviet cybernetics}" V.\,M.~Glushkov in~1968 \cite{glushkov1968completeness},  who obtained the first non-trivial bound $2n^2$ on the diameter. 
(He considered such problems from the viewpoint of his theory of "efficiency of micro-programs", basically aiming to find such generators
which produce all possible states in the least number of operations, so in a sense he was interested in generators with minimal diameter, with the constraint of fixing the number of generators and producing all the states). The question has been studied a lot
(survey: \cite{glukhovzubov1999lengths} pages 18-21), however it seems the precise diameter has not yet been determined. Here we present a conjectural answer: 

\begin{Conj}
The "LX" diameter is $(3n^2-8n+9)/4$ for $n$ odd, and $(3n^2-8n+12)/4$ for $n$ even (\href{https://oeis.org/A039745}{OEIS-A039745}).
\end{Conj}

The conjecture is a demonstration  of our "diameter is quasi-polynomial" conjectural principle.
By efficient algorithms of CayleyPy we computed \href{https://github.com/cayleypy/cayleypy/blob/main/cayleypy/data/lx_cayley_growth.csv}{full growth} and diameters up to $n\le 15$, then fitted them as polynomials.
A beautiful fact is that leading terms in growth (not diameters) precisely coincide with Fibonacci numbers (the reason is simple: ...).
The leading term for diameter $3n^2/4$ has been known before (survey: \cite{glukhovzubov1999lengths}, add: Pogorelov, Zubov)
and is consistent with \href{https://oeis.org/A048200}{OEIS-A048200} (\cite{KuppiliChitturi2020LE}) which conjectures length of NOT the longest
element, but expected close to longest i.e. just the full flip: first exchanges with last, second with next-to-last, etc.... 
The examples of the longest elements for small $n$: \href{https://oeis.org/A378881}{OEIS-A378881},
the conjectural formula for them:  \href{https://oeis.org/A186144}{OEIS-A186144}, (one for even, two for odd, $n>4$).

{\bf LX (Sigma-Tau 1975) hamiltonian path problem (resolved)}. 
Not exactly in the present scope, but worth mentioning: there was a long-standing problem resolved positively recently \cite{SawadaWilliams2019SigmaTau} of whether the directed graph generated by L and X (also denoted Sigma and Tau)
has a hamiltonian path.  D.Knuth assigned to that problem complexity 48 out of 50, it was first articulated in
\cite{NijenhuisWilf1975CombinatorialAlgorithms}.

{\bf Circular permutations (LRX and \cite{alon2025circular})}.
Let us remark that it is important in bioinformatics  to consider circular genomes (bacterial or mitochondrial), that is, to consider permutations up to cyclic circular shifts. It is also interesting in mathematics, and 
from the general perspective it can be rephrased as follows:
add generators $L,R$ to any original set of generators, construct a new Cayley graph and set the weights for $L,R$ equal to zero.
Distances on that new graph exactly represent the circular variant of the original graph. 
For example, in the case of the LRX graph, one can observe that we have $L^i X R^i = (i,i+1) $, and thus the LRX graph with L,R weights
set to zero corresponds to circular permutations generated by $(i,i+1)$. This was the subject of recent study \cite{alon2025circular},
where the diameter has been established. By the argument above it gives the lower bound for the diameter of the LRX graph,
which is the same as obtained in \cite{CayleyPyRL} (and vice versa -- \cite{CayleyPyRL} gives a lower bound for the case studied in \cite{alon2025circular}).

\subsection{Further results and conjectures on growth}

Diameters are not the only interesting features of Cayley graphs, more generally understanding of growth itself is also important.
Here we present the analysis for the LX case (see \cite{CayleyPyRL} for LRX).
It is well-known from works of P.Diaconis that growth for Coxeter generators is well-approximated by the Gaussian, and the
LRX generators are close to Coxeter, e.g. it is easy to express transpositions as $L^i X R^i = (i,i+1) $ (via circular permutations like above),
so the natural question is whether the growth tends to be Gaussian for $n\to \infty$.
Our analysis suggests a negative answer, and moreover we see that the Gumbel distribution is not a bad
fit.

The key feature which we observe is a significant asymmetry (skewness) of the distribution, at least for those $n$ which are available.
And we believe this is true in general, thus a symmetric Gaussian distribution is not an option here. 
It is somewhat natural to expect asymmetry (skewness) for the following reason. We can think of two extremes -- the first is a highly asymmetric distribution with exponential growth and abrupt fall-off
typical for random generators (see e.g. the \href{https://mathoverflow.net/q/322877/10446}{Rubik's cube example}),
and the other extreme is the Gaussian distribution which appears for generators that are close to commutative in some sense
(like Coxeter's ones).
The LRX case is in a sense not so far from Coxeter, but not really that close to it or commutative,
so it shifts to the other extreme -- and therefore asymmetry appears.

Let us formulate our observations in the form of a conjecture:
\begin{Conj} For LX growth distribution for all $n\gg 1$ we have:  asymmetry is less than -0.5, kurtoisis is greater than $0.5$.
Mean and mode can be approximated quadratically as: 
$0.58n^{2} - 2n + 2.4$,  odd modes:
$0.6 n^{2} - 2n + 2.7$,  even modes: $0.5n^{2} - 0.5 - 2.9$. 
\end{Conj}
It might be also that the distribution tends to the Gumbel one when $n\to \infty$, however that is not completely clear
since empiric skew and kurt are not well fitted (however it might be an effect of small $n$).

\begin{figure}[ht]
  \centering
  \includegraphics[width=450pt]{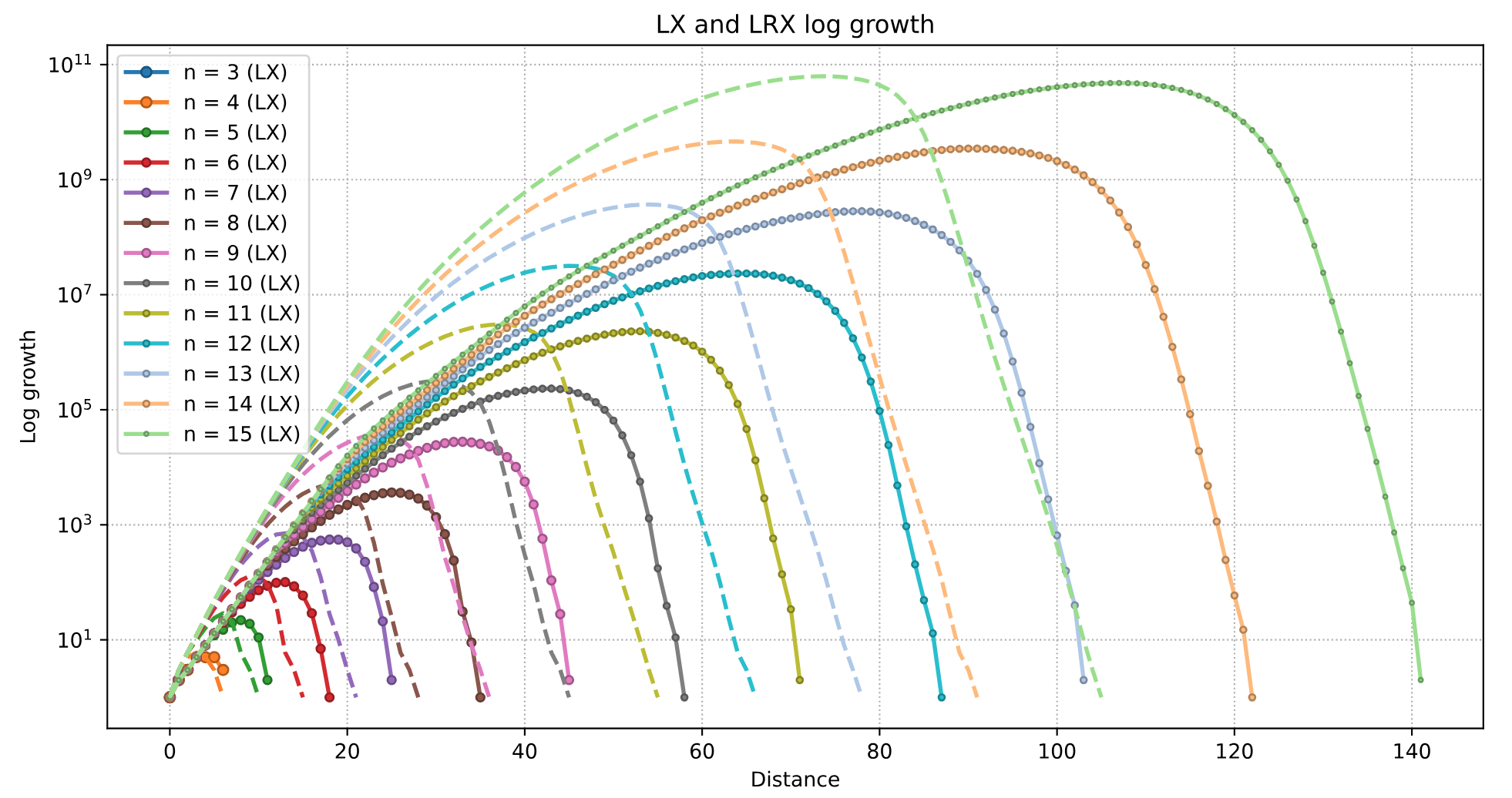}
  \caption{LX and LRX growth log10 scale}\label{fig:lx_growth}
\end{figure}

\begin{figure}[ht]
  \centering
  \includegraphics[width=450pt]{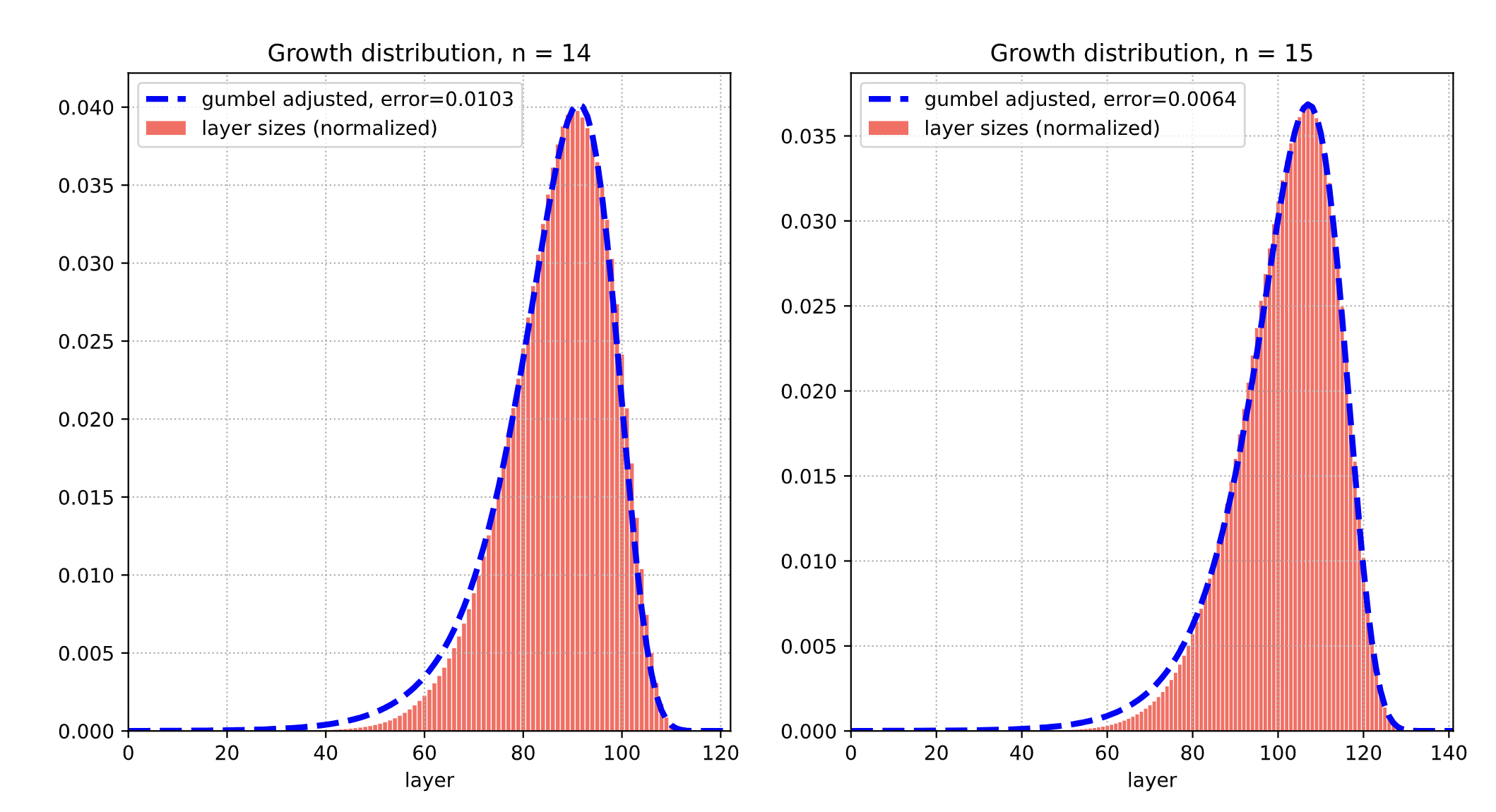}
  \caption{LX growth growth Gumbel distribution fit}\label{fig:lx_gumbel}
\end{figure}

\begin{figure}[ht]
  \centering
  \includegraphics[width=450pt]{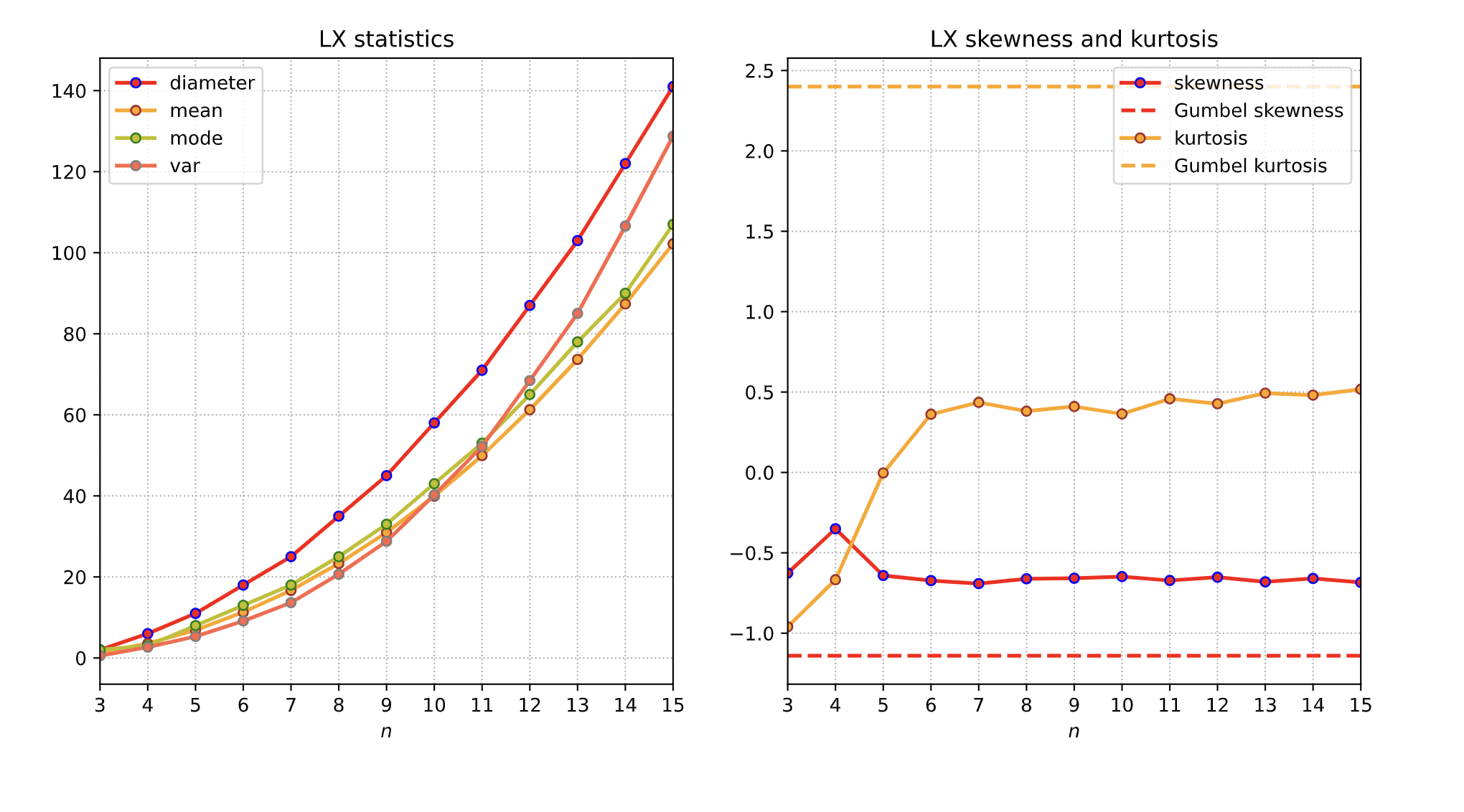}
  \caption{LX growth distribution  statistics }\label{fig:lx_stat}
\end{figure}

Notebooks: 
\href{https://www.kaggle.com/code/fedmug/lx-growth-distribution}{LX}, 
\href{https://www.kaggle.com/code/fedmug/lrx-growth-distribution}{LRX},
\href{https://www.kaggle.com/code/fedmug/lrx-coset-growth}{LRX-coset}.

\clearpage
\section{Case study -- Biologically Relevant Generators} \label{sec:bio-generators} 
\subsection{Definitions and review}
Cayley graph methods have direct applications in computational biology. In comparative genomics, a key task is estimating the evolutionary distance between two genomes. This distance is typically defined as the minimum number of large-scale rearrangement events required to transform one genome into the other.

When genomes are represented as permutations and rearrangements are represented as group generators, this task becomes a distance problem in a Cayley graph.

\subsubsection{Group‐theoretic Preliminaries} 

To model genome rearrangements formally, we represent a genome as a permutation on $\{1,2,\dots,n\}$. When gene orientation is relevant, signed permutations are used instead. These form the symmetric group $S_n$ and the hyperoctahedral group $B_n$, respectively.

Each type of rearrangement defines a specific transformation of the genome and can be treated as a generator of the corresponding group \cite{bhatia2018position}. Let $S$ be a set of allowed rearrangement operations (generators). Then the rearrangement distance $d_S(\pi,\sigma)$ between two genomes $\pi$ and $\sigma$ is the smallest number $k$ of rearrangement operations $s_i$ needed to transform $\sigma$ to $\pi$. Formally:

\[
  d_S(\pi,\sigma)
  \;=\;
  \min\bigl\{\,k : \pi = s_1\,s_2\cdots s_k\,\sigma,\;s_i\in S\bigr\},
\]

This is equivalent to computing the word length of $\pi \sigma^{-1}$ in the Cayley graph of $G$ with respect to $S$. In the biological applications, $S$ typically includes large-scale rearrangement events such as
\[
  S \;=\;\{\text{inversion},\,\text{translocation},\,\text{fusion},\,\text{fission},\dots\}.
\]

To effectively analyze genome rearrangements using Cayley graphs, it's essential to understand how biological events correspond to mathematical operations.

\subsubsection{Core Generators}

The most studied genomic rearrangements correspond to the following types of generators:

\begin{itemize}
  \item \textbf{Reversals (Inversions)}: flipping a contiguous segment of a chromosome, 
    $$\dots\,a\,\underbrace{b\;c\;\dots\;d}_{\text{reversed}}\;e\,\dots \;\longmapsto\;\dots\,a\,d\;\dots\;c\;b\;e\,\dots.$$
    In the signed case, each reversal also toggles the orientation of every gene in the segment \cite{Pevzner1999cabbage2turnip}.
    Biological motivation comes from double sided DNA breaks:  when two breaks of DNA happens
    a fragment of DNA leaves its position, but cell's reparation system brings it back,
    however sometimes mistakes happens and reparation system brings it back in the reversed order,
    thus reversals appears. (Cell's repair machinery doesn’t distinguish orientation — it just ligates available DNA ends. If the piece has been flipped during the breakage–rejoining process, 
    the outcome is a reverse (inversion)).
  \item \textbf{Transposons} ("transpositions" in biological literature - different from mathematical "transpositions"): cutting out a segment and reinserting it elsewhere on the same chromosome (that segment is abstraction of biological \href{https://en.wikipedia.org/wiki/Transposable_element}{transposon} - transposable element of genome):
    $$\dots\,u\;\underbrace{v\;\dots\;w}_{\Delta}\;x\;\dots\;y\;\dots 
      \;\longmapsto\;\dots\,u\;x\;\dots\;y\;\underbrace{v\;\dots\;w}_{\Delta}\;\dots.$$  
  \item \textbf{Translocations}: moving a segment from one chromosome to another,
    $$\underbrace{\dots\,A\,\Delta\,B\,\dots}_{\text{chr.\,1}}
      \quad\longleftrightarrow\quad
      \underbrace{\dots\,C\,\Delta'\,D\,\dots}_{\text{chr.\,2}}.$$
  \item \textbf{Fissions}: splitting a single chromosome into two at one breakpoint,
    $$\dots\,p\;\underbrace{q\;\dots\;r}_{\Delta}\;s\,\dots 
      \;\longmapsto\;(\dots\,p\,)\;(\underbrace{q\;\dots\;r}_{\Delta}\;s\,\dots).$$
  \item \textbf{Fusions}: the inverse of a fission, joining two chromosomes end‐to‐end.
\end{itemize}

Different choices of $S$ and $G$ yield distinct distances and computational complexities \cite{caprara1999sorting, bulteau2012transp, oliveira2019complexity}. In some cases, particularly when the allowed operations do not form a group or when the rearrangement space lacks closure under composition (e.g., unsigned reversals or multi-chromosomal translocations), the problem cannot be modeled via Cayley graphs. Instead, alternative frameworks such as breakpoint graphs, adjacency graphs, or double-cut-and-join graphs are employed \cite{BafnaPevzner1998reversals, yancopoulos2005dcj}.

\subsubsection{Key Rearrangement Models and Algorithms}

In this work, we specifically focus on genome rearrangement problems that can be formulated and solved on Cayley graphs — that is, problems where the allowed operations form a closed generating set of a group. While other operations such as transpositions can also generate \(S_n\) and thus allow a Cayley graph formulation in theory, such representations tend to overlook biologically realistic constraints (e.g., locality or chromosome structure). Depending on the chosen rearrangement model and whether the permutations are signed or unsigned, the computational complexity of finding minimal rearrangement distances varies significantly, as summarized in Table \ref{tab:complexity}.

\begin{table*}[htbp]
\caption{Overview of complexity and approximation in genome rearrangement problems with Cayley-compatible models\label{tab:complexity}}
\begin{tabular*}{\textwidth}
  {@{\extracolsep{\fill}}
   l
   p{2.8cm}  
   p{1.9cm}  
   p{2.8cm}  
   p{2.2cm}  
  @{}}
\toprule
\textbf{Problem} &
\textbf{S Complexity} &
\textbf{S Approx.} &
\textbf{UNS Complexity} &
\textbf{UNS Approx.} \\
\midrule
REV &
P (\(O(n\sqrt{n\log n})\))~\cite{tannier2004sorting} &      
1.375~\cite{berman2002} &                                   
NP-hard~\cite{caprara1999sorting} &                         
1.5 ~\cite{christie1998}                                    
\\
TR &
NP-hard~\cite{bulteau2012transp} &                
– &                                               
NP-hard~\cite{elias2006transp} &                  
1.375~\cite{elias2006transp}                      
\\
REV+TR &
NP-hard~\cite{oliveira2019complexity} &                  
2~\cite{walter1998reversal} &                            
NP-hard~\cite{oliveira2019complexity} &                  
2.83+\(\epsilon\)~\cite{rahman2008revtrans}              
\\
SSOs & P\textsuperscript{*} ~\cite{oliveira2018} & – &
P\textsuperscript{*} ~\cite{oliveira2018} & –
\\
\bottomrule
\end{tabular*}

\vspace{1em}
Approx.: approximation factor (indicates the worst-case ratio between the computed solution and the true minimal rearrangement distance);
P: polynomial time; REV: reversals; TR: transpositions; S: signed; SSO: super short operations; UNS: unsigned.\\
\textsuperscript{*} All known variants of sorting by super short operations are solvable in polynomial time, with time complexity ranging from \( O(n) \) to \( O(n^4) \), depending on the model (see ~\cite{oliveira2018} for references).
\end{table*}

\subsubsection{Reversals }

{\it The Reversal graph} $R_n, n \geqslant 3,$ is the Cayley graph over the symmetric group generated by the reversals from the set
$R=\{r_{i,j}\in Sym_n, 1 \leqslant i < j \leqslant n\}, \ |R|=\binom{n}{2}$,
where any $r_{i,j}$ reverses elements of a permutation $\pi=[\pi_1 \pi_2 \ldots \pi_n]$ within segment $[i,j]$ in a such a way that we have:

$$[\pi_1 \ldots \underline{\pi_i \pi_{i+1} \ldots \pi_{j-1} \pi_j} \ldots \pi_n]r_{i,j}=[\pi_1 \ldots \underline{\pi_j \pi_{j-} \ldots \pi_{i+1} \pi_i} \ldots \pi_n].$$

The distance in this case is defined as the minimal number of reversals transforming one permutation into another. It is known~\cite{BP96,KeSa95} that the diameter of this graph is $(n-1)$, and the only permutations needing this many reversals are the Gollan permutation $\gamma_n$ and its inverse, where the Gollan permutation, in a cycle notation, is defined as follows:

\[ \gamma_n= \left\{
\begin{array} {lr}
(3\,1\,5\,2\,7\,4\,\ldots\,n-3\,n-5\,n-1\,n-4\,n\,n-2), & \mbox{if } n \mbox{ is even} \\
(3\,1\,5\,2\,7\,4\,\ldots\,n-6\,n-2\,n-5\,n\,n-3\,n-1), & \mbox{if } n \mbox{ is
odd.} \end{array} \right. \] 

The analysis of genomes evolving by inversions leads to the combinatorial problem of sorting by reversals. Reversal distance measures the amount of evolution that must have taken place at the chromosome level, assuming evolution proceeded by inversions. There are two algorithmic subproblems. The first one is to find the reversal distance between two permutations. The second one is to reconstruct a sequence of reversals realizing the distance. The last problem is NP--hard~\cite{KeSa94}. There are 1.5--approximation~\cite{Ch98}  
and 1.375--approximation~\cite{BHK02} algorithms for sorting by reversals.

Thus, the reversal graph $R_n, n \geqslant 3,$ is a connected $\binom{n}{2}$-regular 
graph of order $n!$ and diameter $n-1$. It was shown in~\cite{K07} that it does not contain neither triangles nor subgraphs isomorphic to $K_{2,4}$. This gives explicit formulas for the cardinality of the metric spheres $S_1$ and $S_2$ in this graph as 
$|S_1|=\binom{n}{2}$-regular 
and $|S_2|=\frac{1}{6}(n^4-2n^3-n^2-16n+42).$ For any $n \leqslant 14$, the cardinalities $|S_i|, 0\leqslant i \leqslant n-1,$ i.e. growth
can be found at \href{https://oeis.org/A300003}{OEIS-A300003}.

\begin{figure}[H]
  \centering
  \begin{minipage}{0.55\textwidth}
    \centering
    \includegraphics[width=1.3\linewidth]{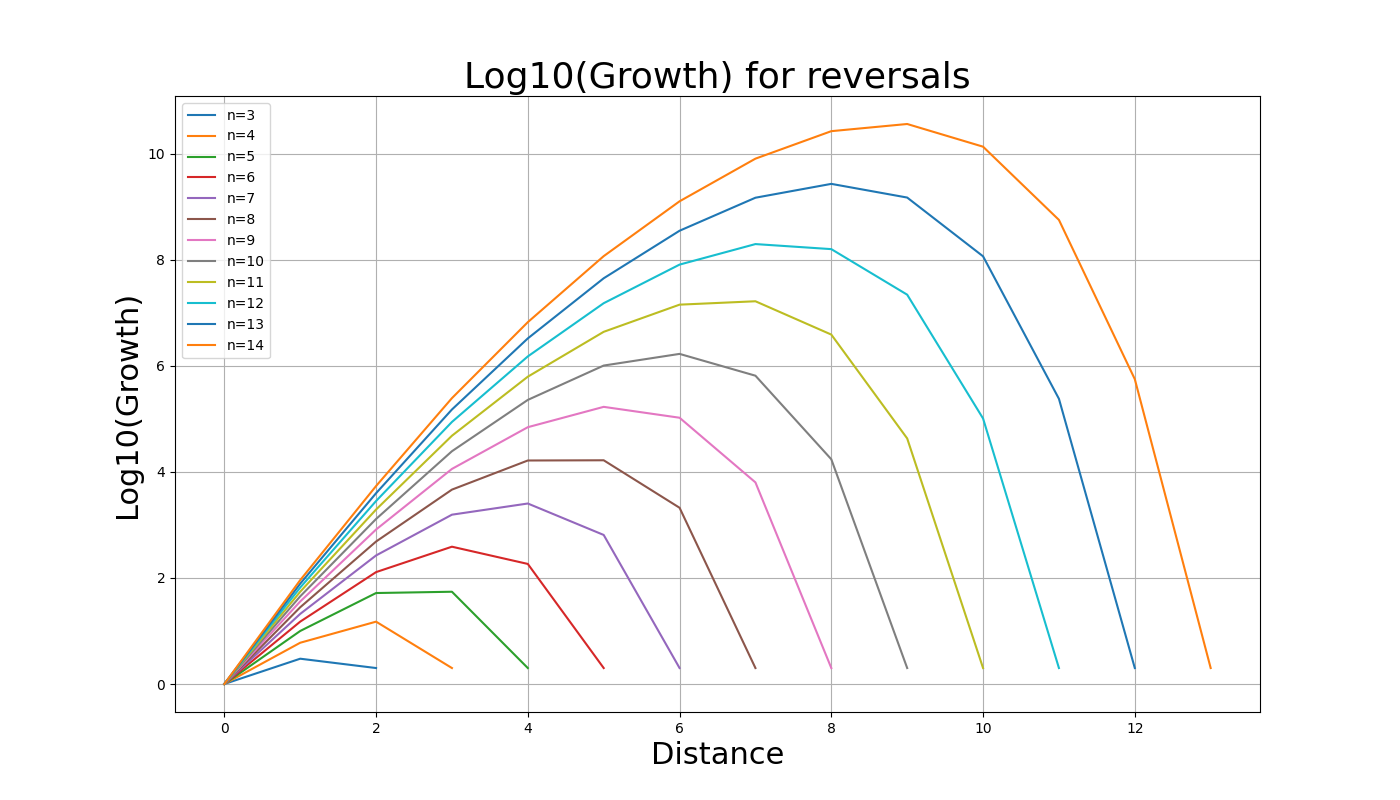}
    \caption{Growth for reversals}
    \label{fig:reversals-growth}
  \end{minipage}\hfill
  \begin{minipage}{0.35\textwidth}
    \centering
    \includegraphics[width=0.7\linewidth]{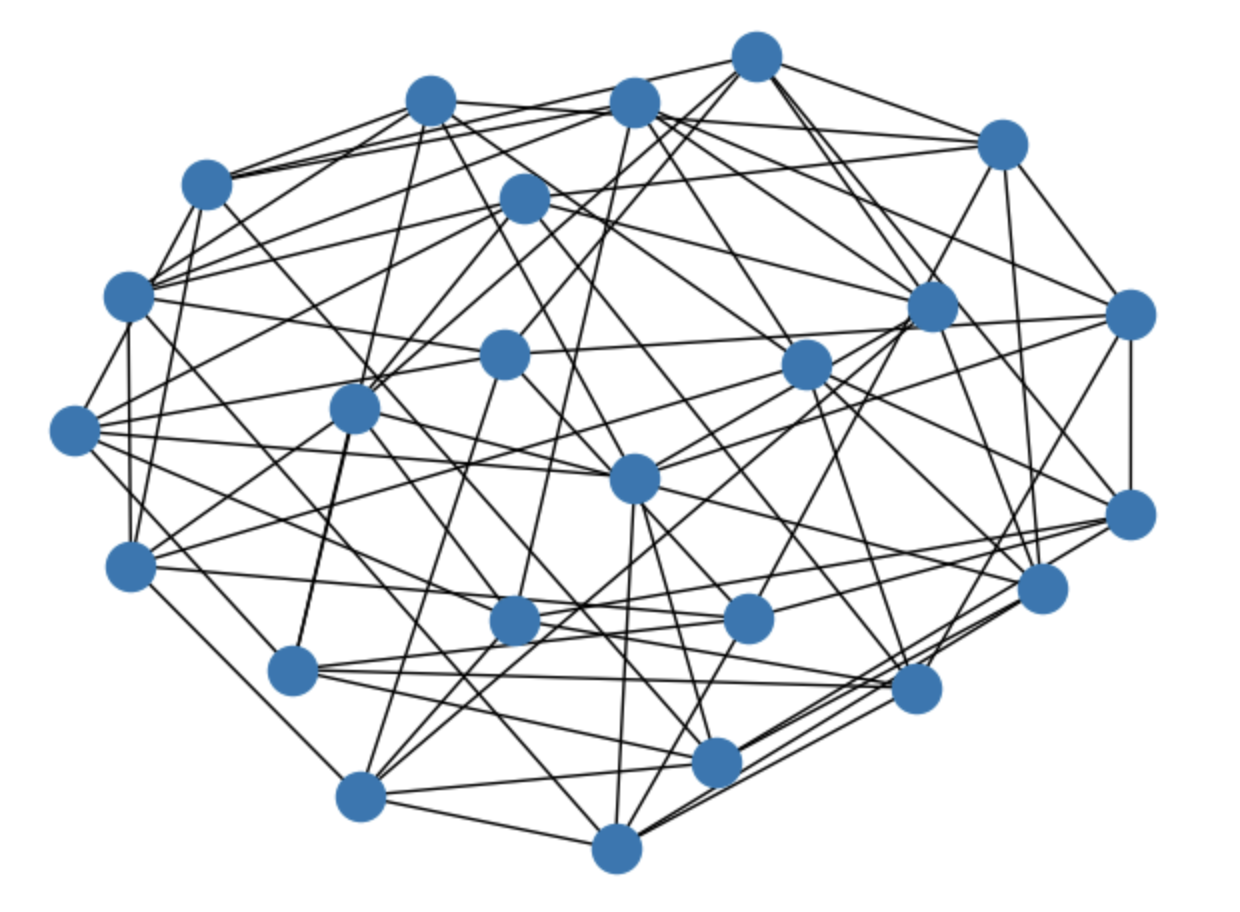}
    \caption{The reversal graph, $n = 4$}
    \label{fig:reversals4}
  \end{minipage}
\end{figure}

\vspace{5mm}



\subsubsection{Transposons (block transpositions) }
{~}\\
The generators called somewhat misleadingly in biological literature "transpositions" not to be confused with standard 
mathematical transpositions, better called "block transpositions", or we will call them just "transposons" - since their biological motivation
is to describe "\href{https://en.wikipedia.org/wiki/Transposable_element}{transposons}" (movig fragments of genome). 
They can be described as follows:  cutting out a segment and reinserting it elsewhere on the same chromosome
(that segment is abstraction of biological \href{https://en.wikipedia.org/wiki/Transposable_element}{transposon}):
    $$\dots\,u\;\underbrace{v\;\dots\;w}_{\Delta}\;x\;\dots\;y\;\dots 
      \;\longmapsto\;\dots\,u\;x\;\dots\;y\;\underbrace{v\;\dots\;w}_{\Delta}\;\dots.$$  

They are supported by CayleyPy under the name "\href{https://cayleypy.github.io/cayleypy-docs/generated/cayleypy.PermutationGroups.html#cayleypy.PermutationGroups.transposons}{transposons}".

{\bf Conjecture.}(\href{https://oeis.org/A065603}{OEIS-A065603})  
Diameter of transposons is $\left\lceil \tfrac{n+1}{2}\right\rceil$ for $n\ne 13,15$. 

The conjecture formulated in \cite{ErikssonErikssonKarlanderSvenssonWastlund2001BridgeHand} (section 3 end), based on earlier
results in \cite{BafnaPevzner1998Transpositions}.

\begin{figure}[H]
  \centering
  \begin{minipage}{0.55\textwidth}
    \centering
    \includegraphics[width=1.3\linewidth]{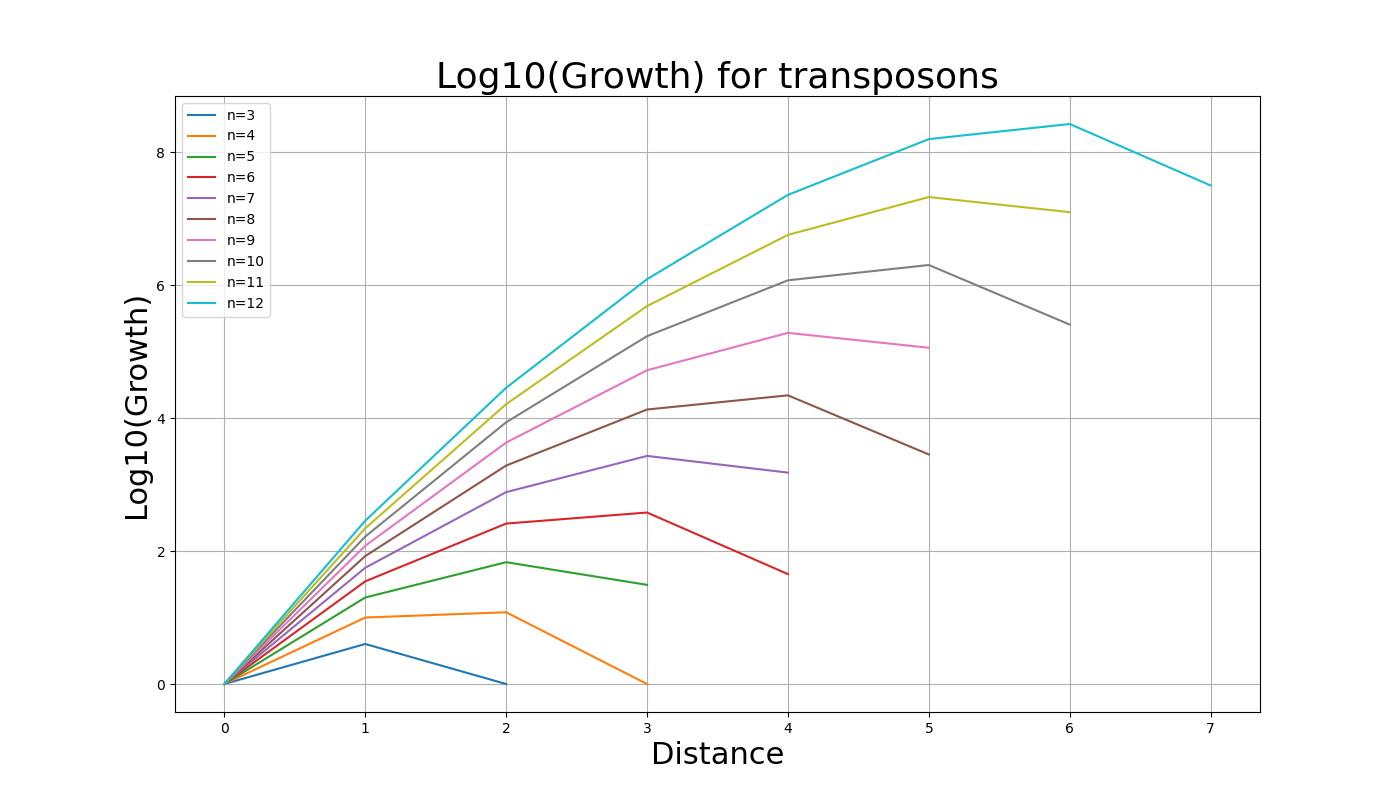}
    \caption{Growth for transposons}
    \label{fig:transposons-growth}
  \end{minipage}\hfill
  \begin{minipage}{0.35\textwidth}
    \centering
    \includegraphics[width=0.7\linewidth]{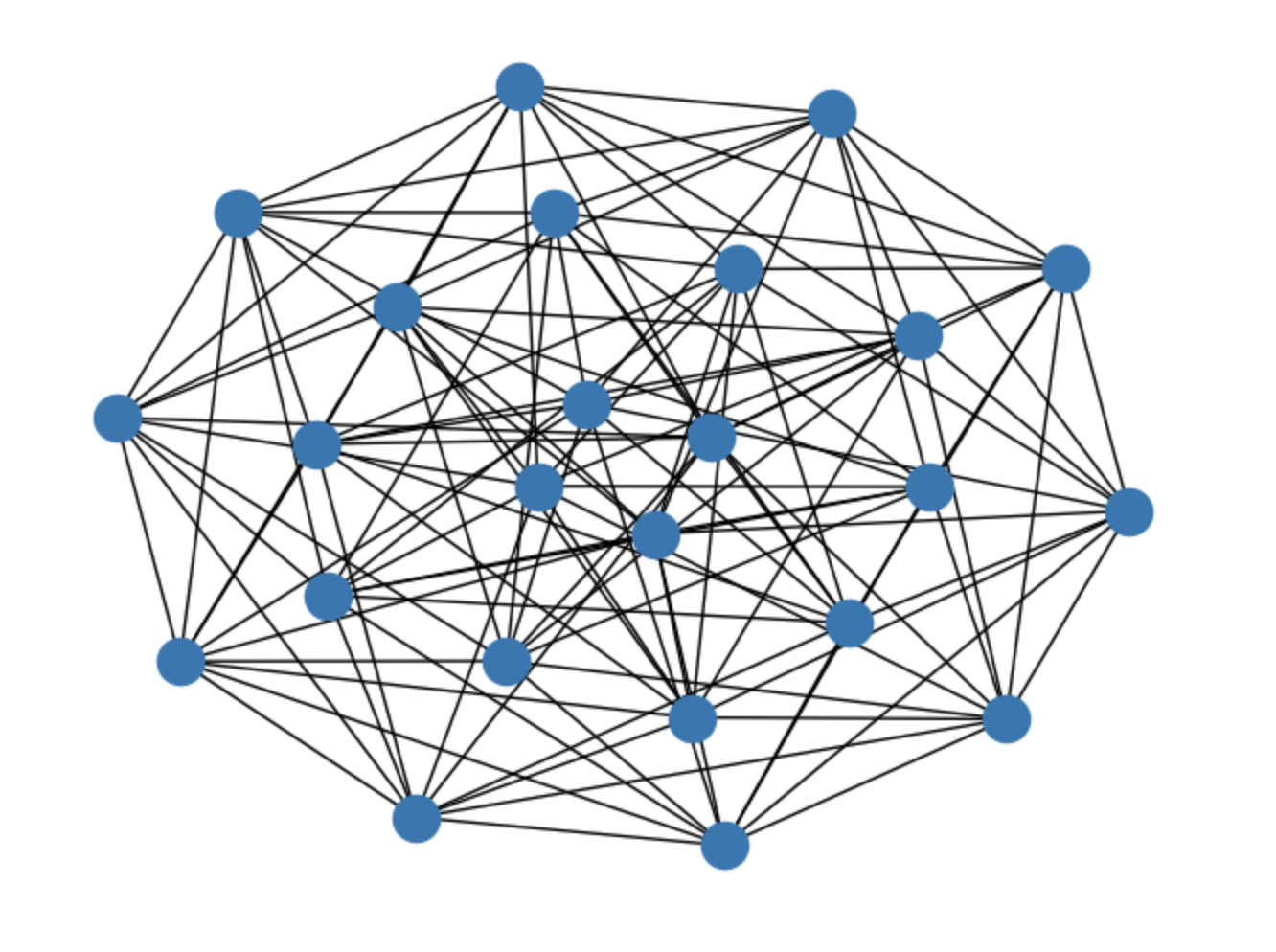}
    \caption{The Transposone graph, $n = 4$}
    \label{fig:transposons4}
  \end{minipage}
\end{figure}


\subsection{Growth of Reversals and Transposons on binary cosets}

Here we describe new findings on  Schreier coset graphs for both reversal and transposons generators.
The graph can be described in rather simple way - consider the action of permutation generators
on sequences of length $n$  such that these sequence contain only 2 distinct elements which can be denoted by 0 and 1.
Moreover assume that there precisely $\lfloor n/2 \rfloor$ zeros, and $n-\lfloor n/2 \rfloor$ ones. 
Biologically that would mean that organisms have only two genes, which have $\lfloor n/2 \rfloor$ and $n-\lfloor n/2 \rfloor$ duplicates each.
Not so meaningful biologically, but of interest from mathematical point of view:
it is clear simplification of the original Cayley graphs and their diameters are not larger than original diameters,
thus provide a way to get lower bounds on full Cayley graphs.
Growth for Schreier graphs may depend on the choice of initial vector, here we take the most simple variant:
the sorted vector i.e. first 0, then 1. 

Surprisingly we can guess almost all desired characteristics of these graphs:

\begin{Conj} For binary cosets above:  
    \begin{enumerate}
        \item The growth of reversals and transposons - coincide
        \item God's number is equal to $\left\lfloor \tfrac{n}{2}\right\rfloor$  
        \item The growth can be described by explicit formulas:\\
        \textbf{Even \(n=2m\):} \(a(n,k)=\binom{m}{k}^2, \; 0\le k\le m\).\\
        \textbf{Odd \(n=2m+1\):} \(a(n,k)=\binom{m}{k}\binom{m+1}{k}, \; 0\le k\le m\).
        \item Mode is equal to $\left\lfloor \tfrac{n+1}{4}\right\rfloor$ 
        \item Mean approximately equals to $\approx n/4$
        \item Antipodes: There is single longest state for even $n$, and $ \tfrac{n+1}{2} $ for odd  \\
        For even $n$ the longest element is [1,0,1,0,1,0,....] - multiple repeat of 1,0\\
        For odd $n$ longest state have similar structure but single 1 is doubled - and since there are 
        $\lfloor n/2 \rfloor+1$ such ones
        we get as much longest states (see examples: \href{https://www.kaggle.com/code/eveelina/fork-of-cayleypy-transposons-analys-diameter-on-n}{notebook})       
        \item Skewness is equal to zero for even $n$, and tends to zero for $n$ odd
        \item Kurtosis tends to zero 
        \item Distribution tends to Gaussian 
    \end{enumerate}
\end{Conj}

Figures below illustrate numerical evidences for the conjecture above.
It is surprising that diameter of full and coset transposons differs so little. So good lower bound for the coset diameter
gives good lower bound for the full diameter, though in that case sharp estimate is known, but might be interesting in general.

To be more precise we are dealing not with the diameter of the graph - but with the maximal distance to the initial state
("God's number"). For Cayley graphs there is no dependence on initial state, 
but it exists for Schreier coset graphs. We always take sorted vector as initial state - first all zeros, than all ones. 

\begin{figure}[H]
  \centering
  \begin{minipage}{0.45\textwidth}
    \centering
    \includegraphics[width=1.1\linewidth]{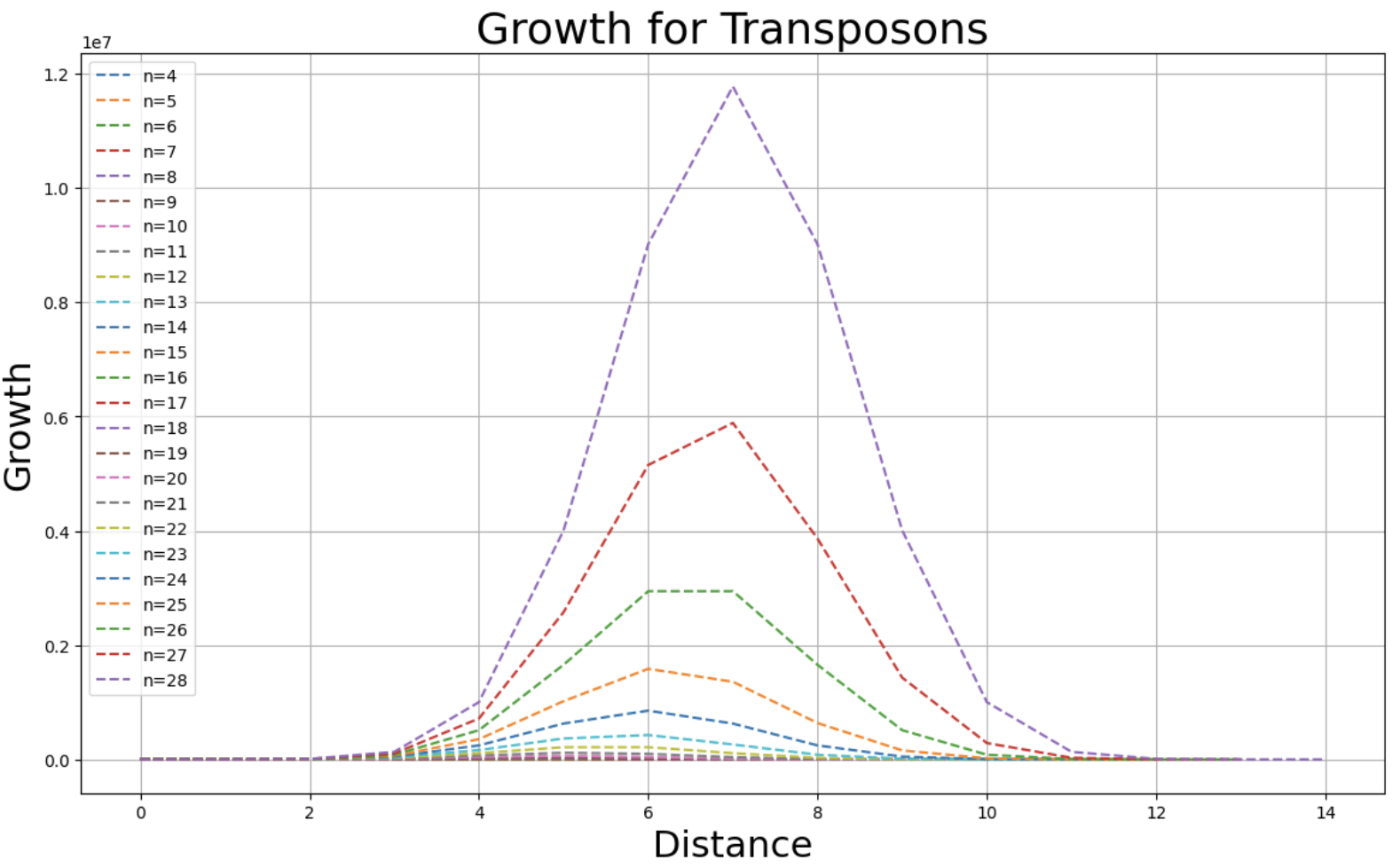}
    \caption{Growth transposone/reversals binary coset}
    \label{fig:transposons-coset-growth}
  \end{minipage}
  \hfill
  \hfill
  \begin{minipage}{0.45\textwidth}
    \centering
    \includegraphics[width=1.1\linewidth]{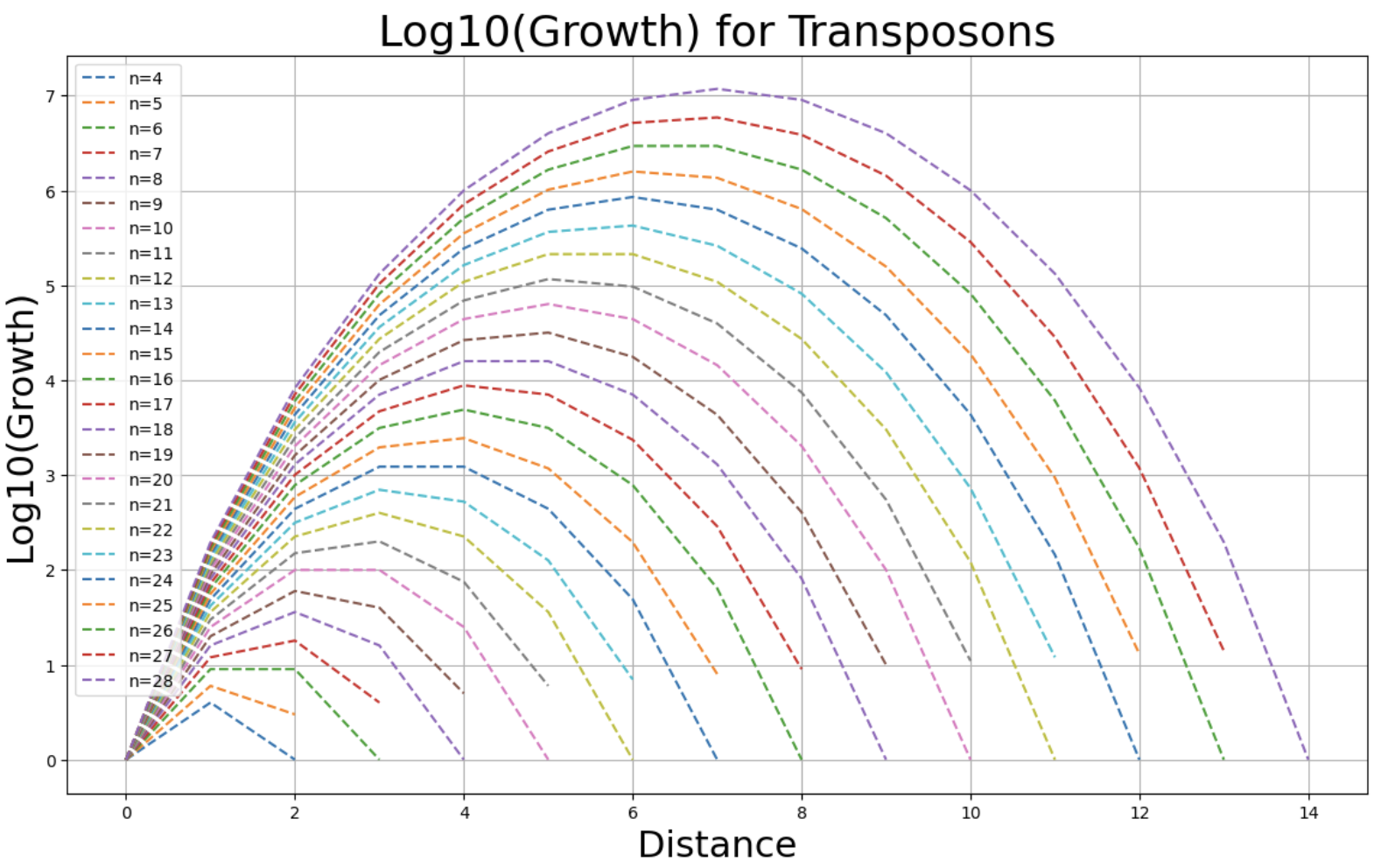}
    \caption{Log10 growth  transposone/reversals binary coset}
    \label{fig:transposons-coset-growth-log}
  \end{minipage}
\end{figure}

\begin{figure}[H]
  \centering
  \begin{minipage}{0.45\textwidth}
    \centering
    \includegraphics[width=1.1\linewidth]{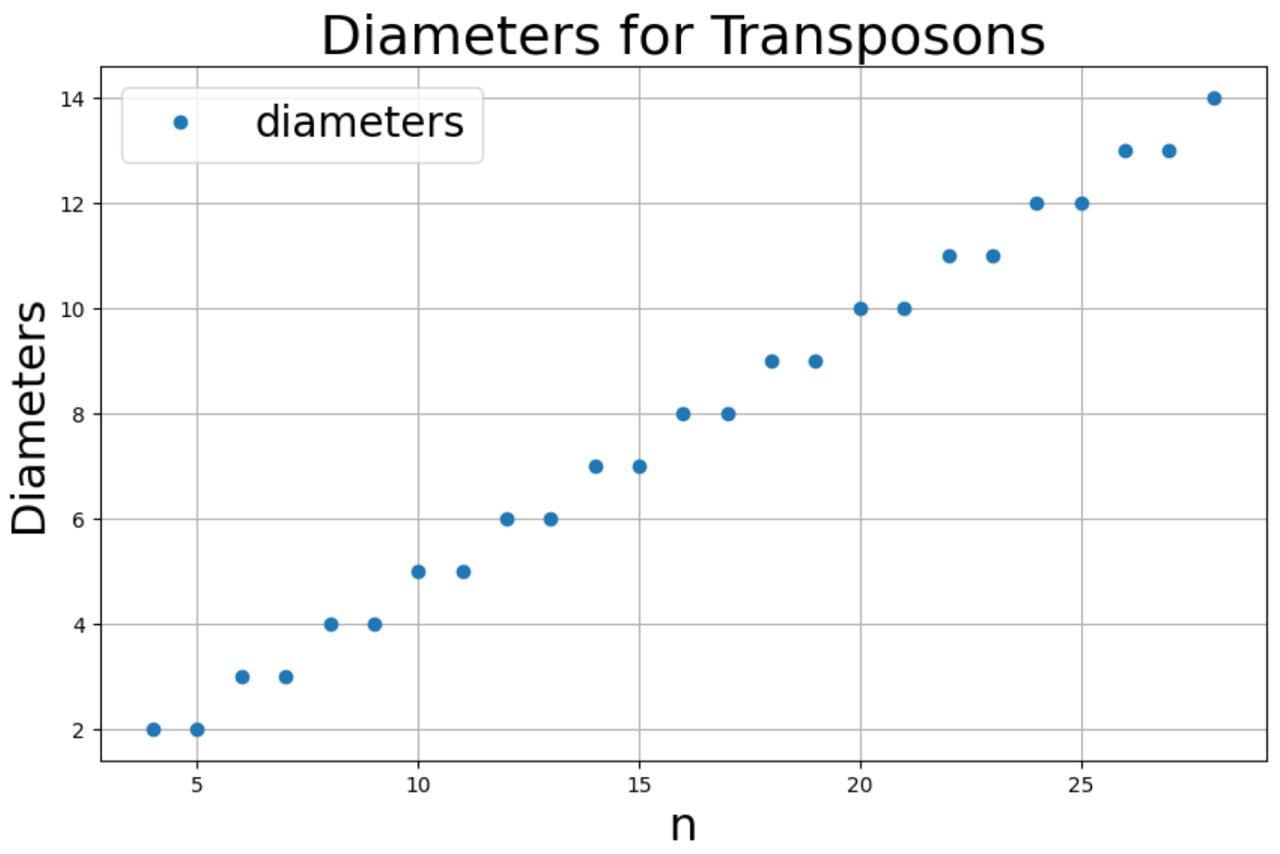}
    \caption{Diameters for transposone/reversals binary coset}
    \label{fig:transposons-coset-diam}
  \end{minipage}
  \hfill
  \hfill
  \begin{minipage}{0.45\textwidth}
    \centering
    \includegraphics[width=1.1\linewidth]{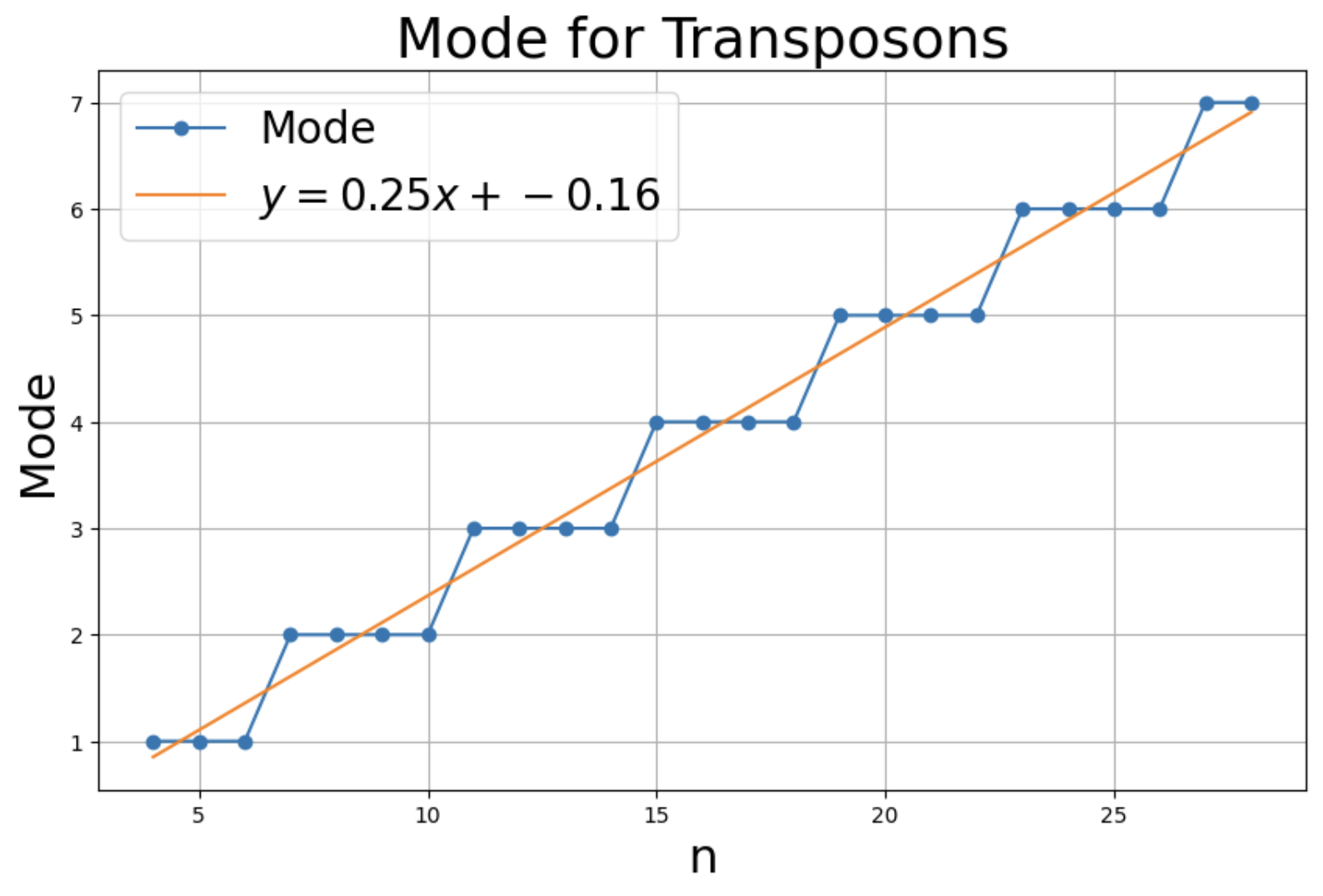}
    \caption{Mode  transposone/reversals binary coset}
    \label{fig:transposons-coset-mode}
  \end{minipage}
\end{figure}

\begin{figure}[H]
  \centering
  \begin{minipage}{0.45\textwidth}
    \centering
    \includegraphics[width=1.1\linewidth]{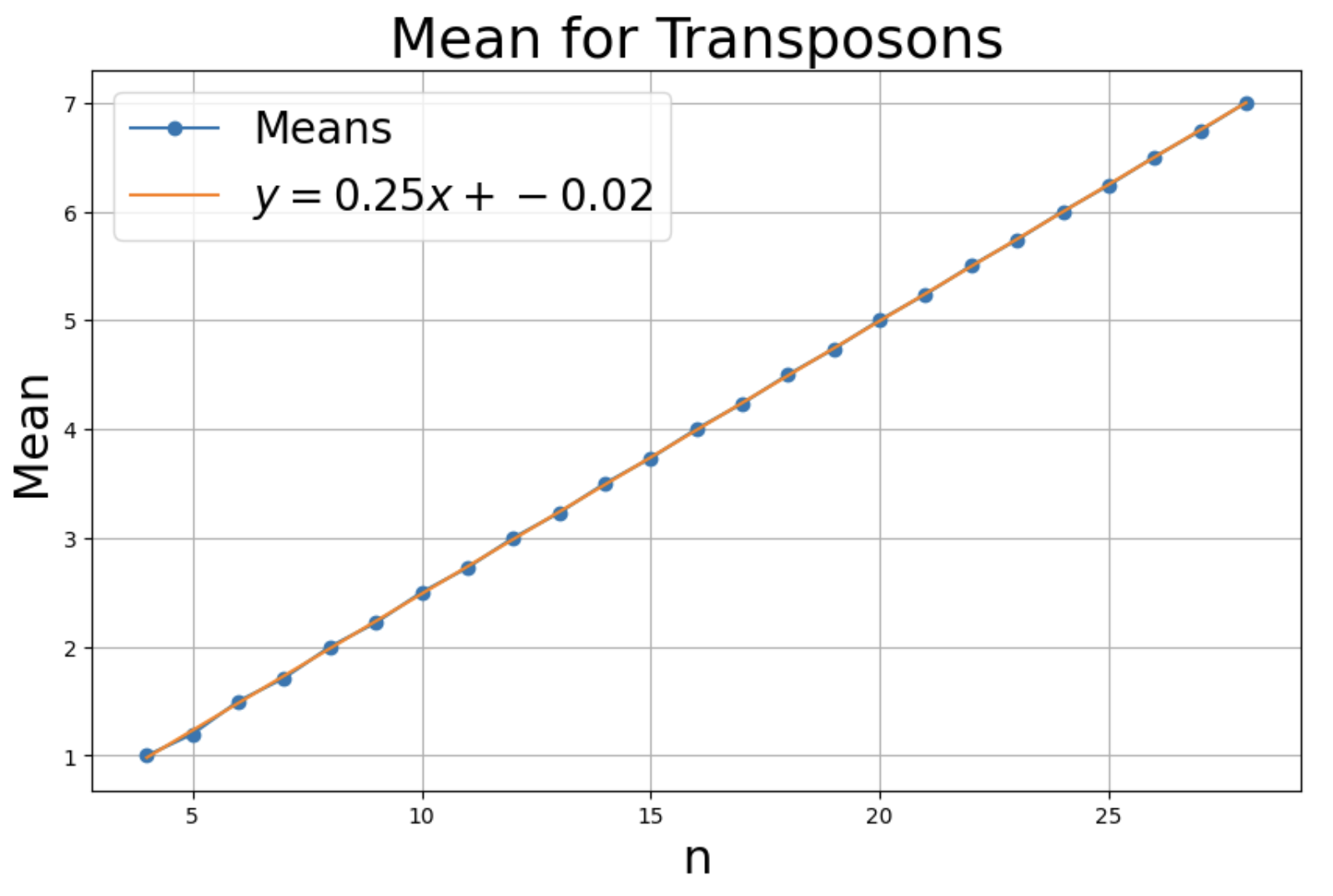}
    \caption{Mean diameter for transposone/reversals binary coset}
    \label{fig:transposons-coset-mean}
  \end{minipage}
  \hfill
  \hfill
  \begin{minipage}{0.45\textwidth}
    \centering
    \includegraphics[width=1.1\linewidth]{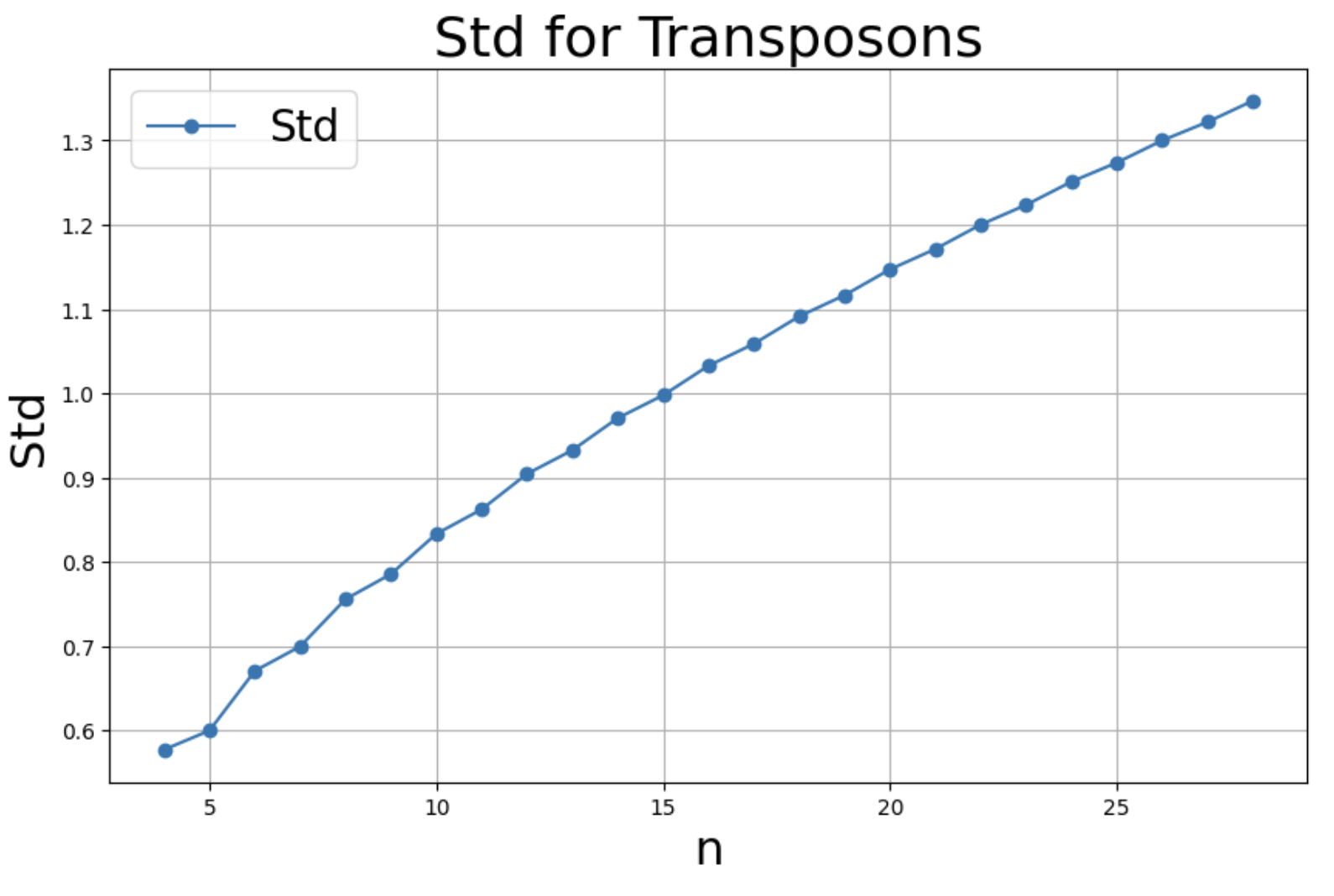}
    \caption{Std transposone/reversals binary coset}
    \label{fig:transposons-coset-std}
  \end{minipage}
\end{figure}

\begin{figure}[H]
  \centering
  \begin{minipage}{0.45\textwidth}
    \centering
    \includegraphics[width=1.1\linewidth]{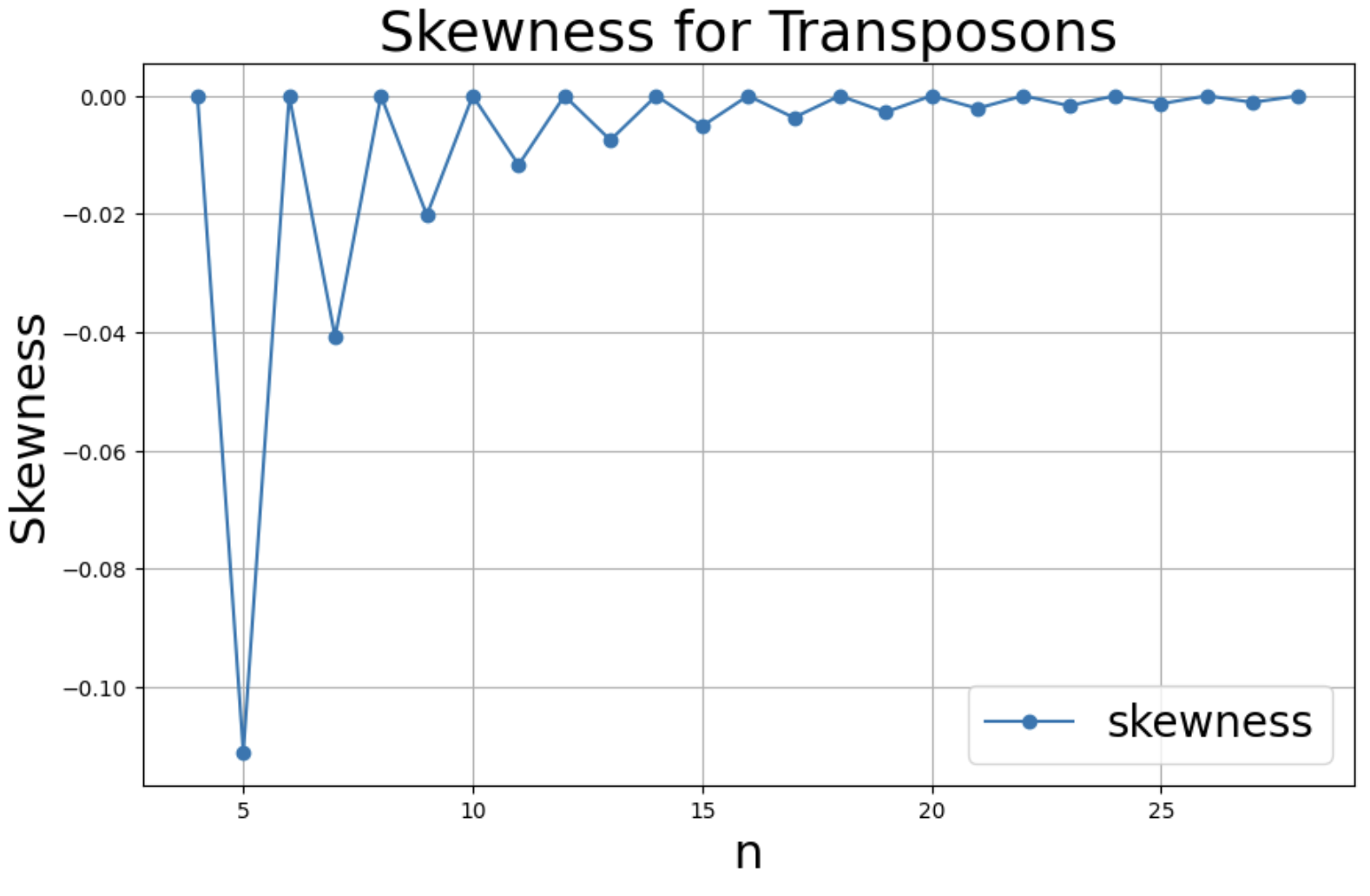}
    \caption{Skew  transposone binary coset}
    \label{fig:transposons-coset-skew}
  \end{minipage}
  \hfill
  \hfill
  \begin{minipage}{0.45\textwidth}
    \centering
    \includegraphics[width=1.1\linewidth]{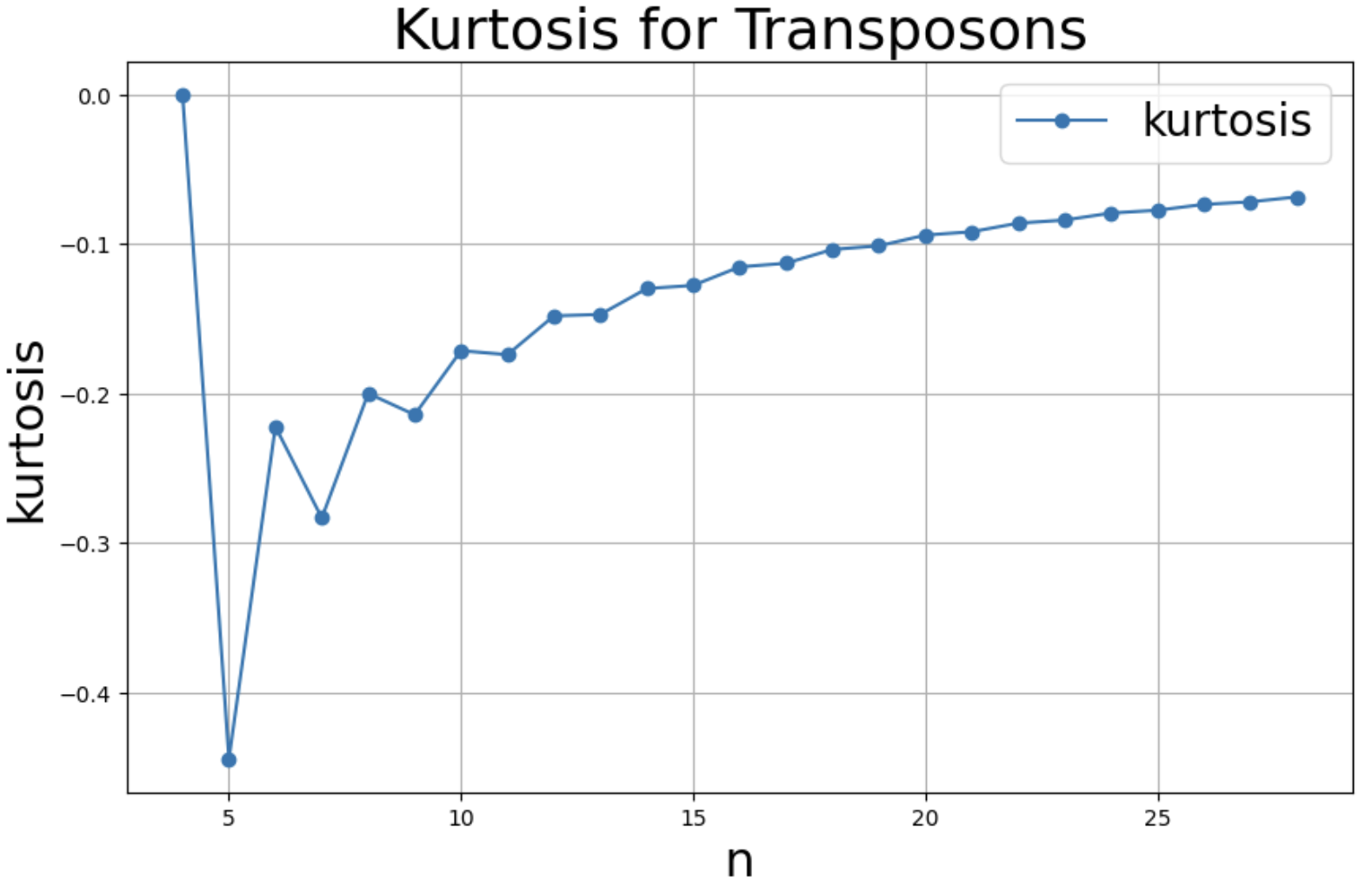}
    \caption{Kurtosis   transposone binary coset}
    \label{fig:kurtosis_transposons-coset}
  \end{minipage}
\end{figure}

\begin{figure}[H]
  \centering
  \begin{minipage}{0.45\textwidth}
    \centering
    \includegraphics[width=1.1\linewidth]{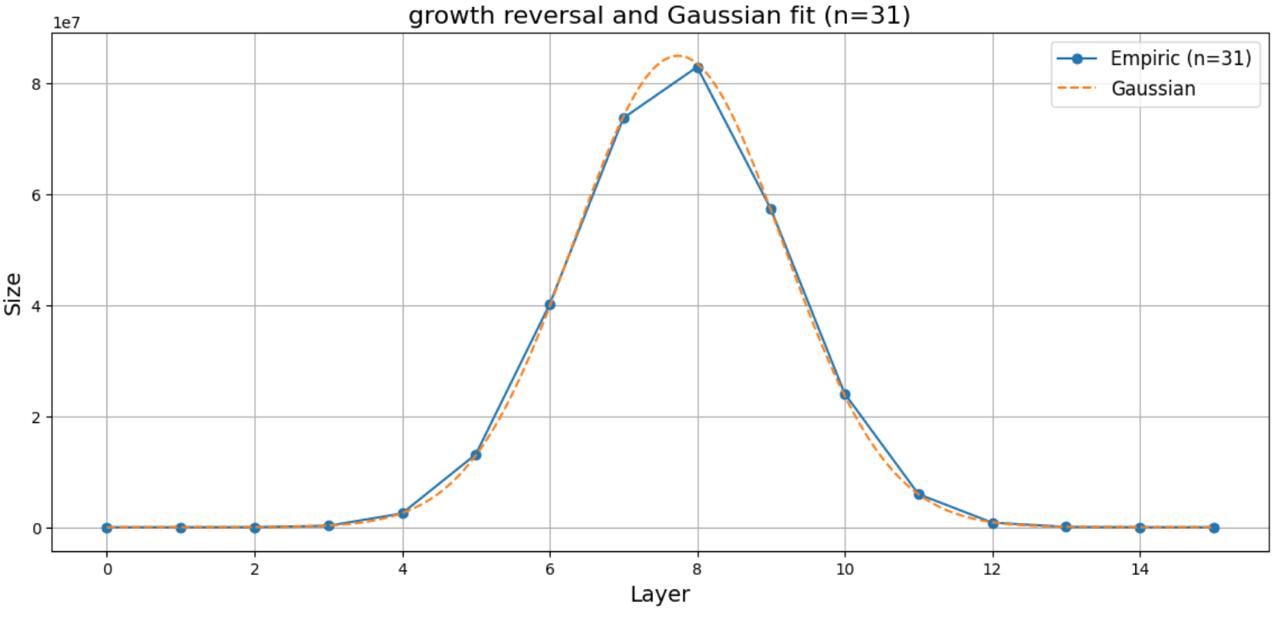}
    \caption{Gaussian fit for growth  transposone/reversals binary coset}
    \label{fig:Gauss-transposons-coset}
  \end{minipage}
  \hfill
  \hfill
  \begin{minipage}{0.45\textwidth}
    \centering
    \includegraphics[width=1.1\linewidth]{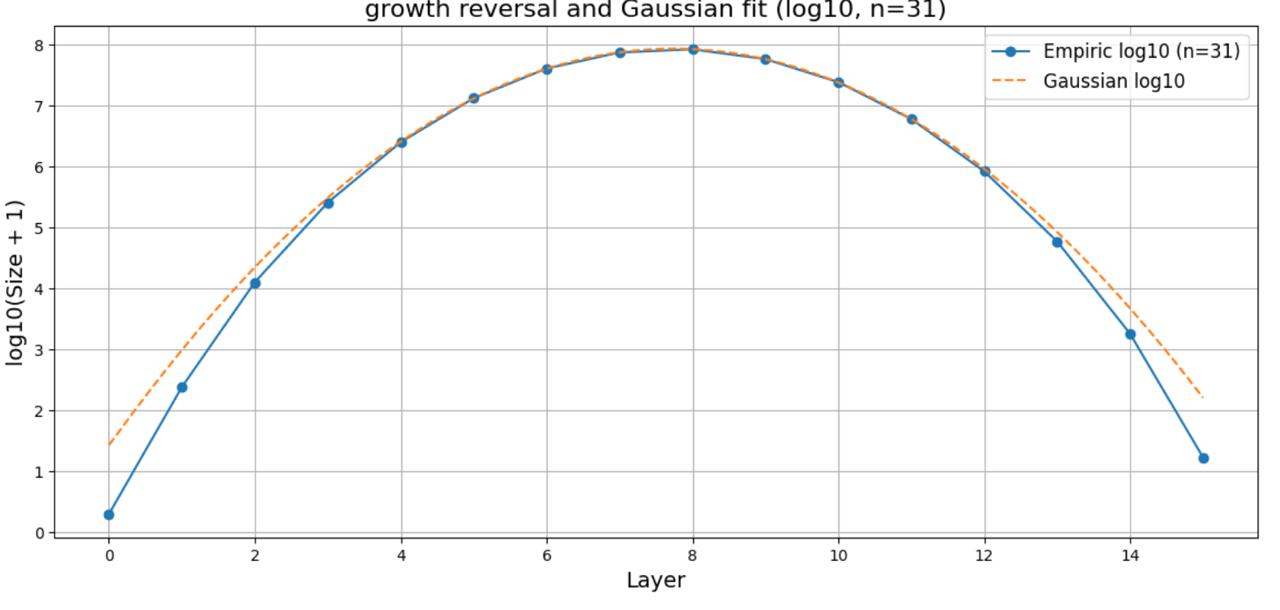}
    \caption{Log10 scale Gaussian fit for growth  transposone/reversals binary coset}
    \label{fig:logGauss_transposons-coset}
  \end{minipage}
\end{figure}

Notebooks \href{https://www.kaggle.com/code/eveelina/fork-of-cayleypy-transposons-analys-diameter-on-n}{Transposons with longest states},
\href{https://www.kaggle.com/code/eveelina/transposons-statistical-parameters}{Transposons main stat},
\href{https://www.kaggle.com/code/eveelina/reversal-statistical-parameters}{Reversals},
\href{https://www.kaggle.com/code/eveelina/cayleypy-reversal-norm-d}{Gaussian}.

\clearpage


\section{Case study -- consecutive 4-cycles \texorpdfstring{$(i,i+1,i+2,i+3)$}{i,i+1,i+2,i+3)}, use of AI} \label{sec:cons-cycles} %
Consider generators of $S_n$ given by consecutive 4-cycles $(i,i+1,i+2,i+4), 1\le i\le n-4 $ and their inverses.
Although these are rather obvious generators, to the best of our knowledge their diameter has not been discussed in the literature. We propose the following conjecture for this set of generators.

\begin{Conj}
    \begin{enumerate}
        \item The diameter is given by $n(n-1)/6 - 1$, for $n=3k, 3k+1, n\ge 6$, and $n(n-1)/6+2/3$, for $n=3k+2$, and $n\ge 8 $
        \item The longest elements (antipodes) include: 
        $[ n, n-1, ..., 3, 1, 0, 2]$, for  $n=3k, 3k+1, n\ge 6$ and 
        $[n, n-1, ..., 3, 2, 0, 1]$, for $n=3k+2$. 
        \item The mean diameter is approximately equal to $0.036n^2$, mode $\approx 0.045n^2$
    \end{enumerate}
\end{Conj}

A direct BFS (breadth first search) computation of growth in that case is feasible for $n< 14$ (so we verified the conjecture in that range).
Although diameters can be guessed already from these values of $n$ by fitting a quadratic quasi-polynomial, this might not be enough data to be confident, while direct computations for much larger $n$ are unfeasible. 
Here the AI-based CayleyPy's path-finding comes into play -- since we were able to guess not only the diameters, but the pattern of the 
longest elements  for all $n$, we take these elements for $n \ge 14$ and run our AI pipeline to decompose them into a product of the generators.
And "bingo" -- we see their lengths are exactly as predicted by the quasi-polynomial formula. 
This does not give a proof, but it gives more support to our conjecture, similar to what has been done for "LRX" generators in \cite{CayleyPyRL}.

We will present more examples like this in the full version of the manuscript and in further publications.


\section*{Acknowledgments}
A.C. is deeply grateful to M. Douglas for his interest in this work, engaging discussions, and invitation to present preliminary results at the Harvard CMSA program on "Mathematics and Machine Learning" in Fall 2024; to A.Hayt, F.Charton invitation to
to Oberwolfach workshop "Machine Learning and AI for mathematics" September 2025, multiple stimulating discussions 
and encouragements; to M. Gromov, S. Nechaev, and V. Rubtsov for their invitation to give a talk (\href{https://youtu.be/RkmBwlSyhfA?si=KgqQtRFaqx5ykd_s}{video}) at "Representations, Probability, and Beyond: A Journey into Anatoly Vershik World" at IHES, as well as for stimulating discussions.

The authors express heartfelt thanks to Tomas Rokicki  for multiple discussions, interest in this work and lots of explanations on his 
remarkable site \href{https://alpha.twizzle.net/explore/?puzzle=pyraminx}{Alpha Twizzle} and \href{https://github.com/cubing/twsearch}{Twsearch} program, and to Darij Grinberg for interest in the present work, stimulating discussions and suggesting various Cayley graphs considered here. 
Z.K. thanks Google Cloud for resources: research partly supported with Cloud TPUs from Google's TPU Research Cloud 
\href{https://sites.research.google/trc/about/}{(TRC)}.

 A.C. is grateful to J. Mitchel for involving  into the Kaggle Santa 2023 challenge, from which this project originated and to M.Kontsevich, Y.Soibelman,  S.Gukov,  T. Smirnova-Nagnibeda,  D.Osin, V. Kleptsyn, G.Olshanskii, A.Ershler, J. Ellenberg, G. Williamson, A. Sutherland,  Y. Fregier, P.A. Melies, I. Vlassopoulos, F.Khafizov, A.Zinovyev, M. Alekseyev
 for the discussions, interest and comments, to his wife A.Chervova and daugther K.Chervova for support, understanding and help with computational experiments.  

We are deeply grateful to the many colleagues who have contributed to the CayleyPy project at various stages of its development, including: N.Bukhal, J.Naghiev, A.Lenin, E.Uryvanov,  A. Abramov, M.Urakov, A.Kuchin,  B.Bulatov,  F.Faizullin,  U.Kniaziuk, D.Naumov,  S.Botman, A.Kostin,
R.Vinogradov, I.Gaiur,  N.Narynbaev, K.Khoruzhii, A.Romanov, A.Korolkova, N. Rokotyan, S.Kovalev, A.Eliseev, F.Petrov, S.Fironov, A.Lukyanenko, A.Ogurtsov, G.Antiufeev, G.Verbii, D.Gorodkov,  A.Rozanov, V.Nelin, S.Ermilov,
K.Yakovlev, V.Shitov, E.Durymanov,  R.Magdiev, M.Krinitskiy, P.Snopov. 

We are deeply indebted to many administrators of telegram channels who helped to widespread a project announcement
and to gather us together, in particular to Professor, Kaggle Grandmaster, ex-Kaggle Top1 Alexander Dyakonov,
who taught and inspired generations of data scientists \href{https://t.me/smalldatascience}{(channel)},
Dr. Grigory Sapunov \href{https://t.me/gonzo_ML}{(channel)}, Dr. Boris Tseytlin
\href{https://t.me/smalldatascience}{(channel)}, Alexander Abramov \href{https://t.me/dealerAI}{(channel)},
Dr. Tatiana Shavrina \href{https://t.me/rybolos_channel}{(channel)}, Dr. Aleksandr Nikolich 
\href{https://t.me/lovedeathtransformers}{(channel)}, to Vadim Gernar \href{https://t.me/gernar228}{(channel)} 
and many others.

\end{document}